\documentclass[12pt]{article}
\usepackage{amsthm,amsfonts,amssymb,amscd}

\begin{document}


\begin{center}
\large \bf Birational geometry \\
of Fano double spaces of index two
\end{center}\vspace{0.5cm}

\centerline{A.V.Pukhlikov}\vspace{0.5cm}

\parshape=1
3cm 10cm \noindent {\small \quad\quad\quad \quad\quad\quad\quad
\quad\quad\quad {\bf }\newline We study birational geometry of
Fano varieties, realized as double covers $\sigma\colon V\to
{\mathbb P}^M$, $M\geq 5$, branched over generic hypersurfaces
$W=W_{2(M-1)}$ of degree $2(M-1)$. We prove that the only
structures of a rationally connected fiber space on $V$ are the
pencils-subsystems of the free linear system $|-\frac12 K_V|$.
The groups of birational and biregular self-maps of the variety
$V$ coincide: $\mathop{\rm Bir} V=\mathop{\rm Aut} V$.}
\vspace{1cm}

\section*{Introduction}

{\bf 0.1. Setting up of the problem and formulation of the main
result.} The integer $M\geq 4$ is fixed throughout the paper. The
symbol ${\mathbb P}$ stands for the projective space ${\mathbb
P}^M$ over an algebraically closed field of characteristic zero
(in the first place, we mean ${\mathbb C}$). Let
$W=W_{2(M-1)}\subset{\mathbb P}$ be a smooth hypersurface of
degree $2(M-1)$. There exists a uniquely determined double cover
$$
\sigma\colon V\to{\mathbb P},
$$
branched over $W$. It can be explicitly defined as the
hypersurface, given by the equation
$$
x^2_{M+1}=f(x_0,\dots,x_M)
$$
in the weighted projective space ${\mathbb
P}^{M+1}(1,\dots,1,M-1)$, where $f(x_*)$ is the equation of the
hypersurface $W$.

The variety $V$ is a Fano variety of index two:
$$
\mathop{\rm Pic}V={\mathbb Z}H,
$$
where $H$ is the ample generator, $K_V=-2H$, the class $H$ is the
pull back via $\sigma$ of a hyperplane in ${\mathbb P}$. On the
variety $V$ there are the following natural structures of a
rationally connected fiber space: let $\alpha_P\colon{\mathbb
P}\dashrightarrow{\mathbb P}^1$ be the linear projection from an
arbitrary linear subspace $P$  of codimension two, then the map
$$
\pi_P=\alpha_P\circ\sigma\colon V\dashrightarrow{\mathbb P}^1
$$
fibers $V$ into $(M-1)$-dimensional Fano varieties of index one.
Now let us formulate the main result.

Recall [1], that a (non-trivial) rationally connected fiber space
is a surjective morphism $\lambda\colon Y\to S$ of projective
varieties, where $\mathop{\rm dim}S\geq 1$ and the variety $S$
and a fiber of general position $\lambda^{-1}(s)$, $s\in S$, are
rationally connected (and the variety $Y$ itself is automatically
rationally connected by the theorem of Graber, Harris and Starr
[2]).\vspace{0.1cm}

{\bf Theorem 1.} {\it Assume that $M\geq 5$ and the branch
hypersurface $W\subset{\mathbb P}$ is sufficiently general. Let
$\chi\colon V\dashrightarrow Y$ be a birational map onto the total
space of the rationally connected fiber space $\lambda\colon Y\to
S$. Then $S={\mathbb P}^1$ and for some isomorphism
$\beta\colon{\mathbb P}^1\to S$ and a subspace $P\subset{\mathbb
P}$ of codimension two we get
$$
\lambda\circ\chi=\beta\circ\pi_P,
$$
that is, the following diagram commutes:}
$$
\begin{array}{ccccc}
& V & \stackrel{\chi}{\dashrightarrow} & W &\\
\pi_P & \downarrow & & \downarrow & \lambda\\
& {\mathbb P}^1 & \stackrel{\beta}{\to} & S. &\\
\end{array}
$$

{\bf Corollary 1.} {\it For a generic double space $V$ of
dimension $\mathop{\rm dim} V\geq 5$ the following claims hold.}

(i) {\it On the variety $V$ there are no structures of a
rationally connected fiber space with the base of dimension $\geq
2$. In particular, on $V$ there are no structures of a conic
bundle and del Pezzo fibration, and the variety $V$ itself is
non-rational.

{\rm (ii)} Assume that there is a birational map $\chi\colon
V\dashrightarrow Y$, where $Y$ is a Fano variety of index $r\geq
2$ with factorial terminal singularities, such that $\mathop{\rm
Pic}Y={\mathbb Z}H_Y$, where $K_Y=-rH_Y$, and moreover, the
linear system $|H_Y|$ is non-empty and base point free. Then
$r=2$ and the map $\chi$ is a biregular isomorphism.

{\rm (iii)} The groups of birational and biregular self-maps of
the variety $V$ coincide:}
$$
\mathop{\rm Bir}V=\mathop{\rm Aut}V={\mathbb Z}/2{\mathbb Z}.
$$

{\bf Proof of the corollary.} The claim (i) and the equality
$r=2$ in (ii) are obvious (any linear subsystem of projective
dimension $\leq r-1$ in the complete linear system $|H_Y|$
defines a structure of a rationally connected fiber space on
$Y$). Furthermore, the $\chi$-preimage of a generic divisor in
the system $|H_Y|$ is by Theorem 1 a divisor in the linear system
$|H|$, which completes the proof of the claim (ii). The part
(iii) follows from (ii) in an obvious way. Q.E.D.

Tha aim of the present paper is to prove Theorem 1. As usual, its
claim will be derived from another fact, a much more technical
and less visual Theorem 2 on the thresholds of canonical
adjunction of movable linear systems on the variety $V$. However,
first of all, let us discuss the position of Theorem 1 in the
context of known results on birational geometry of
higher-dimensional Fano varieties.\vspace{0.3cm}


{\bf 0.2. From birationally rigid varieties to birationally
non-rigid ones.} Recall the principal definitions of the theory
of birational rigidity. Let $X$ be a smooth projective rationally
connected variety. It satisfies the classical condition of
termination of adjunction of the canonical class: for any
effective divisor $D$ the linear system $|D+mK_X|$ is empty for
$m\gg 0$, since $K_X$ is negative on every family of rational
curves sweeping out $X$, whereas an effective divisor is
non-negative on any such family. In order to fix the moment of
termination precisely, let us consider the Picard group
$A^1X=\mathop{\rm Pic}X$, set $A^1_{\mathbb
R}X=A^1X\otimes{\mathbb R}$ and define the cones $A^1_+X\subset
A^1_{\mathbb R}X$ of {\it pseudo-effective} classes and $A^1_{\rm
mov}X\subset A^1_{\mathbb R}X$ of {\it movable} clases as the
closed cones (with respect to the standard real topology of
$A^1_{\mathbb R}X\cong{\mathbb R}^k$), generated by the classes
of effective divisors and movable divisors (that is, divisors in
the linear systems with no fixed components), respectively.

{\bf Definition 0.1.} The {\it threshold of canonical adjunction}
of a divisor $D$ on the variety $X$ is the number
$c(D,X)=\mathop{\rm sup}\{\varepsilon\in {\mathbb
Q}_+|D+\varepsilon K_X\in A^1_+X\}$. If $\Sigma$ is a non-empty
linear system on $X$, then we set $c(\Sigma,X)=c(D,X)$, where
$D\in\Sigma$ is an arbitrary divisor.

{\bf Example 0.1.} (i) Let $X$ be a smooth Fano variety, and
assume that $\mathop{\rm Pic}X={\mathbb Z}H_X$, where $H_X$ is
the ample generator and $K_X=-rH_X$, $r\geq 1$. For any effective
divisor $D$ we have $D\sim nH_X$ for some $n\geq 1$, so that
$$
c(D,X)=\frac{n}{r}.
$$

(ii) Let $\pi\colon X\to S$ be a rationally connected fiber space
with $\mathop{\rm dim}X>\mathop{\rm dim}S\geq 1$, $\Delta$ an
effective divisor on the base $S$. Obviously,
$c(\pi^*\Delta,X)=0$. If $\mathop{\rm Pic}X={\mathbb
Z}K_X\oplus\pi^*\mathop{\rm Pic}S$, that is, $X/S$ is a standard
Fano fiber space, and $D$ is an effective divisor on $X$, which
is not a pull back of a divisor on the base $S$, then $D\in
|-mK_X+\pi^*R|$ for some divisor $R$ on $S$, where $m\geq 1$.
Obviously, $c(D,X)\leq m$, and moreover, if the divisor $R$ is
effective, then $c(D,X)= m$. Indeed, $K_X$ is negative on the
fibers of the morphism $\pi$ (in particular, on the dense
families of rational curves sweeping out fibers of $\pi$),
whereas any divisor, pulled back from the base, is trivial on the
fibers.

{\bf Example 0.2.} Let $Y$ be a Fano variety of index $r\geq 2$,
described in part (ii) of Corollary 1. Consider a movable linear
system $\Sigma$, spanned by the divisors $D_1,\dots,D_r\in |H_Y|$
of general position. Obviously,
$$
Q=D_1\cap\dots\cap D_r\subset Y
$$
is an irreducible subvariety of codimension $r\geq 2$, whereas
$\mathop{\rm Bs}\Sigma =Q$. Let $\varphi\colon Y^+\to Y$ be the
blow up of $Q$ and $\Sigma^+$ the strict transform of the system
$\Sigma$ on $Y^+$. Obviously, the free system $\Sigma^+$ defines
a morphism
$$
\pi_Q\colon Y^+\to{\mathbb P}^{r-1},
$$
the fibers of which are Fano varieties and therefore are
rationally connected. We get the equalities
$$
c(\Sigma,Y)=\frac{1}{r}\quad\mbox{and}\quad c(\Sigma^+,Y^+)=0,
$$
because $\Sigma^+$ is pulled back from the base ${\mathbb
P}^{r-1}$.\vspace{0.1cm}

{\bf Definition 0.2.} For a movable linear system $\Sigma$ on a
variety $X$ define the {\it virtual threshold of canonical
adjunction} by the formula
$$
c_{\rm virt}(\Sigma)=\mathop{\rm inf}\limits_{X^{\sharp}\to X}
\{c(\Sigma^{\sharp},X^{\sharp})\},
$$
where the infimum is taken over all birational morphisms
$X^{\sharp}\to X$, $X^{\sharp}$ is a smooth projective model of
${\mathbb C}(X)$, $\Sigma^{\sharp}$ is the strict transform of the
system $\Sigma$ on $X^{\sharp}$.

The virtual threshold is obviously a birational invariant of the
pair $(X,\Sigma)$: if $\chi\colon X\dashrightarrow X^+$ is a
birational map, $\Sigma^+=\chi_*\Sigma$ is the strict transform of
the system $\Sigma$ with respect to $\chi^{-1}$, then we get
$c_{\rm virt}(\Sigma)=c_{\rm virt}(\Sigma^+)$.

The following obvious claim shows how the information on the
virtual threshold of canonical adjunction makes it possible to
describe birational property of algebraic varieties.

{\bf Proposition 0.1.} (i) {\it Assume that on the variety $X$
there are no movable linear systems with the zero virtual
threshold of canonical adjunction. Then on $X$ there are no
structures of a non-trivial fibration into varieties of negative
Kodaira dimension, that is, there is no rational dominant map
$\rho\colon X\dashrightarrow S$, $\mathop{\rm dim}S\geq 1$, the
generic fiber of which has negative Kodaira dimension.

{\rm (ii)} Let $\pi\colon X\to S$ be a rationally connected fiber
space. Assume that every movable linear system $\Sigma$ on $X$
with the zero virtual threshold of canonical adjunction, $c_{\rm
virt}(\Sigma)=0$, is the pull back of a system on the base:
$\Sigma=\pi^*\Lambda$, where $\Lambda$ is some movable linear
system on $S$. Then any birational map
\begin{equation}\label{o1}
\begin{array}{ccccc}
 & X & \stackrel{\chi}{\dashrightarrow} & X^{\sharp} & \\
\pi & \downarrow & & \downarrow & \pi^{\sharp}\\
& S & & S^{\sharp}, &
\end{array}
\end{equation}
where  $\pi^{\sharp}\colon V^{\sharp}\to S^{\sharp}$ is a
fibration into varieties of negative Kodaira dimension, is
fiber-wise, that is, there exists a rational dominant map
$\rho\colon S \dashrightarrow S^{\sharp}$, making the diagram
(\ref{o1}) commutative,} $\pi^{\sharp}\circ\chi=\rho\circ\pi$.

The only known way to compute the virtual thresholds of canonical
adjunction is by reduction to the ordinary thresholds.

{\bf Definition 0.3.} (i) The variety $X$ is said to be {\it
birationally superrigid}, if for any movable linear system
$\Sigma$ on $X$ the equality
$$
c_{\mathop{\rm virt}}(\Sigma)=c(\Sigma,X)
$$
holds.

(ii) The variety $X$ (respectively, the Fano fiber space $X/S$) is
said to be {\it birationally rigid}, if for any movable linear
system $\Sigma$ on $X$ there exists a birational self-map
$\chi\in\mathop{\rm Bir}X$ (respectively, a fiber-wise birational
self-map $\chi\in\mathop{\rm Bir}(X/S)$), providing the equality
$$
c_{\mathop{\rm virt}}(\Sigma)=c(\chi_*\Sigma,X).
$$

(iii) The variety $X$ is said to be {\it almost birationally
rigid}, if it has a model $\widetilde X$, which is a birationally
rigid variety, that is, the condition (ii) is satisfied for some
smooth projective variety $\widetilde X$, birational to
$X$.\vspace{0.1cm}

{\bf Example 0.3.} (i) A smooth three-dimensional quartic
$X=X_4\subset{\mathbb P}^4$ is birationally superrigid:  this
follows immediately from the arguments of [3].

(ii) A generic smooth hypersurface $X=X_d\subset{\mathbb P}^d$,
$d\geq 5$, is a birationally superrigid variety [4]. A generic
complete intersection $X_{d_1\cdot\dots\cdot d_k}\subset{\mathbb
P}^{M+k}$ of index one (that is, $d_1+\dots+d_k=M+k$) and
dimension $M\geq 4$ is birationally superrigid for $X\geq 2k+1$
[5]. For more examples, see [1].

(iii) Generic smooth complete intersections $X_{2\cdot
3}\subset{\mathbb P}^5$ of a cubic and a quadric hypersurfaces
are birationally rigid, but not superrigid [6,7]. A generic
complete intersection $X_{2\cdot 3}\subset{\mathbb P}^5$ with a
non-degenerate quadratic singularity $o\in X_{2\cdot 3}$ is
almost birationally rigid. For a birationally rigid model
$\widetilde {X_{2\cdot 3}}$ one can take the blow up of the
singular point [8].\vspace{0.1cm}

{\bf Conjecture 0.1.} {\it A smooth Fano complete intersection of
index one and dimension $\geq 4$ in a weighted projective space is
birationally rigid, of dimension $\geq 5$ birationally
superrigid.}

It follows from Example 0.2 that the double space $V$, the main
object of study in the present paper, is neither superrigid, nor
rigid, nor almost rigid.\vspace{0.3cm}


{\bf 0.3. The start and scheme of the proof of Theorem 1.}
Similarly to birationally rigid varieties, Theorem 1 is based on
some claim on the virtual threshold of canonical adjunction of a
movable linear system on $V$. For an arbitrary linear subspace
$P\subset{\mathbb P}$ of codimension two let $V_P$ be the blow up
of the subvariety $\sigma^{-1}(P)\subset V$ (it is irreducible by
the conditions of general position, see Sec. 0.4). For a movable
linear system $\Sigma$ on $V$ the symbol $\Sigma_P$ stands for
its strict transform on $V_P$.

Theorem 1 follows from a more technical fact.

{\bf Theorem 2.} {\it Assume that $M\geq 5$ and for a movable
linear system $\Sigma$ the inequality
\begin{equation}\label{o2}
c_{\rm virt}(\Sigma)<c(\Sigma,V)
\end{equation}
holds. Then there exists a uniquely determined linear subspace
$P\subset{\mathbb P}$ of codimension two, satisfying the
inequality
$$
\mathop{\rm mult}\nolimits_{\sigma^{-1}(P)}\Sigma>c(\Sigma,V),
$$
whereas for the strict transform $\Sigma_P$ the equality
$$
c_{\rm virt}(\Sigma)=c_{\rm virt}(\Sigma_P)=c(\Sigma_P,V_P)
$$
holds.}

Theorem 1 is derived from Theorem 2 in a few lines, see \S 1.
Almost all paper is devoted to proving Theorem 2. Let us fix a
movable linear system $\Sigma$, satisfying the inequality
(\ref{o2}). Taking, if necessary, a symmetric power of $\Sigma$,
we may assume that
$$
\Sigma\subset|2nH|=|-nK_V|,
$$
where $n\geq 1$ is a positive integer. The system $\Sigma$ (and
the integer $n$) are fixed throughout the paper, with the
exception of a few technical sections (in the first place, \S\S
4-5), where the notations are independent; we usually point this
out but it is always clear from the context. Obviously,
$$
c(\Sigma,V)=n.
$$

{\bf Proposition 0.2 (the Noether-Fano inequality).} {\it There
exists a birational morphism $\varphi\colon\widetilde{V}\to V$
and an irreducible exceptional divisor $E\subset\widetilde{V}$,
satisfying the estimate}
\begin{equation}\label{o3}
\mathop{\rm ord}\nolimits_E\varphi^*\Sigma>na(E,V).
\end{equation}

{\bf Proof} is well known (see, for instance, [1]).\vspace{0.1cm}

The divisor $E$ (or the corresponding discrete valuation of the
field of rational functions of the variety $V$) is called a {\it
maximal singularity} of the linear system $\Sigma$. If $\varphi$
is the blow up of an irreducible subvariety $B\subset V$ (and in
that case $E=\varphi^{-1}(B)$), then the latter is called a {\it
maximal subvariety} of the system $\Sigma$. In that case
(\ref{o3}) is equivalent to the inequality
$$
\mathop{\rm mult}\nolimits_B\Sigma>n(\mathop{\rm codim}B-1).
$$
In any case the subvariety $B=\varphi(E)$ is called the {\it
centre} of the maximal singularity $E$, see [1] for the language
of maximal singularities.

An equivalent formulation of Proposition 0.2: the pair
\begin{equation}\label{o4}
\left(V,\frac{1}{n}\Sigma\right)
\end{equation}
is not canonical and the prime divisor $E\subset\widetilde{V}$ is
a non canonical singularity of this pair. If, instead of
(\ref{o3}), the stronger inequality
\begin{equation}\label{o5}
\mathop{\rm ord}\nolimits_E\varphi^*\Sigma>n(a(E)+1),
\end{equation}
holds, then $E$ is a {\it log maximal singularity}, the pair
(\ref{o4}) is not log canonical and the estimate (\ref{o5}) is
the log Noether-Fano inequality. The bigger part of the paper is
devoted to proving the following fact.

{\bf Proposition 0.3.} {\it There exists a unique linear subspace
$P\subset{\mathbb P}$ of codimension two, such that the subvariety
$\sigma^{-1}(P)$ is a maximal subvariety of the system $\Sigma$.}

Proposition 0.3 is proved in \S\S 2-6 in the ``negative''
version: assuming that the system $\Sigma$ has no maximal
subvarieties of the form $\sigma^{-1}(P)$, where
$P\subset{\mathbb P}$ is a linear subspace of codimension two, we
exclude one by one all possibilities for a maximal singularity of
the system $\Sigma$, thus coming to a contradiction with
Proposition 0.2. The arguments of \S 2 exclude also the
possibility that the system $\Sigma$ has two maximal subvarieties,
$\sigma^{-1}(P_1)$ and $\sigma^{-1}(P_2)$, where $P_1\neq P_2$
are distinct linear subspaces of codimension two.

Now we blow up the maximal subvariety $\sigma^{-1}(P)$ and on the
new (generally speaking, singular) variety complete the proof of
Theorem 2, which directly implies Theorem 1. This part of our
work, although it is the concluding one in the sense of the proof
as a whole, is based on Proposition 0.3 only and is independent
of the contents of \S\S 2-6. So it is done in \S 1 (in the
assumption that Proposition 0.3 holds).

\S 7 is devoted to proving the conditions of general position,
which the double space $V$ is supposed to satisfy. There is no
doubt that these conditions are unnecessary, that is, that
Theorems 1 and 2 are true for any smooth double space of index
two. However, the conditions of general position are essentially
used in the proof. Some of those conditions could have been at
least relaxed, however, this would have made our proof, hard and
long as it is, even more complicated.

Note that in this paper we do a considerable part of work for the
double spaces of dimension $M=4$ (we exclude a majority of types
of maximal singularities). For $M=5$, in order to avoid the paper
getting too long, we omit the proof for one of the cases when we
exclude infinitely near maximal singularities (\S 6). For $M\geq
6$ we consider all possible cases.\vspace{0.3cm}


{\bf 0.4. Formulation of conditions of general position.} As we
mentioned above, the main result of the present paper is obtained
under the assumption that the double space $V$ is sufficiently
general (that is, the branch divisor $W\subset {\mathbb P}$ is a
sufficiently general hypersurface of degree $2(M-1)$). We will
use several conditions of general position, the three principal
ones of which are formulated below and proved in \S 7 (we show
that a general hypersurface $W$ satisfies these conditions
indeed). Some other, less significant, conditions are given where
they are used.

The first main condition deals with lines on $V$. As usual, a
curve $C\subset V$ is called a {\it line}, if the equality
$(C\cdot H)=1$ holds. In particular, a line is a smooth
irreducible rational curve. We have\vspace{0.1cm}

{\bf Proposition 0.4.} {\it On a generic variety $V$ there are
finitely many lines through any point.}\vspace{0.1cm}

The second and third condition deal with linear subspaces (planes)
in $\mathbb P$ of codimension two. Consider an arbitrary plane
$P\subset{\mathbb P}={\mathbb P}^M$ of codimension two. The
intersection $P\cap W$, generally speaking, is singular:
$$
p\in\mathop{\rm Sing}P\cap W
$$
if and only if
$$
P\subset T_pW.
$$
It is well known that (without the assumption that the
hypersurface $W$ is generic) the set $\mathop{\rm Sing}P\cap W$
is at most one-dimensional (see, for instance, [7]).

The assumption that $W$ is generic makes it possible to improve
this claim.\vspace{0.1cm}

{\bf Proposition 0.5.} {\it For a generic hypersurface $W$ and an
arbitrary plane $P\subset{\mathbb P}$ of codimension two the set
$\mathop{\rm Sing}P\cap W$ is finite (or empty). In particular,
the closed set $R=\sigma^{-1}(P)$ is irreducible, that is, it is
a subvariety, and the set of its singular points is at most
finite.}\vspace{0.1cm}

The third condition characterizes the singularities of the variety
$\sigma^{-1}(P)$ and the singularities of its blow up on $V$. For
a quadratic singular point, that is, a hypersurface singularity
with a local equation
$$
0=w_2(u_1,\dots,u_N)+w_3(u_*)+\dots,
$$
where $w_i(u_*)$ is a homogeneous polynomial of degree $i$, then
we say that this point is of rank $\mathop{\rm rk} w_2$. When such
a singularity is blown up, the exceptional divisor is the quadric
$\{w_2=0\}\subset{\mathbb P}^{N-1}$ of rank $\mathop{\rm rk} w_2$.

Let $V_P$ be the blow up of the (irreducible by Proposition 0.5)
subvariety $\sigma^{-1}(P)$ on $V$.\vspace{0.1cm}

{\bf Proposition 0.6.} {\it For a generic hypersurface $W$ and an
arbitrary plane $P\subset{\mathbb P}$ of codimension two:

{\rm (i)} for $M\geq 6$ every singular point of the variety $V_P$
is an isolated quadratic point of rank $\geq 4$,

{\rm (ii)} for $M\geq 4$ every singular point of the variety $V_P$
is an isolated quadratic point of rank $\geq 3$,

{\rm (iii)} for $M\geq 6$ every singular point of the variety
$\sigma^{-1}(P)$ is an isolated quadratic point of rank $\geq
2$.}\vspace{0.1cm}

The properties, formulated in Propositions 0.4-0.6, will be
assumed to take place (sometimes we remind about this in the
course of our arguments).\vspace{0.3cm}


{\bf 0.5. Historical remarks.} Up to this day, only very few
papers were describing birational geometry of Fano varieties of
index $r\geq 2$. (We should explain: here we mean dealing with the
problems that give information about the birational type of a
variety, such as the rationality problem, computation of the
group of birational self-maps, description of the set of
structures of rationally connected fiber spaces etc. There are a
lot of papers where {\it particular} birational maps are
constructed and studied, for instance, birational transformations
of the projective space, but this is a completely different area.)
Fano himself pioneered the study of varieties of index $\geq 2$,
he tried to describe the group of birational self-maps of the
three-dimensional cubic [9]. This attempt, as it is clear now,
could not be successful: the problem was too complicated for the
methods of his time. After V.A.Iskovskikh and Yu.I.Manin in 1971
proved birational superrigidity (in the modern terminology) of
three-dimensional quartics, it was natural to try to apply the
new technique of the method of maximal singularities to varieties
of higher index, and such an attempt was immediately made: in
[10] certain auxiliary claims are formulated for varieties of
arbitrary index $r\geq 1$ and some work is done on description of
birational geometry of the Veronese double cone of dimension
three (it is a Fano variety of index two). The paper [11] aimed
at completing that work (in particular, at solving the
rationality problem for that class of varieties). Unfortunately,
the above mentioned paper [11], as it became clear later, was
faulty, see [12], where the proof was completed 20 years later.
However, the fact that there was a mistake in [11], was already
clear enough in the mid-nineties: in order to study the Veronese
double cone successfully, one should be able to exclude maximal
singularities of movable linear systems on the pencils of del
Pezzo surfaces (because there are such pencils on the double
cone), whereas the technique of their exclusion was developed in
[13]. It is impossible to solve this problem by the methods that
were used in the eighties (the test class technique).

It is worth mentioning that the talk [14] at the ICM in Warsaw
announced S.I.Khashin's description of the group of birational
self-maps of the double space of index two and dimension three
(corresponding to the value $M=3$ in the notations of the present
paper), however, this announcement was not confirmed later and no
proof was produced. (Note that the three-dimensional space of
index two is a much more difficult object of study than the
Veronese double cone, so that the announcement, given in [14],
looks somewhat naive.) For the modern techniques this variety
seems already to be within reach, however, the problem is still
very hard. For the up to date description of this problem, see
[15,16].

Thus the Veronese double cone of dimension three up to this day
was the only Fano variety of index two, birational geometry of
which was completely studied. A series of remarkable results were
obtained by means of other methods: non-rationality of the
three-dimensional cubic was proved by Clemens and Griffiths in
[17] (see also [18,19]), non-rationality of ``very general'' Fano
hypersurfaces of arbitrary dimension and of index two and higher
was proved by Koll\' ar [20], non-rationality of double spaces
${\mathbb P}^3$ of index two follows from the fact that they admit
no structures of a conic bundle, as A.S.Tikhomirov showed in
[21-23]. These are just three examples; we do not give a complete
list of those results, since the methods used in the above
mentioned papers are very far from the method of maximal
singularities that makes the basis of the present paper. It should
be mentioned, though, that both the ``transcendent method'' (or
the method of intermediate Jacobian, developed by Clemens and
Griffiths) and Koll\' ar's approach make it possible to obtain
much less information about birational geometry of a given variety
than the method of maximal singularities that gives its almost
exhaustive description. In particular, only the method of maximal
singularities describes all structures of a rationally connected
fiber space (in the case when the work is completed).

The Veronese double cone of dimension three is, in a sense, an
exceptional variety by its numerical characteristics. The double
spaces of arbitrary dimension, considered in this paper, are
already quite typical. Theorems 1 and 2 show that the behaviour
which is natural to expect from higher-dimensional Fano varieties
of index two, really takes place.\vspace{0.3cm}


\section{The structures of rationally connected fiber spaces}

In this section we prove the main results of the paper, Theorems
1 and 2, assuming that Proposition 0.3 on the maximal subvariety
of codimension two holds.\vspace{0.3cm}

{\bf 1.1. Fano fiber space over ${\mathbb P}^1$.} According to
Proposition 0.3, there exist a (unique) linear subspace $P\subset
{\mathbb P}$ of codimension two, satisfying the estimate
$$
\mathop{\rm mult}\nolimits_R \Sigma >n,
$$
where $R=\sigma^{-1}(P)$ is an irreducible variety with at most
zero-dimensional singularities (Proposition 0.5). Let
$\varphi\colon V^+\to V$ be the blow up of the (possibly
singular) subvariety $R=\sigma^{-1}(P)$, $E=\varphi^{-1}(R)$ the
exceptional divisor.\vspace{0.1cm}

{\bf Lemma 1.1.} (i) {\it The variety $V^+$ is factorial and has
at most finitely many isolated double points (not necessarily
non-degenerate).

{\rm (ii)} The linear projection $\pi_{\mathbb P}\colon{\mathbb
P}\dashrightarrow{\mathbb P}^1$ from the plane $P$ generates the
regular projection
$$
\pi=\pi_{\mathbb P}\circ\sigma\circ\varphi\colon V^+\to{\mathbb
P}^1,
$$
the general fiber of which $F_t=\pi^{-1}(t)$, $t\in{\mathbb P}^1$
is a non-singular Fano variety of index one, and finitely many
fibers have isolated double points.

{\rm (iii)} The following equalities hold:
$$
\mathop{\rm Pic}V^+={\mathbb Z}H\oplus{\mathbb Z}E={\mathbb
Z}K^+\oplus{\mathbb Z}F,
$$
where $H=\varphi^*H$ for simplicity of notations, $K^+=K_{V^+}$
is the canonical class of the variety $V^+$, $F$ is the class of
a fiber of the projection $\pi$, whereas}
$$
K^+=-2H+E,\,\,F=H-E.
$$
\vspace{0.1cm}

{\bf Proof.} These claims follow directly from the definition of
the blow up $\varphi$, Proposition 0.5 and the well known fact
that an isolated hypersurface singularity of a variety of
dimension $\geq 4$ is factorial (see [24]).

Let $\Sigma^+$ be the strict transform of the system $\Sigma$ on
the blow up $V^+$ of the subvariety $R$.\vspace{0.1cm}

{\bf Proposition 1.1.} {\it The following equality holds:}
$$
c_{\rm virt}(\Sigma^+)=c(\Sigma^+,V^+).
$$
\vspace{0.1cm}

{\bf Proof of Theorem 2.} This theorem is just the union of
Proposition 0.3 and Proposition 1.1. Q.E.D.\vspace{0.1cm}

{\bf Corollary 1.1.} {\it Assume that $c_{\rm virt}(\Sigma^+)=0$.
Then the system $\Sigma^+$ is composed from the pencil $|H-R|$,
that is,}
$$
\Sigma^+\subset |2nF|.
$$
\vspace{0.1cm}

{\bf Proof of the corollary.} Assume the converse:
$$
\Sigma^+\subset |-mK^++lF|,
$$
where $m\geq 1$. By the part (iii) of Lemma 1.1,
$$
m=2n-\nu,\,\,l=2\nu-2n\geq 2,
$$
so that for the threshold of canonical adjunction we get
$$
c(\Sigma^+,V^+)=m.
$$
Since $c_{\rm virt}(\Sigma^+)=0$, by Proposition 1.1 we get
$m=0$, as we claimed. Q.E.D. for the corollary.\vspace{0.1cm}

{\bf Proof of Theorem 1.} For the linear system $\Sigma$ we take
the strict transform with respect to $\chi$ of any linear system
of the form $\lambda^*\Lambda$, where $\Lambda$ is a movable
system on the base $S$. Applying Theorem 2 (or Corollary 1.1), we
complete the proof.\vspace{0.3cm}


{\bf 1.2. Movable linear systems on the variety $V$.} We start
our proof of Proposition 1.1 with the well known step: assume
that the inequality
$$
c_{\rm virt}(\Sigma^+)< c(\Sigma^+,V^+)=m
$$
holds. Then the pair
\begin{equation}\label{a1}
(V^+,\frac{1}{m}\Sigma^+)
\end{equation}
is not canonical, so that the linear system $\Sigma^+$ has a
maximal singularity, that is, for some birational morphism
$\psi\colon \widetilde{V}\to V^+$ and irreducible exceptional
divisor $E^+\subset \widetilde{V}$ the Noether-Fano inequality
holds:
$$
\nu_E(\Sigma^+)>ma(E^+,V^+).
$$
\vspace{0.1cm}

{\bf Lemma 1.2.} {\it The centre of maximal singularity $E^+$ is
contained in some fiber $F_t=\pi^{-1}(t)$, that is,}
$$
B=\pi\circ\psi(E^+)=t\in{\mathbb P}^1.
$$
\vspace{0.1cm}

{\bf Proof.} Assume the converse: $\pi\circ\psi(E^+)={\mathbb
P}^1$. Restricting the linear system $\Sigma^+$ onto the fiber of
general position $F=F_s$, we get that the pair
$$
(F,\frac{1}{m}\Sigma_F)
$$
is not canonical, where $\Sigma_F\subset |-mK_F|$. However, $F$ is
a smooth double space of index one and it is well known [25], that
this is impossible. Q.E.D. for the lemma.\vspace{0.1cm}

For simplicity of notations, let $F=F_t$ be the fiber, containing
the centre of singularity $E^+$.\vspace{0.1cm}

{\bf Proposition 1.2.} {\it The centre $B$ is a singular point of
the fiber $F$.}\vspace{0.1cm}

{\bf Proof.} Since the anticanonical degree of the divisor
$D_F\in\Sigma_F$ is $2m$, and by genericity of the branch divisor
the anticanonical degree of any subvariety of codimension one on
$F$ is at least 2, we get the inequality $\mathop{\rm
codim}_FB\geq 2$, so that
$$
\mathop{\rm codim}\nolimits_{V^+}B\geq 3.
$$
In the notations of Sec. 1.1 let
$\Pi=\sigma\circ\varphi(F)\subset{\mathbb P}$ be the hyperplane,
corresponding to the fiber $F$. It is easy to see that
$$
\sigma_F=\sigma\circ\varphi\colon F\to\Pi={\mathbb P}^{M-1}
$$
is the double cover, branched over $W_{\Pi}=W\cap\Pi$: the blow
up $\varphi$ does not affect the divisors-elements of the pencil
$|H-R|$. Now we need to consider two cases:\vspace{0.1cm}

1) $\sigma_F(B)\not\subset W_{\Pi}$,

2) $\sigma_F(B)\subset W_{\Pi}$, but the generic point of the
subvariety $B$ is a non-singular point of the fiber
$F$.\vspace{0.1cm}

The case 1) is excluded by the arguments of [26]. Let $o\in B$ be
a point of general position,
$$
\lambda\colon F^{\sharp}\to F
$$
its blow up, $E^{\sharp}=\lambda^{-1}(o)\subset F^{\sharp}$ the
exceptional divisor, $E^{\sharp}\cong{\mathbb P}^{M-2}$. By
inversion of adjunction for a general divisor $D\in\Sigma^+$ we
get: the pair
\begin{equation}\label{a2}
(F,\frac{1}{m}D_F)
\end{equation}
is not log canonical at $B$, so that by [26, Proposition 3] there
is a hyperplane $\Lambda\subset E^{\sharp}$, satisfying the
inequality
$$
\mathop{\rm mult}\nolimits_oD_F+\mathop{\rm
mult}\nolimits_{\Lambda}D^{\sharp}_F>2m,
$$
where $D^{\sharp}_F$ is the strict transform of the divisor $D_F$
on $F^{\sharp}$. Now the arguments of [26, Sec. 2.2] give a
contradiction.\vspace{0.1cm}

Consider the case 2). If $\mathop{\rm dim}B\geq 1$, then for a
point $o\in B$ of general position the intersection of divisors
$$
T_pW_{\Pi}\quad \mbox{and} \quad\sigma_F(D_F),
$$
where $p=\sigma_F(o)$, is of codimension two (by the condition of
general position, for any hyperplane $\Lambda\subset\Pi$ we get
$\mathop{\rm dim}\mathop{\rm Sing}\Lambda\cap W=0$, so that the
tangent hyperplanes $T_pW_{\Pi}$, $p\in B$, form a $\mathop{\rm
dim}B$-dimensional family). In particular, the scheme-theoretic
intersection
$$
(\sigma^{-1}_F(T_pW_{\Pi})\circ D_F)
$$
is an effective cycle of codimension two on $F$, of $H$-degree
$2m$ and of multiplicity at least
$$
2\mathop{\rm mult}\nolimits_oD_F>2m
$$
at the point $o$, which is impossible.

Thus it remains to consider the case when $B=o$ is a smooth point
on the ramification divisor of the morphism $\sigma_F$. Since the
condition of non log canonicity of the pair (\ref{a2}) is linear
in the divisor $D_F\in|-mK_F|$, one may assume that $D_F$ is a
prime divisor. Set $\Lambda=T_pW_{\Pi}$. If
$$
D_F\neq\sigma^{-1}_F(\Lambda),
$$
then we argue as above in the case $\mathop{\rm dim}B\geq 1$. Let
us show that the equality $D_F=\sigma^{-1}_F(\Lambda)$ is
impossible. It can be done by inspection of possible singularities
of the intersection $W_{\Pi}\cap\Lambda$ for a hypersurface $W$ of
general position. We will give a simpler argument suggested by the
anonymous referee of the paper [27], see also [28]. Namely, if the
pair (\ref{a2}) is not log canonical for
$D_F=\sigma^{-1}_F(\Lambda)$, then by [29] (or [30]), the pair
$$
(\Pi,\Lambda+\frac12W_{\Pi})
$$
is not log canonical, either, which, in its turn, implies that
the pair
$$
(\Lambda,\frac12W_{\Lambda})
$$
is not log canonical, $W_{\Lambda}=(W\circ\Lambda)=W\cap\Lambda$.
However, as we pointed out above, the restriction $W_{\Lambda}$
has at most isolated double points as singularities. This
contradiction proves Proposition 1.2.\vspace{0.1cm}

Let $B=o$ be the centre of the maximal singularity
$E^+$.\vspace{0.1cm}

{\bf Proposition 1.3.}  {\it The point $o$ is a singularity of
the variety $V^+$}.\vspace{0.1cm}

{\bf Proof.} Assume the converse: the point $o\in V^+$ is
non-singular. Since the pair (\ref{a1}) is not canonical, we get
the inequality
$$
\mathop{\rm mult}\nolimits_o\Sigma^+>m,
$$
whence by Proposition 1.2 it follows that
$$
\mathop{\rm mult}\nolimits_oD_F>2m
$$
(since $o\in F$ is a singular point of the fiber). As we pointed
out above, this is impossible, which proves the
proposition.\vspace{0.3cm}


{\bf 1.3. The centre of the maximal singularity is a singular
point of the variety $V^+$.} We have shown above that the centre
of the maximal singularity $E^+$ is a singular point $o\in V^+$,
which we will assume from now on. Let
$$
\lambda\colon V^{\sharp}\to V^+
$$
be the blow up of the point $o$,
$E^{\sharp}=\lambda^{-1}(o)\subset V^{\sharp}$ the exceptional
divisor, which can be seen as a quadratic hypersurface in
${\mathbb P}^M$.

Recall (Proposition 0.6), that for $M\geq 6$ we may assume that
for a generic hypersurface $W\subset{\mathbb P}$, arbitrary plane
$P\subset{\mathbb P}$ of codimension two and any singularity
$o\in V^+$ the quadric $E^{\sharp}$ is of rank at least 4.

Define the integer $\beta\in{\mathbb Z}_+$ by the formula
$$
D^{\sharp}\sim\lambda^*D-\beta E^{\sharp},
$$
where $D\in\Sigma^+$ is a generic divisor, $D^{\sharp}$ its
strict transform on $V^{\sharp}$. By Proposition 1.4, which we
prove below, Proposition 0.6 implies the inequality
$$
\beta>m.
$$
Furthermore, the divisor
$$
\lambda^*_FD_F-\beta E^{\sharp}_F
$$
on the strict transform $F^{\sharp}\subset V^{\sharp}$ is
effective (the symbols $\lambda_F$ and $E^{\sharp}_F$ stand for
the blow up of the point $o\in F$ and for the exceptional divisor
$\lambda^{-1}_F(o)$, respectively). This implies the inequality
$$
\mathop{\rm mult}\nolimits_oD_F\geq 2\beta>2m,
$$
which is impossible. Proof of Proposition 1.1 for $M\geq 6$ is
complete.\vspace{0.3cm}


{\bf 1.4. Maximal singularities over quadratic points.} Consider
the following local situation. Let $o\in X$ be a germ of a
quadratic singularity, $\mathop{\rm dim}X\geq 3$. Let us blow up
the point $o$:
$$
\lambda\colon X^+\to X,
$$
and denote by the symbol $E$ the exceptional divisor
$\lambda^{-1}(o)$, which we consider as a quadric hypersurface
$$
E\subset{\mathbb P}^{\mathop{\rm dim} X}
$$

Let, furthermore, $D$ be an effective ${\mathbb Q}$-Cartier
divisor on the variety $X$, $D^+$ its strict transform on $X^+$.
Assuming the exceptional quadric $E$ to be irreducible, define
the number $\beta\in{\mathbb Q}_+$ by the relation
$$
D^+\sim\lambda^*D-\beta E.
$$
\vspace{0.1cm}

{\bf Proposition 1.4.} {\it Assume that the rank of the quadric
hypersurface $E$ is at least 4 and the pair
$$
(X,D)
$$
has the point $o$ as an isolated centre of a non canonical
singularity, that is, it is non canonical, but canonical outside
the point $o$. Then the following inequality holds:}
$$
\beta>1.
$$
\vspace{0.1cm}

{\bf Proof.} If $\mathop{\rm dim}X=3$, then by assumption the
point $o\in X$ is a non-degenerate quadratic singularity, and this
fact is well known [31]. (If $\beta\leq 1$, then the pair
$(X^+,D^+)$ is non canonical, so that by inversion of adjunction
the pair $(E,D^+_E)$ is not log canonical, but $E\cong{\mathbb
P}^1\times{\mathbb P}^1$ and $D^+_E$ is an effective curve of
bidegree $(\beta,\beta)$, which is impossible [32].) If
$\mathop{\rm dim}X\geq 4$, then, restricting $D$ onto a generic
hyperplane section $Y\ni o$ of the variety $X$ with respect to
some embedding $X\hookrightarrow{\mathbb P}^N$, and repeating this
procedure $\mathop{\rm dim}X-3$ times, we reduce the problem (by
inversion of adjunction) to the already considered case
$\mathop{\rm dim}X=3$. Proof of the proposition is
complete.\vspace{0.3cm}


{\bf 1.5. Double spaces of dimension five.} Assume now that
$M=\mathop{\rm dim}V=5$. Let the singular point $o\in V^+$ be an
isolated centre of non log canonical singularities of the pair
$(V^+,\frac{1}{m}\Sigma^+)$. The fiber $F\ni o$ is a Cartier
divisor, so that the point $o$ is the centre of a non log
canonical singularity of the pair
\begin{equation}\label{a3}
(F,\frac{1}{m}D_F),
\end{equation}
where $D_F\in\Sigma^+|_F$ is a general divisor, $D_F\sim -mK_F$.
By the arguments of Sec. 1.2, the point $o$ is an isolated centre
of non log canonical singularities of that pair. If the quadratic
singularity $o\in F$ is of rank 4 or 5, we argue as for $M\geq 6$
and come to a contradiction. Since by the conditions of general
position the rank of the quadratic point $o\in F$ is at least 3,
we assume that it is equal to 3.

The variety $F$ is realized as the double cover
$$
\sigma_F\colon F\to{\mathbb P}^4,
$$
generated by the morphism $\sigma$. By Proposition 1.8, proved
below, the conditions of general position imply (see Proposition
1.7) that the pair (\ref{a3}) is log canonical for $m=1$. Now
Proposition 1.1 for $M=5$ comes from the following
claim.\vspace{0.1cm}

{\bf Proposition 1.5.} {\it For any effective Cartier divisor
$D\sim-mK_F$ on $F$ the pair
\begin{equation}\label{a4}
(F,\frac{1}{m}D),
\end{equation}
is log canonical.}\vspace{0.1cm}

{\bf Proof} will be given by induction on $m\geq 2$ (as we
mentioned above, for $m=1$ the claim of the proposition is true).
It is sufficient to show that non log canonicity of the pair
(\ref{a4}) implies non log canonicity of a similar pair with a
smaller value of the parameter $m\geq 2$.

We may assume that the point $o$ is an isolated centre of non log
canonical singularities of the pair (\ref{a4}). For a generic
surface $S\ni o$ (a section of the germ $o\in F$ by two generic
hyperplanes with respect to some projective embedding) the pair
$$
(S,\frac{1}{m}D_S),
$$
where $D_S=D|_S$, is not log canonical at the point $o$ by
inversion of adjunction. On the other hand, the point $o$ is an
isolated centre of non log canonical singularities of that pair:
in the opposite case on $F$ there is a divisor $T$ such that
$$
D=aT+D_1,
$$
where $a>m$ and $D_1$ is effective, which is impossible. The
singularity $o\in S$ is a non-degenerate quadratic point.

Let
$$
\psi\colon F^{\sharp}\to F\quad \mbox{and}\quad
\bar{\psi}\colon{\mathbb P}^{\sharp}\to{\mathbb P}^4
$$
be the blow ups of the points $o\in F$ and
$p=\sigma_F(o)\in{\mathbb P}^4$. Denote the exceptional divisors
of the blow ups $\psi$ and $\bar{\psi}$ by the symbols
$E^{\sharp}$ and $\bar{E}^{\sharp}$, respectively. Obviously,
$\bar{E}^{\sharp}\cong{\mathbb P}^3$, and $\sigma_F$ extends to a
double cover
$$
\sigma_{\sharp}\colon F^{\sharp}\to{\mathbb P}^{\sharp},
$$
which on the level of exceptional divisors gives a double cover
$$
\sigma_E=\sigma_{\sharp}|_{E^{\sharp}}\colon E^{\sharp}
\to\bar{E}^{\sharp}.
$$

Set
$$
\psi^*D=D^{\sharp}+\nu E^{\sharp},
$$
where $D^{\sharp}$ is the strict transform. If $\nu>m$, then, as
above, we get $\mathop{\rm mult}_oD>2m$, which is impossible. For
this reason we assume that $\nu\leq m$. Applying Proposition 1.6,
which is proved below, to the pair $(S,\frac{1}{m}D_S)$, we
conclude that on the quadric $E^{\sharp}$ there is a plane
$P\cong{\mathbb P}^2$, such that the centre of any non log
canonical of the pair $(S,\frac{1}{m}D_S)$ on the strict
transform $S^{\sharp}\subset F^{\sharp}$ is a point $P\cap
S^{\sharp}$. Obviously, $\sigma_E(P)$ is a plane in
$\bar{E}^{\sharp}\cong{\mathbb P}^3$.

Let $Q\subset{\mathbb P}^4$ be the only hyperplane, such that
$$
Q^{\sharp}\cap\bar{E}^{\sharp}=\sigma_E(P)
$$
(as always, $Q^{\sharp}\subset{\mathbb P}^{\sharp}$ is the strict
transform). Set
$$
\Pi=\sigma^{-1}_F(Q)\subset F.
$$
The divisor $\Pi$ is irreducible, and moreover,
$$
\Pi^{\sharp}\cap E^{\sharp}\supset P.
$$
Now write down
$$
D=a\Pi+D^*,
$$
where $a\in{\mathbb Z}_+$ and $D^*$ is effective and does not
contain $\Pi$ as a component.\vspace{0.1cm}

{\bf Lemma 1.3.} {\it The inequality $a\geq 1$ holds. The pair
\begin{equation}\label{a5}
(F,\frac{1}{m^*}D^*)
\end{equation}
is not log canonical at the point $o$, where
$m^*=m-a$.}\vspace{0.1cm}

{\bf Proof.} By the conditions of general position, the pair
$(F,\Pi)$ is log canonical. Now by linearity we conclude that the
pair (\ref{a5}) is not log canonical at the point $o$, and by the
arguments of Sec. 1.2 this point is an isolated centre of non log
canonical singularities. This proves the second claim of the
lemma.

However, it is true for $a=0$ in a trivial way: $m^*=m$ and
$D^*=D$. Let us show that in fact $a\geq 1$. Indeed, by
Proposition 1.6, the inequality
\begin{equation}\label{a6}
\mathop{\rm mult}\nolimits_PD^{\sharp}+2\nu>2m
\end{equation}
holds. If $a=0$, then $\Pi$ is not contained in the support of
the divisor $D$, so that the effective cycle
$$
Y=(\Pi\circ D)
$$
of codimension two on $F$ is well defined. By (\ref{a6}) we get
$$
\mathop{\rm mult}\nolimits_oY>2m,
$$
but at the same time $\mathop{\rm deg}Y=2m$. Contradiction.
Q.E.D. for the lemma.\vspace{0.1cm}

Since for the pair (\ref{a5}), where $D^*\sim-m^*K_F$, the point
$o\in F$ is an isolated centre of a non log canonical singularity
and $m^*<m$, we apply the induction hypothesis and complete the
proof of Proposition 1.5.\vspace{0.1cm}

Now Proposition 1.1 is proven for $M=5$ as well.\vspace{0.3cm}


{\bf 1.6. Non log canonical singularities over a singular point
of the surface.} Let us consider the following local situation.
Let $o\in S$ be a germ of a non-degenerate double point on a
surface $S$ (that is, a germ, analytically isomorphic to the germ
$(0,0,0)\in\{x^2+y^2+z^2=0\}\subset{\mathbb C}^3)$. Let
$$
\varphi\colon S^+\to S
$$
be the blow up of the double point $o$, $E=\varphi^{-1}(o)$ the
exceptional conic. Assume that $C$ is an effective 1-cycle on $S$,
and for some positive $m$ the pair
\begin{equation}\label{a7}
(S,\frac{1}{m}C)
\end{equation}
is not log canonical at the point $o$, but log canonical outside
this point. Define the number $\nu\in{\mathbb Z}_+$ by the
relation
$$
C^+\sim\varphi^*C-\nu E,
$$
where $C^+$ is the strict transform of the 1-cycle $C$ on $S^+$.
Similar to Proposition 3 in [26], for the double point we
have\vspace{0.1cm}

{\bf Proposition 1.6.} {\it There exists a point $q\in E$ such
that}
\begin{equation}\label{a8}
2\nu+\mathop{\rm mult}\nolimits_q C^+>2m.
\end{equation}
\vspace{0.1cm}

{\bf Proof.} Note that if the exceptional conic $E$ makes itself a
non log canonical singularity of the pair $(S,\frac{1}{m}C)$, then
the inequality $\nu>m$ holds, that is, (\ref{a8}) holds for any
point $q\in E$. If $\nu\leq m$, then the connectedness principle
implies that the centre of any non log canonical singularity of
the pair (\ref{a7}) is some uniquely defined (by the pair) point
$q\in E$. We will prove that the inequality (\ref{a8}) holds for
that point.

Let
\begin{equation}\label{a9}
\varphi_i\colon S_i\to S_{i-1},
\end{equation}
$i=1,\dots,N$, be the sequence of blow ups of the points which are
the centres of a fixed non log canonical singularity of the
(\ref{a7}). More precisely, let
$$
\beta\colon\widetilde{S}\to S^+
$$
be some birational morphism, $E^+\subset\widetilde{S}$ an
irreducible exceptional curve, realizing the non log canonical
singularity of the pair (\ref{a7}), that is, the log Noether-Fano
inequality holds:
$$
\nu_{E^+}(C)>m(a(E^+,S)+1).
$$
By assumption, $\beta(E^+)$ is the point $q\in E$. Let us define a
sequence of blow ups (\ref{a9}), setting
$$
S_0=S,\,\,S_1=S^+,
$$
$\varphi_i$ blows up the point
$$
x_{i-1}=\mathop{\rm centre}(E^+,S_{i-1}),
$$
$i=1,\dots,N$ (so that $x_0=o$, $x_1=q$), and the last exceptional
curve
$$
E_N=\varphi^{-1}_N(x_{N-1})\subset S_N
$$
realizes $E^+$. As usual, we denote the exceptional curves by the
symbols
$$
E_i=\varphi^{-1}_i(x_{i-1})\subset S_i,
$$
so that $E_1=E$. Set
$$
\nu_i=\mathop{\rm mult}\nolimits_{x_{i-1}}C^{i-1}\in{\mathbb Z}_+,
$$
where $C^{i-1}$ is the strict transform of the cycle $C$ on
$S_{i-1}$, $i=2,\dots,N$, and $\nu_1=\nu$. The log Noether-Fano
inequality now is re-writed in the traditional form:
\begin{equation}\label{a10}
p_1\nu+\sum^N_{i=2}p_i\nu_i>m\left(\sum^N_{i=2}p_i+1\right),
\end{equation}
where $p_i=p_{Ni}$ is the number of paths from the vertex $N$ to
the vertex $i$ in the graph of the sequence of blow ups
(\ref{a9}). Note that in (\ref{a10}) in the right hand past there
is no component with $i=1$, because the discrepancy of $E$ is
zero. The multiplicities $\nu_i$ satisfy the system of linear
inequalities
\begin{equation}\label{a11}
\nu_i\geq\sum_{j\to i}\nu_j
\end{equation}
for $i=2,\dots,N$ and, besides,
\begin{equation}\label{a12}
2\nu\geq\sum_{j\to 1}\nu_j.
\end{equation}
Finally, $\nu_N\geq 0$. Non-negativity of the other numbers
$\nu_1,\nu_2,\dots,\nu_{N-1}$ follows from (\ref{a11},\ref{a12}).
For simplicity let us denote by the symbol the $(\ref{a10})^*$
{\it non-strict} log Noether-Fano inequality, that is, the
inequality (\ref{a10}), in which the sign $>$ is replaced by
$\geq$. Finally, by the symbol ${\cal L}$ denote the system of
non-strict linear inequalities $(\ref{a10})^*$,
(\ref{a11},\ref{a12}) and $\nu_N\geq 0$.

Let us show that if the set of real numbers
$$
\nu_1,\dots,\nu_N
$$
satisfies the system  ${\cal L}$, then the estimate
\begin{equation}\label{a13}
2\nu_1+\nu_2\geq 2m
\end{equation}
holds. This immediately implies the inequality (\ref{a8}).

Set $\Lambda\subset{\mathbb R}^N$ to be the convex subset defined
by the system ${\cal L}$. Obviously, the linear function
$2\nu_1+\nu_2$ is bounded from below on $\Lambda$, and moreover,
the infimum is a monimum, attained at some point
$$
v=(\theta_1,\dots,\theta_N)\in{\mathbb R}^N.
$$
We may assume that the point $v$ is one of the vertices of the set
$\Lambda$, that is, that $N$ inequalities from the system ${\cal
L}$ become equalities at that point.\vspace{0.1cm}

{\bf Lemma 1.4.} {\it Assume that $\theta_N=0$. Then there exists
$K\in\{2,\dots,N-1\}$ such that the inequality}
\begin{equation}\label{a14}
\sum^K_{i=1}p_{Ki}\theta_i>m\left(\sum^K_{i=2}p_{Ki}+1\right)
\end{equation}
{\it holds.} \vspace{0.1cm}

Arguing by induction on the length $N$ of the resolution of the
singularity $E^+$, we obtain the estimate (\ref{a13}) in the case
$\theta_N=0$.\vspace{0.1cm}

{\bf Proof of the lemma 1.4.} Since the curve
$$
\bigcup^N_{i=1}E^N_i
$$
is a normal crossing divisor on a non-singular surface, the vertex
$N$ is connected by arrows with one or two vertices: always
$$
N\to N-1
$$
and, possibly, $N\to L$ for some $L\leq N-2$. The first case is
trivial; we will consider the second one (our arguments, with
simplifications, prove the claim of the lemma in the first case as
well). Setting $p_{ij}=0$ for $i<j$, by definition of the
incidence graph we get
$$
p_{Ni}=p_{N-1,i}+p_{L,i}
$$
for any $i\leq N-1$. Now, taking into account that $\theta_N=0$,
we can re-write the inequality $(\ref{a10})^*$ in the following
way:
$$
\left(\sum^{N-1}_{i=1}p_{N-1,i}\theta_i
-m\left(\sum^{N-1}_{i=2}p_{N-1,i}+1\right)\right)+
$$
$$
+\left(\sum^L_{i=1}p_{Li}\theta_i-
m\left(\sum^L_{i=2}p_{Li}+1\right)\right)\geq0,
$$
which implies that the inequality (\ref{a14}) holds either for
$K=N-1$, or for $K=L$ (or for both these values). Q.E.D. for the
lemma.\vspace{0.1cm}

Thus we may assume that $\theta_N>0$ and therefore for the vector
$v$ the inequalities $(\ref{a10})^*$, (\ref{a11}) and (\ref{a12})
are equalities. It follows that for $\theta=\theta_N$
$$
\theta_i=p_i\theta
$$
for $i=2,\dots,N$, and $\theta_1=\frac12p_1\theta$, so that
$\theta$ can be found from the equation
$$
\left(\frac12p^2_1+\sum^N_{i=2}p^2_i\right)\theta=
\left(\sum^N_{i=2}p_i+1\right)m.
$$
The value of the linear function $2\nu_1+\nu_2$ at the vector $v$
is $(p_1+p_2)\theta$, so that the inequality (\ref{a13}) comes
from the following combinatorial fact.\vspace{0.1cm}

{\bf Lemma 1.5.} {\it The following inequality holds:}
\begin{equation}\label{a15}
(p_1+p_2)\left(\sum^N_{i=2}p_i+1\right)\geq
p^2_1+2\sum^N_{i=2}p^2_i.
\end{equation}
\vspace{0.1cm}

{\bf Proof} will be given by induction on the number $N$ of
vertices in the incidence graph. If $N=2$, then $p_1=p_2=1$ and
(\ref{a15}) holds. Furthermore, the inequality
$$
(p_2+p_3)\left(\sum^N_{i=3}p_i+1\right)\geq
p^2_2+2\sum^N_{i=3}p^2_i
$$
holds by the induction hypothesis. In order to obtain (\ref{a15}),
it is sufficient to show the estimate
\begin{equation}\label{a16}
(p_1+p_2)p_2+(p_1-p_3)\left(\sum^N_{i=3}p_i+1\right)\geq
p^2_1+p^2_2.
\end{equation}
In [33, Lemma 1.6] it was proved that the inequality
$$
\sum^N_{i=3}p_i+1\geq p_1
$$
holds, so that (\ref{a16}) follows from the estimate
$$
(p_1+p_2)p_2+(p_1-p_3)p_1\geq p^2_1+p^2_2,
$$
which is obvious, since $p_2\geq p_3$.

Q.E.D. for Lemma 1.5 and Proposition 1.6.\vspace{0.3cm}


{\bf 1.7. Additional conditions of general position for $M=5$.}
Here we assume that $M=5$. For an arbitrary point $p\in W$ set
$$
T(p)=\sigma^{-1}(T_pW).
$$
Let
$$
\varphi\colon T^+(p)\to T(p)
$$
be the blow up of the isolated double point $o=\sigma^{-1}(p)$
with the exceptional divisor $E(p)$, a three-dimensional quadric
in ${\mathbb P}^4$. Set
$$
Y_i=\{p\in W|\mathop{\rm rk}E(p)=i\}\subset W.
$$
\vspace{0.1cm}

{\bf Proposition 1.7.} {\it For a generic variety $V$ we have
$\mathop{\rm rk}E(p)\geq 3$, that is, $Y_1=Y_2=\emptyset$.
Furthermore, $\mathop{\rm dim}Y_3=1$. For a point $p\in Y_3$ the
singularities of the variety $T(p)$ are of the following form:

1) for a point $p\in Y_3$ of general position on the line
$L=\mathop{\rm Sing}E(p)$ there are three distinct singular
points of the variety $T^+(p)$, which are non-degenerate
quadratic points, and on $E(p)$ the variety $T^+(p)$ has no other
singular points;

2) for a finite set of points $p\in Y_3$, which do not satisfy
the condition 1), on the line $L=\mathop{\rm Sing}E(p)$ there are
two distinct singular points $p_1$ and $p_2$ of the variety
$T(p)$. On $E(p)$ the variety $T(p)$ has no other singular
points. On of these points (say, $p_1$) is a non-degenerate
quadratic singularity. The point $p_2$ is an isolated quadratic
point of rank 4. Its blow up
$$
\varphi_{\sharp}\colon T^{\sharp}(p)\to T^+(p)
$$
has a unique singular point $p_3$ on the exceptional divisor
$E^{\sharp}=\varphi^{-1}_{\sharp}(p_2)$, which is the vertex of
the cone $E^{\sharp}$, and moreover, $p_3\in T^{\sharp}(p)$ is a
non-degerate quadratic point.}\vspace{0.1cm}

{\bf Proof:} an easy dimension count for the local equation
$$
y^2=q_2(z_1,z_2,z_3,z_4)+q_3(z_*)+\dots
$$
of the variety $T(p)$ at the point $o$. For a point $p\in Y_3$
singularities of the variety $T^+(p)$ correspond to the zeros of
the polynomial
$$
q_3(z_*)|_L
$$
on the vertex line $L$ of the quadric $E(p)$. If the three roots
are distinct, we get the case 1). If one of the roots is a double
root, we get the case 2), where the point $p_2\in L$ corresponds
to the double root. Simple calculations are left to the reader.
Q.E.D.\vspace{0.1cm}

{\bf Proposition 1.8.} {\it For a generic variety $V$, an
arbitrary point $p\in Y_3$ and an arbitrary hyperplane $R\subset
T_pW$, $R\ni p$, the pair
$$
(T(p),\Pi=\sigma^{-1}(R))
$$
is canonical at the point $o$.}\vspace{0.1cm}

{\bf Proof.} Since
$$
\mathop{\rm mult}\nolimits_o \Pi=a(E,T)=2,
$$
where $E=E(p)$, $T=T(p)$ for simplicity of notations, it is
sufficient to prove that the pair
$$
(T^+,\Pi^+),
$$
is canonical, where $\Pi^+\subset T^+$ is the strict transform of
the divisor $\Pi$.

Let
$$
y^2=q^*_2(z_1,z_2,z_3)+q^*_3(z_*)+\dots
$$
be the local equation of the variety $\Pi$. If $\mathop{\rm
rk}q^*_2=3$, then $o\in\Pi$ is an ordinary double point and there
is nothing to prove. If $\mathop{\rm rk}q^*_2=2$, then it is easy
to check that the singularities of $\Pi$ over the point $o$ are
resolved by a sequence of $k\leq 6$ blow ups of isolated quadratic
points of rank $\geq 3$. In that case it is also obvious that the
pair $(T,\Pi)$ is canonical.

Now assume that $\mathop{\rm rk}q^*_2=1$. It imposes 3 independent
conditions on $\Pi$ (more precisely, on the polynomial $f|_R$), so
that there is a 4-dimensional family of subspaces
$R\subset{\mathbb P}$, for which $\Pi$ satisfies this property.
Let
$$
E_{\Pi}=\Pi^+\cap E
$$
be the exceptional quadric of rank 2 in ${\mathbb P}^3$, that is,
a pair of planes, $L=\mathop{\rm Sing}E_{\Pi}$ the line of their
intersections.

If
$$
q^*_3|_L\equiv 0
$$
(this imposes 4 additional independent conditions on $\Pi$, so
that there are only finitely many such pairs), then let
$$
\varphi_L\colon\Pi_L\to\Pi^+
$$
be the blow up of the line $L$. It is easy to check that $\Pi_L$
has finitely many isolated double points, resolved by one blow up.
It is now obvious that the pair $(T,\Pi)$ is canonical, taking
into account that
$$
\mathop{\rm mult}\nolimits_L\Pi^+=a(L,T^+)=2.
$$

Assume that $q^*_3|_L\not\equiv 0$, that is, at a general point of
the line $L$ the variety $\Pi^+$ is non-singular. It is sufficient
to check that the pair $(T^+,\Pi^+)$ is canonical at the singular
points of the variety $\Pi^+$ on the line $L$. The explicit
computations in local coordinates (they are elementary but
tiresome and we omit them) show that singularities of the divisor
$\Pi^+$ are resolved by a sequence of blow ups of isolated
quadratic points (of rank $\geq 2$), which implies canonicity of
the pair $(T^+,\Pi^+)$.

Assume finally that $q^*_2\equiv 0$. This imposes 6 independent
conditions on $\Pi$, so that there is a one-dimensional family of
such subvarieties on $V$. The quadric $E_{\Pi}$ is the double
plane $2\Lambda=\{y=0\}$, however, the arguments of general
position imply that $\mathop{\rm mult}_{\Lambda}\Pi^+=1$ and,
moreover,
$$
C=\mathop{\rm Sing} \Pi^+=\{q^*_3|_{\Lambda}=0\}
$$
(we mean the singularities over the point $o$) is an irreducible
cubic curve, and moreover, if it is singular, then the only
singular point of that curve lies outside the line $L=\mathop{\rm
Sing}E(p)$. We get
$$
\mathop{\rm mult}\nolimits_C\Pi^+=a(C,T^+)=2.
$$
Blowing up of the cubic $C$ gives the variety $\Pi_C$, which is
non-singular over a generic point of $C$. It is easy to check that
for any point $p\in C\backslash L$ the pair $(T^+,\Pi^+)$ is
canonical (even terminal) at that point, either. Finally, the
variety $\Pi_C$ has only isolated quadratic points of rank $\geq
3$, so that it is easy to check that the pair $(T^+,\Pi^+)$ is
canonical also over the points $q\in C\cap L$. Note that for the
blow up
$$
\varphi_q\colon T_q\to T^+
$$
of a point $q\in C\cap L$, with $E_q=\varphi^{-1}_q(q)$ the
exceptional divisor, we get
$$
\mathop{\rm mult}\nolimits_q\Pi^+=a(E_q,T^+)=2,
$$
$E_q$ is a quadric of rank $\geq 4$, so that $E_q$ realizes one
more canonical, but not terminal singularity of the pair
$(T,\Pi)$. This completes our proof of Proposition
1.8.\vspace{0.1cm}

{\bf Remark 1.1.} Since the multiplicities of singular points and
subvarieties are equal to 2, the proof of Proposition 1.8 reduces
to checking that if the strict transform of the divisor $\Pi$ has
a curve of singular points then the strict transform $T$ is
non-singular at the generic point of this curve, and that the
singularities of the strict transform of $\Pi$ are at most
one-dimensional.\vspace{0.3cm}


\section{Exclusion of maximal subvarieties\\ of codimension two}

In this section, we start to prove Proposition 0.3: we show that,
except for the preimage $\sigma^{-1}(P)$, where $P\subset{\mathbb
P}$ is a linear subspace of codimension two, no subvariety of
codimension two can make a maximal subvariety of the system
$\Sigma$.\vspace{0.3cm}

{\bf 2.1. Set up of the problem.} The following claim is true.
\vspace{0.1cm}

{\bf Proposition 2.1.} {\it If an irreducible subvariety
$B\subset V$ of codimension two is maximal for the movable linear
system $\Sigma\subset|\,2nH|$, that is, the inequality
$\mathop{\rm mult}_B\Sigma>n$ holds, then
$B=\sigma^{-1}(\bar{B})$, where $\bar{B}\subset{\mathbb P}$ is a
linear subspace of codimension two.}\vspace{0.1cm}

{\bf Proof.} The self-intersection $Z=(D_1\circ D_2)$,
$D_i\in\Sigma$, of the linear system $\Sigma$ is of $H$-degree
$8n^2$ and contains the subvariety $B$ with multiplicity strictly
higher than $n^2$. Therefore, $\mathop{\rm deg}B\leq 7$. It is
necessary to show that only one of these possibilities realize:
$\mathop{\rm deg}B=2$, and moreover, $\bar{B}=\sigma(B)$ is a
$(M-2)$-plane in ${\mathbb P}$, that is, the double cover
$\sigma^{-1}(\bar{B})\to\bar{B}$ is irreducible.

Note that for $\mathop{\rm dim}V=M\geq 5$ we have $A^2V={\mathbb
Z}H^2$, so that only three possibilities occur, $B\sim H^2$ or
$2H^2$ or $3H^2$. In particular, $\mathop{\rm deg}B\in\{2,4,6\}$.
However, we exclude below maximal subvarieties of codimension two
for $M=4$, either.

Let us exclude, first of all, the case $\mathop{\rm deg}B=1$. Here
$M=4$, so that $\bar{B}\subset{\mathbb P}={\mathbb P}^4$ is a
2-plane, and moreover, the double cover
$\sigma^{-1}(\bar{B})\to\bar{B}$ is reducible. Therefore, the
curve $\bar{B}\cap W$ is everywhere non-reduced (it is a cubic
curve with multiplicity two). This is impossible by generality of
the hypersurface $W$, see Sec. 0.4.

If $B=\sigma^{-1}(\bar{B})$, then $\mathop{\rm
deg}B\in\{2,4,6\}$. Assume that $\mathop{\rm deg}B\in\{4,6\}$,
that is, $\mathop{\rm deg}\bar{B}\in\{2,3\}$. Let us show that
these cases do not realize. Indeed, let $L\subset{\mathbb P}$ be
a generic secant line of the subvariety $\bar{B}\subset{\mathbb
P}$. By generality, the curve $C=\sigma^{-1}(L)$ is non-singular
and irreducible, and such curves sweep out at least a divisor on
$V$, so that $C\not\subset\mathop{\rm Bs}\Sigma$. For a general
divisor $D\in\Sigma$ we get $C\not\subset D$ and $(C\cdot D)=4n$.
On the other hand, let $p_1\neq p_2$ be the points of intersection
$L\cap\bar{B}$, then (by generality)
$\sigma^{-1}(p_i)=\{p_{i1},p_{i2}\}$, $i=1,2$, where $p_{ij}$ are
four distinct points on $B$. For this reason,
$$
4n=(C\cdot D)\geq\sum\limits_{i,j}(C\cdot D)_{p_{ij}}>4n.
$$
Contradiction.

Thus if $B=\sigma^{-1}(\bar{B})$, then $\bar{B}\subset{\mathbb P}$
is a $(M-2)$-plane, which is exactly what we need.

Starting from this moment, we assume that
$\sigma^{-1}\bar{B}=B\cup B'$ breaks into two irreducible
components and
$$
\mathop{\rm deg}B=\mathop{\rm deg}\bar{B}\in\{2,3,4,5,6,7\}.
$$
We show below that none of these cases realizes. Let us describe
first of all the main technical tools that will be used for their
exclusion.\vspace{0.3cm}


{\bf 2.2. Conics on the variety $\bar B$.} Let $C\subset\bar{B}$
be an irreducible conics, $P=<C>$ its linear span (a 2-plane).
Assume that $C\not\subset W$, the curve $W\cap P$ is reduced and
the two finite sets
$$
C\cap W \quad \mbox{and} \quad \mathop{\rm Sing}(W\cap P)
$$
are disjoint. Set $S=\sigma^{-1}(P)$, this is an irreducible
surface with a finite set of singular points $\sigma^{-1}(W\cap
P)$. Let $C_+$ and $C_-$ be components of the curve
$\sigma^{-1}(C)=C_+\cup C_-$, where $C_+\subset B$, $C_-\subset
B'$. \vspace{0.1cm}

{\bf Lemma 2.1.} {\it The surface $S$ is contained in the base
set $\mathop{\rm Bs}\Sigma$.}\vspace{0.1cm}

{\bf Proof.} Assume the converse. Then for a general divisor
$D\in\Sigma$ we get $S\not\subset D$, so that $(D\circ S)$ is an
effective curve on $S$, containing $C_+$ with some multiplicity
$$
\nu_+\geq\mathop{\rm mult}\nolimits_B\Sigma> n
$$
and $C_-$ with some multiplicity $\nu_-\in{\mathbb Z}_+$. Let
$H_S=H\,|\,_S$ be the class of a hyperplane section of $S$. By
what we have said,
\begin{equation}\label{b1}
((2nH_S-\nu_+C_+-\nu_-C_-)\cdot C_{\pm})\geq 0.
\end{equation}
Note that by assumption the curves $C_{\pm}$ do not contain
singular points of the surface $S$, so that the local
intersection numbers $(C_+\cdot C_-)_x$ are equal to
$\frac12(C\cdot W)_{\sigma(x)}$ and therefore
$$
(C_+\cdot C_-)=\frac12(C\cdot W)=2(M-1).
$$
Furthermore, $C_++C_-\sim 2H_S$, whence we obtain
$$
(C^2_+)=(C^2_-)=2(3-M).
$$
Therefore, the inequalities (\ref{b1}) take the form of linear
inequalities
\begin{equation}\label{b2}
\begin{array}{c}
4n+2(M-3)\nu_+-2(M-1)\nu_-\geq 0,\\
4n-2(M-1)\nu_++2(M-3)\nu_-\geq 0,
\end{array}
\end{equation}
whence we get $\nu_{\pm}\leq n$. Contradiction. Q.E.D. for the
lemma.\vspace{0.1cm}

{\bf Corollary 2.1.} {\it The following inequality holds:}
$\mathop{\rm deg}B\geq 4$.\vspace{0.1cm}

{\bf Proof.} We have to exclude two cases: $\mathop{\rm deg}B=2$
and $\mathop{\rm deg}B=3$. First assume that $\mathop{\rm
deg}B=2$. Applying Lemma 2.1 to the irreducible conic
$C=\bar{B}\cap P$, where $P\subset<\bar{B}>$ is a generic 2-plane
in the linear span of $\bar{B}$, we get that
$\sigma^{-1}(P)\subset\mathop{\rm Bs}\Sigma$. Therefore,
$$
\sigma^{-1}(<\bar{B}>)\subset\mathop{\rm Bs}\Sigma,
$$
which is impossible, since $<\bar{B}>$ is a divisor in ${\mathbb
P}$. If $\mathop{\rm deg}B=3$, the arguments are similar: the
variety $\bar{B}$ is swept out by conics, and moreover a generic
conic $C\subset\bar{B}$ satisfies the assumptions of Lemma 2.1.
The linear spans $P=<C>$ of those conics sweep out at least a
divisor in ${\mathbb P}$, which again contradicts the fact that
the linear system $\Sigma$ is movable. Q.E.D. for the
corollary.\vspace{0.3cm}


{\bf 2.3. The secant lines of the variety $\bar B$.} Let
$C\subset\bar{B}$ be an irreducible curve, not contained in $W$.
Let $x\in{\mathbb P}$ be a point, satisfying the following
conditions of general position:
\begin{itemize}
\item $x\not\in C$,
\item for any point $p\in C\cap W$ the line $L=<x,p>$, connecting
the points $x$ and $p$, intersects the hypersurface $W$
transversally at the point $p$ and is not a secant line of the
curve $C$, that is, $C\cap L=\{p\}$.
\end{itemize}

Consider the cone $\Delta=\Delta(x,C)$ with the vertex at the
point $x$ and the base $C$. Set $S=\sigma^{-1}(\Delta)$, it is an
irreducible surface. Let $C_+$ and $C_-$ once again be the
components of the curve $\sigma^{-1}(C)=C_+\cup C_-$, where
$C_+\subset B$, $C_-\subset B'$. By the assumptions above, all
the points of intersection of the curves $C_+$ and $C_-$ are
smooth points of the surface $S$. Obviously,
$$
(C_+\cdot C_-)=\frac12(C\cdot W)=(M-1)\mathop{\rm deg}C.
$$
Furthermore, it is well known [1], that on the cone $\Delta$ the
curve $C$ is numerically equivalent to the hyperplane section.
Thus on the surface $S$
$$
C_++C_-\equiv H_S=\sigma^*H_{\Delta},
$$
where $H_{\Delta}$ is the hyperplane section of the cone $\Delta$.
From here we get that
$$
(C^2_+)=(C^2_-)=-(M-2)\mathop{\rm deg}C.
$$
The restriction $\Sigma_S=\Sigma\,|\,_S$ of the system $\Sigma$
on $S$ is a non-empty linear system of curves, containing
$C_{\pm}$ with the multiplicity $\nu_{\pm}$, respectively, where
$\nu_+\geq\mathop{\rm mult}_B\Sigma>n$. Therefore,
$$
((2nH_S-\nu_+C_+-\nu_-C_-)\cdot C_{\pm})\geq 0,
$$
which yields the system of linear inequalities
\begin{equation}\label{b3}
\begin{array}{c}
2n-(M-1)\nu_++(M-2)\nu_-\geq 0,\\
2n+(M-2)\nu_+-(M-1)\nu_-\geq 0.
\end{array}
\end{equation}

From here we immediately get\vspace{0.1cm}

{\bf Proposition 2.2.} {\it The following estimate holds}
\begin{equation}\label{b4}
\mathop{\rm
mult}\nolimits_{B'}\Sigma>\frac{M-3}{M-2}n\geq\frac{n}{2}.
\end{equation}
\vspace{0.1cm}

{\bf Proof.} For a general choice of the vertex $x$ of the cone
$\Delta$ we obtain $\nu_+=\mathop{\rm mult}_B\Sigma$,
$\nu_-=\mathop{\rm mult}_{B'}\Sigma$, whereas the inequality
$$
\nu_->\frac{M-3}{M-2}n
$$
follows directly from (\ref{b3}). Q.E.D. for the
proposition.\vspace{0.1cm}

Note that the estimate (\ref{b4}) is the stronger, the higher is
$M$. Proposition 2.2 makes it possible to exclude the case
$\mathop{\rm deg}B=7$ straightaway.\vspace{0.1cm}

{\bf Proposition 2.3.} {\it The case $\mathop{\rm deg}B=7$ does
not realize.}\vspace{0.1cm}

{\bf Proof.} Assume the converse: $\mathop{\rm deg}B=7$. Let
$D_1,D_2\in\Sigma$ be general divisors, $Z=(D_1\circ D_2)$ the
self-intersection of the system $Z$. We obtain the inequality
$$
8n^2=\mathop{\rm deg}Z\geq 7((\mathop{\rm
mult}\nolimits_B\Sigma)^2+(\mathop{\rm
mult}\nolimits_{B'}\Sigma)^2)> 7 \cdot\frac{5}{4}n^2,
$$
which is impossible. Contradiction. Q.E.D. for the
proposition.\vspace{0.1cm}

Note that, repeating this argument word for word, we exclude the
case $\mathop{\rm deg}B=6$ for $M\geq 5$: for the multiplicity of
the subvariety $B'$ Proposition 2.2 gives the estimate
$\mathop{\rm mult}_{B'}\Sigma>2n/3$, so that
$$
8n^2>6\cdot(1+\frac{4}{9})n^2=\frac{26}{3}n^2,
$$
which is impossible once again. \vspace{0.3cm}


{\bf 2.4. Three-secant lines of the variety $\bar B$.} Thus it
remains to exclude three cases: $\mathop{\rm deg}B=4,5,6$,
whereas in the two latter cases $\mathop{\rm dim}V=M=4$. We will
need another simple construction. Let $L\subset{\mathbb P}$ be a
3-secant line of the variety $\bar{B}$, that is, a line that
intersects $\bar{B}$ at (at least) three points outside
$W$.\vspace{0.1cm}

{\bf Proposition 2.4.} {\it If the curve $\sigma^{-1}(L)=C$ is
irreducible, then $C\subset\mathop{\rm Bs}\Sigma$. If $C=C_+\cup
C_-$ is reducible, then at least one of the components $C_{\pm}$
is contained in $\mathop{\rm Bs}\Sigma$.}\vspace{0.1cm}

{\bf Proof.} Let $D\in\Sigma$ be an arbitrary divisor. The curve
$C$ intersects $D$ at at least 6 points. The total multiplicity
of $D$ at those points is at least
$$
3(\mathop{\rm mult}\nolimits_B\Sigma+\mathop{\rm
mult}\nolimits_{B'}\Sigma)>\frac{9}{2}n,
$$
whereas $C\cdot D =4n$. Therefore $L\subset\sigma(D)$, which is
what we need. Q.E.D. for the proposition.\vspace{0.1cm}

Therefore, if the subvariety $\bar{B}\subset{\mathbb P}$ has
sufficiently many 3-secant lines (more precisely, if they sweep
out at least a divisor on ${\mathbb P}$), then the subvariety
$B\subset V$ can not be maximal since the linear system $\Sigma$
is movable.\vspace{0.1cm}

{\bf Remark 2.1.} The claim of Proposition 2.4 (and its proof)
remain true if the line $L$ intersects $\bar{B}$ at two distinct
points outside $W$, and in one of them, say, $\bar{x}\in L\cap
\bar{B}$, is tangent to $\bar{B}$. In that case the curve
$C=\sigma^{-1}(L)$ is tangent to $B$ and $B'$ at the points $x$,
$x'$, respectively, where $\sigma^{-1}(\bar{x})=\{x,x'\}$, $x\in
B$, $x'\in B'$, and it is easy to see that the local intersection
numbers satisfy the inequalities
$$
(C\cdot D)_x\geq 2\mathop{\rm mult}\nolimits_B \Sigma,\quad
(C\cdot D)_{x'}\geq 2\mathop{\rm mult}\nolimits_{B'} \Sigma,
$$
which makes it possible to argue in word for word the same way as
in the case of three distinct points. In the sequel, when
speaking about 3-secant lines, we will include the limit case of
tangency without special reservations.

As a first application of the construction of Proposition 2.4 we
exclude the case $\mathop{\rm deg}B=4$ (the dimension $M\geq 4$
is arbitrary).\vspace{0.1cm}

{\bf Proposition 2.5.} {\it The case $\mathop{\rm deg}B=4$ does
not take place.}\vspace{0.1cm}

{\bf Proof.} Assume the converse. Let $P\subset{\mathbb P}$ be a
generic 3-plane. For the irreducible curve $B_P=\bar{B}\cap P$ in
${\mathbb P}^3$ the four cases are possible:\vspace{0.1cm}

1) $B_P\subset R$ is a plane curve, $R={\mathbb P}^2$ is a plane
in $P$;

2) $B_P=Q_1\cap Q_2$ is a smooth elliptic curve, the intersection
of quadrics $Q_1$ and $Q_2$;

3) $B_P$ is a smooth rational curve;

4) $B_P$ has a double point.\vspace{0.1cm}

The case 1) does not realize, because any line $L\subset R$ is a
4-secant line. Proposition 2.4 implies that the entire surface
$\sigma^{-1}(R)$ is contained in the base set $\mathop{\rm
Bs}\Sigma$. This is impossible, since $P$ is a generic 3-plane.

In the case 2) we come to a contradiction in exactly the same way
as in the proof of Lemma 2.1. Namely, let $Q$ be a generic
quadric, containing the curve $B_P$. On the surface $Q$ we get
$B_P\sim 2H_Q$, where $H_Q$ is the plane section. Set
$$
S=\sigma^{-1}(Q),\quad\sigma^{-1}(B_P)=C_+\cup C_-,\quad
C_+\subset B,\quad C_-\subset B',
$$
so that on $S$ we have $C_++C_-\sim 2H_S$, where $H_S=H\,|\,_S$
is the class of a hyperplane section. Now we argue in exactly the
same way as in Lemma 2.1 and obtain the inequalities (\ref{b1}),
which give the system of linear inequalities (\ref{b2}). This
contradiction excludes the case 2). Note that of key importance
(as in Lemma 2.1) is the fact that the curve $B_P$ is equivalent
to {\it two} hyperplane sections of the surface $Q$. In a general
case, a curve can be embedded into a surface as a hyperplane
section (with multiplicity one), which gives just some estimate
for the multiplicity of the second component $B'$, but does not
allow to get a contradiction in one step.

Consider the case 3). Let $x\in B_P$ be a point of general
position, $\pi_x\colon B_P\to{\mathbb P}^2$ the projection from
the point $x$. The image $\pi_x(B_P)\subset{\mathbb P}^2$ is a
rational cubic curve with a double point. Therefore, the curve
$B_P$ has a 3-secant line, passing through the point $x$. Since
$P$ is a 3-plane and $x\in B_P$ is a general point, we apply
Proposition 2.4 and obtain a contradiction.

Consider the case 4). The curve $B_P$ has a unique double point.
This implies that the variety $\bar{B}$ contains a $(M-3)$-plane
$\Pi$ of double points. Let $L\subset \Pi$ be a generic line,
$\Lambda\supset L$ a generic 3-plane, containing $L$. Now the
curve $B_{\Lambda}=B\cap \Lambda$ is a quartic in ${\mathbb
P}^3$, containing the line $L$ with multiplicity 2. Therefore,
$$
B_{\Lambda}=C_{\Lambda}+2L,
$$
where $C_{\Lambda}$ is a (in the general case irreducible) conic.
The variety $\bar{B}$ is swept out by the conics $C_{\Lambda}$.
Now we apply Lemma 2.1 and obtain a contradiction.

Q.E.D. for Proposition 2.5.\vspace{0.3cm}


{\bf 2.5. Exclusion of the cases $\mathop{\rm deg} B=5\,
\mbox{and}\, 6$.} Recall that we may assume that $\mathop{\rm
dim}V=M=4$ (although our arguments work in arbitrary dimension).
Let $P\subset{\mathbb P}$ be a generic hyperplane (that is, a
3-plane), $B_P=\bar{B}\cap P$ an irreducible curve. We may assume
that the linear span of the curve $B_P$ is $P={\mathbb P}^3$
(otherwise we argue as in the case 1) for $\mathop{\rm deg}B=4$).
Besides, the curve $B_P$ does not contain singular points of
multiplicity $\geq 3$ if $\mathop{\rm deg}B=5$ and of multiplicity
$\geq 4$ if $\mathop{\rm deg}B=6$ (otherwise we argue as in the
case 4) for $\mathop{\rm deg}B=4$).

Assume that $\mathop{\rm deg}B=5$. It is easy to check that there
is 3-secant line through a generic point $x\in B_P$. Indeed, if
the curve $B_P$ is smooth, then the projection from the point $x$
realizes $B_P$ as a plane quartic $Q\subset{\mathbb P}^2$, which
can not be smooth: if the curve $Q$ were smooth, by Riemann-Roch
we would have got
$$
h^0(l_Q+x)-h^1(l_Q+x)=5+1-3=3,
$$
where $l_Q=L\cap Q$ is the section of $Q$ by a line
$L\subset{\mathbb P}^2$. Furthermore,
$$
h^1(l_Q+x)=h^0(-x)=0,
$$
whence $h^0(l_Q+x)=3$, but at the same time $l_Q+x$ is a plane
section of the smooth curve $B_P\subset{\mathbb P}^3$ and for that
reason $h^0(l_Q+x)\geq 4$. Contradiction. Therefore, the quartic
$Q$ is singular and $B_P$ has a 3-secant line, passing through the
point $x$. Now Proposition 2.4 gives a contradiction.

Therefore, the curve $B_P$ has $\delta\geq 1$ double points. Let
$p\in\mathop{\rm Sing}B_P$ be a double point. The projection
$\pi_p\colon{B_P}\to{\mathbb P}^2 $ from the point $p$ realizes
$B_P$ as a plane cubic with $\geq(\delta-1)\geq 0$ double points,
that is, a curve of genus $\leq 2-\delta$. On the other hand, the
projection $\pi_x\colon B_P\to{\mathbb P}^2$ from a generic point
$x\in B_P$ realizes $B_P$ as a plane quartic with
$\delta^*\geq\delta$ double points, that is, a curve of genus
$3-\delta^*$. Therefore, we get the inequality
$\delta^*\geq\delta+1$, that is, there is a 3-secant line
$L\subset P$ through the point $x$, that {\it does not contain}
the double points of the curve $B_P$. Now we can apply Proposition
2.4 and obtain a contradiction. This excludes the case
$\mathop{\rm deg}B=5$.

Assume that $\mathop{\rm deg}B=6$. The fact that the curve $B_P$
can not be smooth is proved as the similar fact for $\mathop{\rm
deg}B=5$. Assume that the point $p\in B_P$ is of multiplicity 3.
Comparing the curves
$$\pi_p(B_P)\subset{\mathbb P}^2 \quad \mbox{и}\quad
\pi_x(B_P)\subset{\mathbb P}^2,
$$
where $x$ is a generic point, we exclude this case by the same
arguments as in the case of a curve of degree 5 with
singularities. Thus we may assume that the curve $B_P$ contains
$\delta\geq 1$ double points and does not contain points of higher
multiplicity.

Now we argue in word for word the same way as for $\mathop{\rm
deg}B=5$: we compare the curve $\pi_p(B_P)\subset{\mathbb P}^2$ of
degree 4 with $\geq(\delta-1)\geq 0$ double points
($p\in\mathop{\rm Sing}B_P$ is one of the singular points) with
the curve $\pi_x(B_P)\subset{\mathbb P}^2$ of degree 5 with
$\delta^*\geq\delta$ double points. We get that $B_P$ has a
3-secant line, passing through the point $x\in B_P$ of general
position and not containing singular points of the curve $B_P$.
(If the curve $\pi_p(B_P)\subset{\mathbb P}^2$ is a conic, that
is, $\mathop{\rm deg}\pi_p=2$, then $B_P$ is contained in a
quadric cone with the vertex $p$, and moreover, we may assume that
$B_P$ has no other double points. In that case the genus of the
curve $B_P$ is easy to compute and we can show that there exists a
3-secant line through a point of general position.)  Applying
Proposition 2.4, we get a contradiction. The case $\mathop{\rm
deg}B=6$ is excluded.

This completes the proof of Proposition 2.1.\vspace{0.3cm}


\section{Exclusion of maximal singularities\\ with the centre
of codimension three}

In this section we continue the proof of Proposition 0.3: we
prove that the linear system $\Sigma$ has no maximal
singularities, the centre of which is a subvariety of codimension
three on $V$.\vspace{0.3cm}

{\bf 3.1. Set up of the problem. Exclusion of centres of degree
$\geq 2$.} Recall that the linear system $\Sigma$ has a maximal
singularity, the centre of which is an irreducible subvariety
$B\subset V$. In Sec. 2 we proved that if $\Sigma$ has no maximal
subvariety of the form $\sigma^{-1}(P)$, where $P\subset{\mathbb
P}$ is a linear subspace of codimension two, then $\Sigma$ has no
maximal subvarieties of codimension two at all. Therefore we may
assume that $\mathop{\rm codim} B\geq 3$. \vspace{0.1cm}

{\bf Proposition 3.1.} {\it The subvariety $B$ is of codimension}
$\geq 4$.\vspace{0.1cm}

{\bf Proof.} Assume the converse. By Proposition 2.1 then
$\mathop{\rm codim}B=3$. The case $\mathop{\rm deg} B=1$ is
excluded below in Sec. 3.2 and 3.3. Therefore we may assume that
$\mathop{\rm deg} B\geq 2$.

Note that the morphism
$$
\sigma|_B\colon B\to\sigma(B)=\bar{B}
$$
is birational. Indeed, let $Z=(D_1\circ D_2)$ be the
self-intersection of the system $\Sigma$, then
$$
\mathop{\rm mult}\nolimits_B Z>4n^2,
$$
whence it follows that if $\sigma|_B$ is a double cover, then
$$
\mathop{\rm mult}\nolimits_{\bar{B}}\sigma_*Z>8n^2,
$$
however, $\sigma_*Z$ is an effective cycle of codimension two on
${\mathbb P}$ of degree $8n^2$. We get a contradiction.

Therefore, $\mathop{\rm deg}{\bar{B}}=\mathop{\rm deg}_HB\geq 2$
(since the branch divisor $W$ does not contain linear subspaces
of codimension three). Let $p, q\in\bar{B}$ be points of general
position, $L\subset{\mathbb P}$ the line, connecting these
points, $\Pi\supset L$ a generic (two-dimensional) plane,
$\Lambda=\sigma^{-1}(\Pi)$ an irreducible surface on $V$. If
$L\not\subset\mathop{\rm Supp}\sigma_*Z$, then the intersection
$\Pi\cap\mathop{\rm Supp}\sigma_*Z$, and therefore, also the
intersection $\Lambda\cap\mathop{\rm Supp}Z$, is
zero-dimensional, so that we get
$$
8n^2=(\Lambda\cdot Z)\geq \sum\limits_{x\in\sigma^{-1}(L)\cap
B}(\Lambda\cdot Z)_x\geq
$$
$$
\geq \sum_{x\in\sigma^{-1}(L)\cap B}\mathop{\rm
mult}\nolimits_xZ>8n^2,
$$
a contradiction. Therefore, $L\subset\mathop{\rm Supp}\sigma_*Z$.

Let $Q\subset{\mathbb P}$ be the irreducible subvariety, swept
out by all secant lines of the variety $\bar{B}$. By what we have
proved, $\mathop{\rm codim}Q=2$, so that $Q$ is a subspace of
codimension two and $\bar{B}\subset Q$ is some hypersurface.

Now let us write down
$$
Z=a\sigma^{-1}(Q)+Z^{\sharp},
$$
where $Z^{\sharp}$ does not contain the subvariety
$\sigma^{-1}(Q)$ as a component and $a\geq 1$. The cycle $Z$
satisfies the linear inequality
$$
2\mathop{\rm mult}\nolimits_BZ>\mathop{\rm deg}Z.
$$
It is easy to see that any effective cycle of codimension two,
satisfying this inequality, contains the subvariety
$\sigma^{-1}(Q)$ as a component: as above,
$$
\mathop{\rm deg}Z=(\Lambda\cdot Z)\geq\sum_{x\in\sigma^{-1}(L)\cap
B}\mathop{\rm mult}\nolimits_BZ>\mathop{\rm deg}Z
$$
for every secant line $L$ of the variety $\bar{B}$, which is not
contained in the support of the cycle $\sigma_*Z$ (and a generic
plane $\Pi\supset L$). However,
$$
\mathop{\rm mult}\nolimits_B\sigma^{-1}(Q)=1\,\,\mbox{and}
\,\,\mathop{\rm deg}\sigma^{-1}(Q)=2
$$
(recall that for a general hypersurface $W$ the intersection
$Q\cap W$ has at most zero-dimensional singularity, so that
$\sigma^{-1}(Q)$ is an irreducible set), whence it follows that
the cycle $Z^{\sharp}$ satisfies the inequality
$$
2\mathop{\rm mult}\nolimits_BZ^{\sharp}>\mathop{\rm deg}Z^{\sharp}
$$
and therefore contains the subvariety $\sigma^{-1}(Q)$ as a
component. Contradiction.

This excludes the case $\mathop{\rm deg} B\geq 2$. \vspace{0.3cm}


{\bf 3.2. Exclusion of infinitely near singularities with
$\mathop{\rm deg} B=1$.} Starting from this moment and up to the
end of the section we assume that $\mathop{\rm deg} B=1$. By the
conditions of general position this case can realize for the
double spaces of dimension 4 only. Let $X$ be the
$\sigma$-preimage of a generic 3-plane in ${\mathbb P}$ (in
particular, intersecting $\bar B$ at exactly one point). Then
$\sigma_X\colon X\to{\mathbb P}^3$ is a double cover branched
over a smooth hypersurface $W_X\subset{\mathbb P}^3$ of degree
$2m_X\geq 8$, $o=X\cap B$ a point lying outside the ramification
divisor:
$$
p=\sigma_X(o)\not\in W_X,
$$
where $H_X$ is the pull back via $\sigma_X$ of the class of a
plane in ${\mathbb P}^3$. To simplify the notations, we write $H$
instead of $H_X$. By Proposition 0.4 we may assume that on $X$
there are no lines passing through the point $o$, that is, for
any line $L\subset{\mathbb P}^3$, $L\ni p$, the curve
$\sigma^{-1}_X(L)$ is irreducible.

By the symbol $\Sigma_X$ we denote the restriction of the system
$\Sigma$ onto $X$. The movable linear system
$\Sigma_X\subset|2nH|$ has a maximal singularity with the centre
at the point $o$, that is, for the pair $(X,\frac{1}{n}\Sigma_X)$
the point $o$ is a centre of a non canonical singularity. Assume
that the inequality
$$
\nu=\mathop{\rm mult}\nolimits_o\Sigma_X\leq 2n
$$
holds, that is, the point $o$ itself is not maximal (see Lemma
3.2 which is proved below). Let us blow up this point:
$$
\varphi\colon\widetilde{X}\to X,
$$
$E=\varphi^{-1}(o)\cong{\mathbb P}^2$ is the exceptional divisor.
\vspace{0.1cm}

{\bf Proposition 3.2.} {\it The centre of the maximal singularity
on $\widetilde{X}$ is a line in} $E\cong{\mathbb
P}^2$.\vspace{0.1cm}

{\bf Proof.} If the centre of the maximal singularity is a curve
$C\subset E$ of degree $d_C\geq 1$, then the inequality
$$
\nu>nd_C
$$
holds, whence by the assumptions above we get $d_C=1$, that is,
$C$ is a line. Therefore, it is sufficient to exclude the case
when the centre of the singularity is a point $y\in E$. Let us
assume that this is the case and show that this assumption leads
to a contradiction.\vspace{0.1cm}

{\bf Lemma 3.1.} {\it For any irreducible curve $C\subset X$ the
inequality
\begin{equation}\label{c1}
\mathop{\rm mult}\nolimits_oC+\mathop{\rm
mult}\nolimits_y\widetilde{C}\leq\mathop{\rm deg}C=(C\cdot H),
\end{equation}
holds, where $\widetilde{C}\subset\widetilde{X}$ is the strict
transform.} \vspace{0.1cm}

{\bf Proof.} Let
$$
\bar{\varphi}\colon\widetilde{\mathbb P}\to{\mathbb P}^3
$$
be the blow up of the point $p=\sigma_X(o)$,
$\bar{E}=\bar{\varphi}^{-1}(o)$ the exceptional divisor. The
morphism $\sigma_X$ induces an isomorphism
$$
\sigma_E\colon E\to\bar{E}.
$$
Set $\bar{y}=\sigma_E(y)\in\bar{E}$. For any plane $P\ni p$ such,
that its strict transform $\widetilde{P}\subset\widetilde{\mathbb
P}$ contains the point $\bar{y}$, its inverse image
$H=\sigma^{-1}_X(P)$ contains the point $o$ and $\widetilde{H}\ni
y$. Let us denote by the symbol
$$
|H-o-y|
$$
the linear subsystem of the system $H$, defined by that condition.
Obviously,
$$
\mathop{\rm Bs}|H-o-y|=\sigma^{-1}(L),
$$
where $L\ni p$ is the line in ${\mathbb P}^3$ with the tangent
direction $\bar{y}$ at the point $p$. Let $C\subset X$ be an
irreducible curve, $C\ni p$.

Recall that by assumption there are no lines through the point
$o$, that is, the curve $\sigma^{-1}(L)$ is irreducible. We get
$$
\mathop{\rm mult}\nolimits_o\sigma^{-1}(L)=\mathop{\rm
mult}\nolimits_y\widetilde{\sigma^{-1}(L)}=1
$$
and $(H\cdot\sigma^{-1}(L))=2$, so that for the curve
$\sigma^{-1}(L)$ the inequality (\ref{c1}) holds.

Assume that $C\neq \sigma^{-1}(L)$. For a generic divisor
$R\in|H-o-y|$ we have
$$
(C\cdot R)=(C\cdot H)\geq(C\cdot R)_o\geq
$$
$$
\geq\mathop{\rm
mult}\nolimits_oC+(\widetilde{C}\cdot\widetilde{R})_y\geq\mathop{\rm
mult}\nolimits_oC+\mathop{\rm mult}\nolimits_y\widetilde{C},
$$
which is what we need ($\widetilde{R}$ is the strict transform of
$R$ on $\widetilde{X}$). Q.E.D. for the lemma.\vspace{0.1cm}

Now we complete the proof of Proposition 3.2 by word for word the
same arguments as the proof of the $8n^2$-inequality (Lemma 4.2).
Indeed, let
$$
\varphi_{i,i-1}\,\colon X_i\to X_{i-1},
$$
$i=1,\dots,N$, be the resolution of the maximal singularity, that
is, $\varphi_{i,i-1}$ blows up its centre $B_{i-1}$ on $X_{i-1}$,
$E_i=\varphi^{-1}_{i,i-1}(B_{i-1})$ is the exceptional divisor.
For $i=1,\dots,L$ the centres of the blow ups are points, for
$i=L+1,\dots,N$ they are curves, and moreover, it follows from the
inequality $\nu\leq 2n$ that all these curves are smooth and
rational: $B_L\subset E_L\cong{\mathbb P}^2$ is a line,
$B_i\subset E_i$ is a section of the ruled surface $E_i\to
B_{i-1}$ for $i=L+1,\dots,N-1$. By the same inequality $\nu\leq
2n$ we have $N\geq L+1$ and
$$
B_L\not\subset E^L_{L-1},
$$
that is, $L+1\nrightarrow L-1$ in the oriented graph of the
sequence of blow ups $\varphi_{i,i-1}$. Finally, by assumption
$L\geq 2$: more precisely, $B_0=o$ and $B_1=y\in E_1$. Now
repeating the proof of Lemma 4.2 word for word, we get the
inequality
$$
\mathop{\rm mult}\nolimits_oZ+\mathop{\rm
mult}\nolimits_y\widetilde{Z}>8n^2
$$
for the self-intersection $Z=(D_1\circ D_2)$ of the movable linear
system $\Sigma$. However, $Z$ is an effective 1-cycle of degree
$\mathop{\rm deg}Z=(Z\cdot H)=8n^2$. We obtained a contradiction
with Lemma 3.1 which completes the proof of Proposition 3.2.
\vspace{0.3cm}


{\bf 3.3. Exclusion of the last case: preliminary constructions.}
In order to complete the proof of Proposition 3.1, it remains to
exclude the situation described in Proposition 3.2. We assume
that $M\geq 4$.

Let $L\subset{\mathbb P}^4$ be the line generating a line on $V$,
that is, $\sigma^{-1}(L)=C_+\cup C_-$, where $C=C_+$ and $C_-$
are smooth rational curves. Let
$$
\varphi\colon\widetilde{V}\to
V\quad\mbox{and}\quad\varphi_{\mathbb P}\colon\widetilde{\mathbb
P}\to{\mathbb P}^4
$$
be the blow ups of the curve $C$ and the line $L$, respectively,
with the exceptional divisors
$$
E=\varphi^{-1}(C)\subset\widetilde{V}\quad\mbox{and}\quad
E_{\mathbb P}=\varphi^{-1}_{\mathbb
P}(L)\subset\widetilde{\mathbb P}.
$$
The morphism $\sigma$ induces a rational map
$$
\sigma_E\colon E\dashrightarrow E_{\mathbb P},
$$
which is a birational isomorphism, mapping
$E\backslash\varphi^{-1}(C\cap\sigma^{-1}(W))$ isomorphically onto
$E_{\mathbb P}\backslash\sigma^{-1}_{\mathbb P}(L\cap W)$. In
particular, for any irreducible surface $S\subset E$, covering
$C$, its image
$$
\sigma_E(S)\subset E_{\mathbb P}\cong L\times{\mathbb P}^2
$$
is well defined. \vspace{0.1cm}

{\bf Proposition 3.3.} {\it The movable linear system
$\Sigma\subset|2nH|$ can not have a maximal singularity, the
centre of which on $V$ is the curve $C$, and on $\widetilde{V}$
it is some surface} $S\subset E$.\vspace{0.1cm}

{\bf Proof.} Assume the converse: such a maximal singularity
exists.\vspace{0.1cm}

{\bf Lemma 3.2.} {\it The following inequality holds:
$$
\mathop{\rm mult}\nolimits_C\Sigma\leq 2n,
$$
that is, the curve $C$ itself is not a maximal subvariety of the
system $\Sigma$.}\vspace{0.1cm}

{\bf Proof.} Let $P\subset{\mathbb P}^4$ be a generic plane,
containing the line $L$. It is easy to see that the intersection
$P\cap W$ is a non-singular curve, so that $Q=\sigma^{-1}(P)$ is a
non-singular $K3$ surface. The restriction $\Sigma|_Q=\Sigma_Q$ is
a linear system of curves that has, generally speaking, two fixed
components, $C_+$ and $C_-$, of multiplicity $\nu_+$ and $\nu_-$,
respectively. Therefore, on $Q$ the inequalities
$$
((2nH_Q-\nu_+C_+-\nu_-C_-)\cdot C_{\pm})\geq 0
$$
hold, where $H_Q=H|_Q$, or, explicitly,
$$
\begin{array}{c}
2n+2\nu_+-3\nu_-\geq 0,\\
2n-3\nu_++2\nu_-\geq 0\\
\end{array}
$$
(since $(C^2_{\pm})=-2$ and $(C_+\cdot C_-)=3$). Multiplying the
first inequality by 2, the second one by 3 and putting them
together, we obtain
$$
10n-5\nu_+\geq 0,
$$
which is what we need. Q.E.D. for the lemma.\vspace{0.1cm}

{\bf Corollary 3.1.} {\it For a point $x\in L$ we have:
$$
S\cap\sigma^{-1}(x)
$$
is a line in $\sigma^{-1}(x)\cong{\mathbb P}^2$. The graph of the
resolution of the maximal singularity is a chain.}\vspace{0.1cm}

{\bf Proof:} the inequality
$$
\mathop{\rm mult}\nolimits_C\Sigma\geq n\mathop{\rm
deg}(S\cap\sigma^{-1}(x))
$$
holds, which implies the first claim. The second is
obvious.\vspace{0.1cm}

Let us continue our proof of Proposition 3.3.

The surface $\sigma_E(S)$ in $E_{\mathbb P}\cong L\times{\mathbb
P}^2={\mathbb P}^1\times{\mathbb P}^2$ is of bidegree $(d,1)$.
\vspace{0.3cm}


{\bf 3.4. The hard case $d=0$.} Имеет место\vspace{0.1cm}

{\bf Proposition 3.4.} {\it The case $d=0$ is
impossible.}\vspace{0.1cm}

{\bf Proof.} This is the hardest case and to consider it, we have
to inspect quite a few possible cases.

First of all, note, that there exists a unique hyperplane
$\Pi\subset{\mathbb P}^4$, cutting out $S$ on $E_{\mathbb P}$:
$$
S=\widetilde{\Pi}\cap E_{\mathbb P},
$$
where $\widetilde{\Pi}\subset\widetilde{\mathbb P}$ is the strict
transform. Since the linear system $\Sigma$ is movable, its
restriction
$$
\Sigma_{\Pi}=\Sigma|_{\sigma^{-1}(\Pi)}
$$
is a non-empty linear system (possibly with fixed components).
Now let us consider a generic plane $P\supset L$, $P\subset\Pi$,
and argue in exactly the same way as in the proof of Lemma 2.1:
we restrict an effective divisor in the system $\Sigma_{\Pi}$
onto the surface $Q=\sigma^{-1}(P)$ and show that the effective
curve $\Sigma_Q$ obtained in this way cannot contain the curve
$C$ with a multiplicity strictly higher than $2n$.

Unfortunately, we can not argue in word for word the same way as
in the proof of Lemma 2.1, since the surface $Q$, generally
speaking, has singular points. We need to resolve the
singularities, and for that purpose, in its turn, to list the
possible cases for a generic branch divisor $W$. The singularities
can appear because the hyperplane $\Pi$ (which is uniquely
determined by the system $\Sigma$) may turn out to be the tangent
hyperplane to $W$ at one or more points of intersection of the
line $L$ and $W$.

By the symbol $T_xW$ for a point $x\in W$ we denote the hyperplane
in ${\mathbb P}^4$, which is tangent to $W$ at the point $x$.

For the scheme-theoretic intersection $(L\circ W)$ there are three
possible cases:\vspace{0.1cm}

1. $(L\circ W)=2x_1+2x_2+2x_3$, where $x_1,x_2,x_3$ are distinct
points on the line $L$,

2. $(L\circ W)=4x_1+2x_2$, where $x_1\neq x_2$ are two distinct
point,

3. $(L\circ W)=6x$, $x\in L$.\vspace{0.1cm}

The first case takes place for a line of general position (a
three-dimensional family). The case 2 takes place for a
two-dimensional, the case 3 for a one-dimensional family of lines
$C\subset V$.

Furthermore, the hyperplane $\Pi$ is tangent to the divisor $W$ at
the points $x,y$ if and only if
$$
\Pi=T_xW=T_yW.
$$
Therefore we can detailise the cases 1 and 2 in the following
way:\vspace{0.1cm}

1.1. The three hyperplanes $T_{x_i}W$ are distinct (the case of
general position),

1.2. $T_{x_1}W=T_{x_2}W\neq T_{x_3}W$, this case takes place for a
one-dimensional family of lines, since the coincidence of two
tangent hyperplanes imposes two independent conditions on the line
$C$,

1.3. $T_{x_1}W=T_{x_2}W=T_{x_3}W$, this case does not take place
for a general divisor $W$, however we consider it to make the
picture complete,

2.1. $T_{x_1}W\neq T_{x_2}W$ are distinct hyperplanes,

2.2. $T_{x_1}W=T_{x_2}W$, this case takes place for a finite
number of lines.\vspace{0.1cm}

The strategy of the further arguments is as follows. Assuming that
$\Pi$ is tangent to $W$ at at least one point (otherwise we repeat
the proof of Lemma 2.1 without modifications), we resolve
singularities of the surface $Q=\sigma^{-1}(P)$ for a general
plane $P$, where $L\subset P\subset\Pi$. Let $T_i$, $i\in I$, be
the set of irreducible exceptional curves (they are $(-2)$-curves
on a $K3$-surface $\widetilde{Q}$). Let $\widetilde{C_{\pm}}$ be
strict transforms of the curves $C_{\pm}$ on $\widetilde{Q}$,
$\widetilde{\Sigma}_Q$ the strict transform of the linear system
$\Sigma_Q$ on $\widetilde{Q}$. For some no-negative integers
$a_i\in{\mathbb Z}_+$, $i\in I$, we get
$$
\Sigma_Q\subset|2nH_Q-\sum_{i\in I}a_iT_i|,
$$
where $H_Q=H|_Q$ (the pull back of this class on $\widetilde Q$ we
denote by the same symbol). The linear system
$\widetilde{\Sigma}_Q$ contains the curves $C_{\pm}$ as fixed
components of multiplicities $\nu_{\pm}$, whereas at least one of
these two multiplicities by construction is strictly higher than
$2n$. We assume that $\nu_+>2n$. Therefore, the $\sharp I+2$
linear inequalities hold:
$$
((2nH_Q-\nu_+\widetilde{C}_+-\nu_-\widetilde{C}_--\sum_{i\in
I}a_iT_i)\cdot R)\geq 0,
$$
where $R\in\{\widetilde {C}_+,\widetilde{C}_-\}\cup\{T_i|i\in
I\}$. In each of the possible cases this system of linear
inequalities gives the estimate $\nu_+\leq 2n$, contradicting the
initial assumption. This would complete the proof of Proposition
3.4.

Let us realize the program that was described. For that purpose,
we list all possible types of singularities of the surface $Q$ for
a {\it generic} plane $P\subset\Pi$, $P\supset L$ (in the
assumption that the divisor $W$ is generic). For an arbitrary
plane there are many more cases, but we do not need them. The
types of singularities, listed below, are obtained by elementary
computations in the affine coordinates $z_1,z_2,z_3,z_4$ on
${\mathbb P}^4$, in which the hyperplane $\Pi$ is given by the
equation $z_4=0$, and the line $L$ is given by the system of
equations $z_1=z_2=0$, so that the plane $P$ is given by the
equation $z_1+\beta z_2=0$, where $\beta\in{\mathbb C}$ is some
number. Direct coordinate computations show that if singularities
of the surface $Q$ are worse than in the cases listed below, for
at least one line $C\subset V$, then $W$ is not a hypersurface of
general position. The computations are absolutely elementary and
we omit them.

Here is the list of possible types of singularities.\vspace{0.1cm}

{\bf Type A.} One ordinary double point, on $\widetilde{Q}$ there
is one exceptional curve $E$. The multiplication table:
$$
\begin{array}{cccc}
& \widetilde{C}_+ & \widetilde{C}_- & E\\
\widetilde{C}_+ & -2 & 2 & 1\\
\widetilde{C}_- & 2 & -2 & 1\\
E & 1 & 1 & -2\\
\end{array}
$$
This type takes place in the cases 1.1, 2.1 and 3.\vspace{0.1cm}

{\bf Type B.} One degenerate double point, resolved by one blow
up, on $\widetilde{Q}$ there are two exceptional lines $E_+$ and
$E_-$. The multiplication table:
$$
\begin{array}{ccccc}
& \widetilde{C}_+ & \widetilde{C}_- & E_+ &
E_-\\
\widetilde{C}_+ & -2 & 2 & 1 & 0\\
\widetilde{C}_- & 2 & -2 & 0 & 1\\
E_+ & 1 & 0 & -2 & 1\\
E_- & 0 & 1 & 1 & -2\\
\end{array}
$$
This type takes place in the cases 1.1, 2.1.\vspace{0.1cm}

{\bf Type C.} A degenerate double point on $Q$, the exceptional
divisor of its blow up is a pair of lines, the point of their
intersection is an ordinary double point of the surface (resolved
by one blow up). On $\widetilde{Q}$ there are three exceptional
curves $E_+$, $E_-$, $E$ with the multiplication table
$$
\begin{array}{cccccc}
& \widetilde{C}_+ & \widetilde{C}_- & E_+ &
E_- & E\\
\widetilde{C}_+ & -2 & 2 & 1 & 0 & 0\\
\widetilde{C}_- & 2 & -2 & 0 & 1 & 0\\
E_+ & 1 & 0 & -2 & 0 & 1\\
E_- & 0 & 1 & 0 & -2 & 1\\
E & 0 & 0 & 1 & 1 & -2\\
\end{array}
$$
This type takes place in the case 1.1.\vspace{0.1cm}

{\bf Type D.} The surface $Q$ has two ordinary double points
(resolved by one blow up). On $\widetilde {Q}$ there are two
exceptional curves $E_1$ and $E_2$. The multiplication table:
$$
\begin{array}{ccccc}
& \widetilde{C}_+ & \widetilde{C}_- & E_1 &
E_2\\
\widetilde{C}_+ & -2 & 1 & 1 & 1\\
\widetilde{C}_- & 1& -2 & 1 & 1\\
E_1 & 1 & 1 & -2 & 0\\
E_2 & 1 & 1 & 0 & -2\\
\end{array}
$$
This type takes place in the cases 1.2 and 2.2.\vspace{0.1cm}

These four types complete the list of possible singularities of
the surface $Q$ for a variety $V$ of general position. However the
condition of general position is not essential. The author
considered examples of more complicated singularities and our
method in all cases gives a proof of Proposition 3.4. As an
illustration, in addition to the types A-D, we will give two more
examples (they do not take place on a variety of general
position). \vspace{0.1cm}

{\bf Type E.} Two singular points: a non-degenerate one and a
degenerate one. Both are resolved by one blow up. On
$\widetilde{Q}$ there are three exceptional curves: $E$
(corresponds to the non-denerate point) and $E_{\pm}$ (they
correspond to the exceptional lines on the blow up of the
degenerate point). The multiplication table:
$$
\begin{array}{cccccc}
& \widetilde{C}_+ & \widetilde{C}_- & E_+ &
E_- & E\\
\widetilde{C}_+ & -2 & 1 & 1 & 0 & 1\\
\widetilde{C}_- & 1 & -2 & 0 & 1 & 1\\
E_+ & 1 & 0 & -2 & 1 & 0\\
E_- & 0 & 1 & 1 & -2 & 0\\
E & 1 & 1 & 0 & 0 & -2\\
\end{array}
$$
This type takes place on a variety of non-general position in the
case 1.2.\vspace{0.1cm}

{\bf Type F.} On the surface $Q$ there are three non-degenerate
double points. On $\widetilde{Q}$ there are three exceptional
curves $E_1$, $E_2$, $E_3$. The multiplication table:
$$
\begin{array}{cccccc}
& \widetilde{C}_+ & \widetilde{C}_- & E_1 &
E_2 & E_3\\
\widetilde{C}_+ & -2 & 0 & 1 & 1 & 1\\
\widetilde{C}_- & 0 & -2 & 1 & 1 & 1\\
E_1 & 1 & 1 & -2 & 0 & 0\\
E_2 & 1 & 1 & 0 & -2 & 0\\
E_3 & 1 & 1 & 0 & 0 & -2\\
\end{array}
$$
This type takes place on a variety of non-general position in the
case 1.3. \vspace{0.1cm}

It remains to realize the program described above for each type of
singularities. We will consider two examples, A and C, in the
other cases the computations are similar. After that we explain
the essence of the computations.

Consider the type A. Let $\nu_+$, $\nu_-$ and $\alpha$ be the
multiplicities of the curves $\widetilde{C}_+$, $\widetilde{C}_-$
and $E$ in the linear system $\Sigma_Q$, pulled back on
$\widetilde{Q}$. Multiplying the class
$$
2nH_Q-\nu_+\widetilde{C}_+-\nu_-\widetilde{C}_--\alpha E
$$
by $\widetilde{C}_+$, $\widetilde{C}_-$ and $E$, we obtain a
system of linear inequalities:
$$
\begin{array}{cr}
2n & +2\nu_+-2\nu_--\alpha\geq 0,\\
2n & -2\nu_++2\nu_--\alpha\geq 0,\\
& -\nu_+-\nu_-+2\alpha\geq 0.
\end{array}
$$
Adding to the first and second inequalities one half of the third
one, we get the system
$$
\begin{array}{c}
4n+3\nu_+-5\nu_-\geq 0,\\
4n-5\nu_++3\nu_-\geq 0.
\end{array}
$$
Multiplying the first inequality by 3, the second one by 5 and
putting them together, we get
$$
32n-16\nu_+\geq 0,
$$
which is precisely what we need.

Let us consider the type C. Denoting the multiplicities of the
components $\widetilde{C}_+$, $\widetilde{C}_-$, $E_+$, $E_-$, $E$
by the symbols $\nu_+$, $\nu_-$, $\alpha_+$, $\alpha_-$, $\alpha$,
respectively, multiply the effective class
$$
2nH_Q-\nu_+\widetilde{C}_+-\nu_-\widetilde{C}_-
-\alpha_+E_+-\alpha_-E_--\alpha E
$$
by $\widetilde{C}_+$, $\widetilde{C}_-$, $E_+$, $E_-$, $E$ and
obtain the system of inequalities
$$
\begin{array}{ccccccc}
2n & +2\nu_+ & -2\nu_- & -\alpha_+ & & & \geq 0\\
2n & -2\nu_+ & +2\nu_- &  & -\alpha_- & & \geq 0\\
&  -\nu_+ & & +2\alpha_+ &   & -\alpha &   \geq 0\\
& &  -\nu_- & & +2\alpha_- &  -\alpha & \geq 0\\
& & & -\alpha_+ & -\alpha_- & +2\alpha & \geq 0.
\end{array}
$$
By means of the fifth inequality eliminate $\alpha$ in the third
and fourth inequalities, which take the form
$$
\begin{array}{ccccc}
-\nu_+ &  & +\frac32\alpha_+ & -\frac12\alpha_- & \geq 0,\\
 & -\nu_- & -\frac12\alpha_+ & +\frac32\alpha_- &\geq 0.\\
\end{array}
$$
Multiplying one of these inequalities by $\frac32$, another one by
$\frac12$ and putting them together, we obtain the inequalities
$$
\begin{array}{c}
-3\nu_+-\nu_-+4\alpha_+\geq 0,\\
-\nu_+-3\nu_-+4\alpha_-\geq 0,\\
\end{array}
$$
which make it possible to eliminate $\alpha_+$, $\alpha_-$ in the
first two inequalities and obtain the system of inequalities
$$
\begin{array}{c}
8n+5\nu_+-9\nu_-\geq 0,\\
8n-9\nu_++5\nu_-\geq 0,\\
\end{array}
$$
whence, similar to the case A, we get
$$
112n-56\nu_+\geq 0,
$$
which is precisely what we need.

The other types of singularities are considered in a similar way.
Now let us explain, why all the types listed above lead to the
inequality $2n\geq\nu_+$. Let us consider the space
$$
{\cal L}={\mathbb R}[\widetilde{C}_+]\oplus{\mathbb
R}[\widetilde{C}_-]\oplus{\cal E},
$$
where
$$
{\cal E}=\bigoplus^k_{i=1}T_i,
$$
where $\{T_i|i=1,\dots,k\}$ is the set of irreducible exceptional
curves on $\widetilde{Q}$. On ${\cal L}$ there is a natural
bilinear form, generated by intersection of curves. Set
$$
\Theta=\|(T_i\cdot T_j)\|_{1\leq i,j\leq k}
$$
to be the negative definite matrix of the intersection form on
${\cal E}$. Set $\Theta^{-1}=\|\lambda_{ij}\|$ to be the inverse
matrix. It is easy to check that in each of the cases A-F (and in
all other cases of singularities of non-general position, studied
by the author) the matrix $\Theta$ satisfies the following
condition:\vspace{0.1cm}

all coefficients $\lambda_{i,j}$, $1\leq i,j\leq k$, of the
inverse matrix ${\Theta}^{-1}$ are negative.\vspace{0.1cm}

Let ${\cal E}_+=\{a_iT_i\,|\, a_i\in{\mathbb R}_+\}$ be the
positive coordinate cone. Since the matrix $\Theta$ is
non-degenerate, there exist uniquely determined vectors
$e_{\pm}\in{\cal E}$, such that
$$
R_{\pm}=\widetilde{C}_{\pm}+e_{\pm}\in{\cal E}^{\bot}.
$$
\vspace{0.1cm}

{\bf Lemma 3.3.} $e_{\pm}\in{\cal E}_+$.\vspace{0.1cm}

{\bf Proof.} This follows immediately from the inequalities
$(C_{\pm}\cdot T_i)\geq 0$ and the properties of the matrix
$\Theta$ ($\lambda_{ij}< 0$).

Since the intersection form on ${\cal L}$ is non-degenerate, the
class $H_Q$ can be considered as an element of the space ${\cal
L}$, $H_Q\in{\cal E}^\bot$.\vspace{0.1cm}

{\bf Lemma 3.4.} {\it The following inequality holds:}
$H_Q=R_++R_-$.\vspace{0.1cm}

{\bf Proof.} $H_Q$ is the class of a Cartier divisor on $Q$, which
consists of the curves $C_+$ and $C_-$. Therefore, for some
$e\in{\cal E}$ we get the equality
$$
H_Q=\widetilde{C}_++\widetilde{C}_-+e.
$$
Since $H_Q\in{\cal E}^\bot$, the fact that the quadratic
intersection form on ${\cal L}$ is non-degenerate, implies the
claim of the lemma.\vspace{0.1cm}

Obviously, $(R_{\pm}\cdot H_Q)=1$.

It is easy to check that $(R^2_{\pm})=-a<0$, so that the
intersection form on the two-dimensional space
$$
{\cal R}={\mathbb R}[R_+]\oplus{\mathbb R}[R_-]
$$
is given by the matrix
$$
\left(\begin{array}{cc}
-a & 1+a\\
1+a & -a\\
\end{array}\right),
$$
the inverse matrix for which is
$$
\frac{1}{1+2a}\left(\begin{array}{cc}
a & 1+a\\
1+a & a\\
\end{array}\right).
$$
By what was said above, non-negativity of the intersections of the
class
$$
2nH_Q-\nu_+\widetilde{C}_+-\nu_-\widetilde{C}_-
-\sum^k_{i=1}b_iT_i
$$
(where $b_i\in {\mathbb Z}_+$) with the classes $\widetilde{C}_+$,
$\widetilde{C}_-$, $T_i$ implies non-negativity of the
intersections of the class
$$
\beta_+R_++\beta_-R_-=2nH_Q-\nu_+R_+-\nu_-R_-
$$
with the classes $R_+$ and $R_-$, which implies that
$\beta_{\pm}\geq 0$. However,
$$
\beta_{\pm}=2n-\nu_{\pm},
$$
which is what we need.

This general argument works for any type of singularities of the
surface $Q$, satisfying the two properties: the elements
$\lambda_{ij}$ of the matrix $\Theta^{-1}$ are all negative and
$(R^2_{\pm})=-a< 0$. These properties are to be checked directly
(for instance, for the type C the matrix $\Theta^{-1}$ is
$$
\left(\begin{array}{ccc}
-\frac34 & -\frac14 & -\frac12\\
-\frac14 & -\frac34 & -\frac12\\
-\frac12 & -\frac12 & -1\\
\end{array}\right),
$$
as we need). It seems, however, that this is a consequence of some
general fact.

Proof of Proposition 3.4 is complete.

Let us get back to the proof of Proposition 3.3. \vspace{0.3cm}


{\bf 3.5. The case $d\geq 1$: end of the proof of Proposition
3.3.} We assume that $d\geq 1$. Let
$$
\varphi_{i,i-1}\colon V_i\to V_{i-1},
$$
$i=1,\dots,k+1$, be the resolution of the maximal singularity.
Here $V_0=V$,
$$
\varphi_{1,0}\colon V_1\to V
$$
is the blow up of the curve $C$, that is, $\widetilde{V}\cong
V_1$, and $B_1\subset E_1=E$ is the surface $S$. Setting, as
usual,
$$
\nu_i=\mathop{\rm mult}\nolimits_{B_{i-1}}\Sigma^{i-1},
$$
write down the Noether-Fano inequality:
$$
\nu_1+\dots+\nu_{k+1}>(k+2)n.
$$
Consider a generic plane $P\supset L$ and a non-singular surface
$Q=\sigma^{-1}(P)$. By the inequality $d\geq 1$ the section
$\widetilde{P}\cap E_{\mathbb P}$ intersects transversally the
surface $\sigma_E(S)$ at at least one point of general position,
and the surface $\widetilde{P}$ intersects transversally
$\sigma_E(S)$ at a point of general position. Therefore, the
surface $Q^1=\widetilde{Q}$ intersects transversally the surface
$S$ at at least one point of general position, say
$$
x\in S\cap Q^1,
$$
and we may assume that  $x\not\in C^1_-$.

Furthermore, by generality of $P$ the restriction
$\Sigma^1|_{Q^1}$ can have only one fixed component, the curve
$C^1_-$. Set
$$
\nu_-=\mathop{\rm mult}\nolimits_{C^1_-}(\Sigma^1|_{Q^1})=
\mathop{\rm mult}\nolimits_{C_-}(\Sigma|_Q)=\mathop{\rm
mult}\nolimits_{C_-}\Sigma.
$$
Since the intersection $S\cap Q^1$ is transversal, the surfaces
$B_2,\dots,B_k$ generate infinitely near base points
$$
x_i\in B_i\cap Q^i
$$
of the linear system of curves $\Sigma^1|_Q$ on the non-singular
surface $Q^1=Q$, lying over the point $x=x_1$. Since the point $x$
lies outside the fixed component $C_-$, the self-intersection of
the movable part of the linear system $\Sigma^1|_Q$ is not less
than
$$
\nu^2_2+\dots+\nu^2_{k+1}.
$$
This gives the inequality
$$
f(\nu_-,\nu_1,\dots,\nu_{k+1})= (2nH_Q-\nu_1C_+-\nu_-C_-)^2-\sum^{
k+1}_{i=2}\nu^2_i=
$$
$$
=8n^2-4n\nu_1-4n\nu_--2\nu^2_1-2\nu^2_-+
6\nu_1\nu_--\sum^{k+1}_{i=2}\nu^2_i\geq 0
$$
(see the proof of Lemma 2.1 for the intersection numbers).
Besides, as we have shown in the proof of Lemma 3.2, the
inequality
$$
\nu_-\geq\frac{3\nu_1-2n}{2}>\frac{n}{2}
$$
holds. Now let us estimate the function
$f(\nu_-,\nu_1,\dots,\nu_{k+1})$ from above. Let us replace the
Noether-Fano inequality by the equality
\begin{equation}\label{c2}
\nu_1+\nu_2+\dots+\nu_{k+1}=(k+2)n,
\end{equation}
which can only increase the value $f(\cdot)$. Furthermore, let us
fix $\nu_-$ and $\nu_1$ and consider $f$ as a function of
$\nu_2,\dots,\nu_{k+1}$ under the constraint (\ref{c2}).
Obviously its maximum is attained at
$$
\nu_2=\dots=\nu_{k+1}=\frac{(k+2)n-\nu_1}{k}.
$$
On the other hand, the maximum of $f$ as a function of one
argument $\nu_-$ is attained, as it is easy to check, at
$$
\nu_-=\frac{3\nu_1-2n}{2}.
$$
Substituting these values of $\nu_2,\dots,\nu_{k+1},\nu_-$ into
$f(\cdot)$, we obtain the following expression:
$$
\frac{1}{k}[-(k^2+4)n^2+(4-2k)n\nu_1].
$$
For $k\geq 2$ it is obviously negative. Let $k=1$, then we have
$$
-5n^2+2n\nu_1.
$$
We have shown above that $\nu_1\leq 2n$. Therefore, for any
$k\geq 1$ we obtain the estimate
$$
f(\nu_-,\nu_1,\dots,\nu_{k+1})<-kn^2.
$$
We get a contradiction which completes the proof of Proposition
3.3.

Q.E.D. for Proposition 3.1.\vspace{0.3cm}


\section{A local inequality for the self-intersection\\ of a
movable system}

We give a proof of the so called $8n^2$-inequality for the
self-intersection of a movable linear system, correcting the
mistake in the papers [34-36]. The notations of this section are
independent of the rest of the paper.\vspace{0.3cm}

{\bf 4.1. Set up of the problem and start of the proof.} Let $o\in
X$ be a germ of a smooth variety of dimension $\mathop{\rm dim}
X\geq 4$. Let $\Sigma$ be a movable linear system on $X$, and the
effective cycle
$$
Z=(D_1\circ D_2),
$$
where $D_1,D_2\in \Sigma$ are generic divisors, its
self-intersection. Blow up the point $o$:
$$
\varphi\colon X^+\to X,
$$
$E=\varphi^{-1}(o)\cong{\mathbb P}^{\mathop{\rm dim}X-1}$ is the
exceptional divisor. The strict transform of the system $\Sigma$
and the cycle $Z$ on $X^+$ we denote by the symbols $\Sigma^+$ and
$Z^+$, respectively.\vspace{0.1cm}

{\bf Proposition 4.1 ($8n^2$-inequality).} {\it Assume that the
pair
$$
(X,\frac{1}{n}\Sigma)
$$
is not canonical, but canonical outside the point $o$, where $n$
is some positive number. There exists a linear subspace $P\subset
E$ of codimension two (with respect to $E$), such that the
inequality
$$
\mathop{\rm mult}\nolimits_oZ+\mathop{\rm mult}\nolimits_PZ^+>8n^2
$$
holds.} \vspace{0.1cm}

An equivalent claim, but formulated in a rather cumbersome way,
was several times published by Cheltsov [34-36], however his proof
is essentially faulty (see [37]).

{\bf Proof.} The first part of our arguments follows [34,36]. Note
that if $\mathop{\rm mult}_oZ>8n^2$, then for $P$ we may take any
subspace of codimension two in $E$. However, if $\mathop{\rm
mult}_oZ\leq 8n^2$, then the subspace $P$ is uniquely determined:
it follows easily from the connectedness principle of Shokurov and
Koll\' ar [29,38].

Restricting $\Sigma$ onto a germ of a generic smooth subvariety,
containing the point $o$, we may assume that $\mathop{\rm
dim}X=4$. Moreover, we may assume that $\nu=\mathop{\rm
mult}_o\Sigma\leq 2\sqrt{2}n< 3n$, since otherwise
$$
\mathop{\rm mult}\nolimits_oZ\geq\nu^2>8n^2
$$
and there is nothing to prove.\vspace{0.1cm}

{\bf Lemma 4.1.} {\it The pair
\begin{equation}\label{d1}
(X^+,\frac{1}{n}\Sigma^++\frac{(\nu-2n)}{n}E)
\end{equation}
is not log canonical, and the centre of any of its non log
canonical singularities is contained in the exceptional divisor
$E$.}\vspace{0.1cm}

{\bf Proof.} Let $\lambda\colon\widetilde{X} \to X$ be a
resolution of singularities of the pair $(X,\frac{1}{n}\Sigma)$
and $E^*\subset \widetilde{X}$ a prime exceptional divisor,
realizing a non-canonical singularity of that pair. Then
$\lambda(E^*)=o$ and the Noether-Fano inequality holds:
$$
\nu_{E^*}(\Sigma)>na(E^*).
$$
For a generic divisor $D\in\Sigma$ we get $\varphi^*D=D^++\nu E$,
so that
$$
\nu_{E^*}(\Sigma)=\nu_{E^*}(\Sigma^+)+\nu\cdot\nu_{E^*}(E)
$$
and
$$
a(E^*,X)=a(E^*,X^+)+3\nu_{E^*}(E).
$$
From here we get
$$
\nu_{E^*}\left(\frac{1}{n}\Sigma^++\frac{\nu-2n}{n}E\right)=
\nu_{E^*}\left(\frac{1}{n}\Sigma\right)-2\nu_{E^*}(E)>
$$
$$
>a(E^*,X^+)+\nu_{E^*}(E)\geq a(E^*,X^+)+1,
$$
which proves the lemma.

Let $R\ni o$ be a generic three-dimensional germ, $R^+\subset X^+$
its strict transform on the blow up of the point $o$. For a small
$\varepsilon>0$ the pair
$$
\left(X^+,\frac{1}{1+\varepsilon}\frac{1}{n}\Sigma^+
+\frac{\nu-2n}{n}E+R^+\right)
$$
still satisfies the connectedness principle (with respect to the
morphism $\varphi\colon X^+\to X$), so that the set of centres of
non log canonical singularities of that pair is connected. Since
$R^+$ is a non log canonical singularity itself, we obtain, that
there is a non log canonical singularity of the pair (\ref{d1}),
the centre of which on $X^+$ is of positive dimension, since it
intersects $R^+$.

Let $Y\subset E$ be a centre of a non log canonical singularity of
the pair (\ref{d1}) that has the maximal dimension.

If $\mathop{\rm dim}Y=2$, then consider a generic two-dimensional
germ $S$, intersecting $Y$ transversally at a point of general
position. The restriction of the pair (\ref{d1}) onto $S$ is not
log canonical at that point, so that, applying Proposition 4.2,
which is proven below, we see that
$$
\mathop{\rm mult}\nolimits_Y(D^+_1\circ
D^+_2)>4\left(3-\frac{\nu}{n}\right)n^2,
$$
so that
$$
\mathop{\rm mult}\nolimits_oZ\geq\nu^2+\mathop{\rm
mult}\nolimits_Y(D^+_1\circ D^+_2)\mathop{\rm deg}Y>
$$
$$
>(\nu-2n)^2+8n^2,
$$
which is what we need.

If $\mathop{\rm dim}Y=1$, then, since the pair
\begin{equation}\label{d2}
\left(R^+,\frac{1}{1+\varepsilon}\frac{1}{n}\Sigma^+_R
+\frac{\nu-2n}{n}E_R\right),
\end{equation}
where $\Sigma^+_R=\Sigma^+|_{R^+}$ and $E_R=E|_{R^+}$, satisfies
the condition of the connectedness principle and $R^+$ intersects
$Y$ at $\mathop{\rm deg}Y$ distinct points, we conclude that
$Y\subset E$ is a line in ${\mathbb P}^3$.

Now we need to distinguish between the following two cases: when
$\nu\geq 2n$ and when $\nu<2n$. The methods of proving the
$8n^2$-inequality in these two cases are absolutely different.
Consider first the case $\nu\geq 2n$.

Let us choose as $R\ni o$ a generic three-dimensional germ,
satisfying the condition $R^+\supset Y$. Since the pair (\ref{d1})
is effective (recall that $\nu\geq 2n$), one may apply inversion
of adjunction [29, Chapter 17] and conclude that the pair
(\ref{d2}) is not log canonical at $Y$.

Applying now to the pair (\ref{d2}) (where $R^+\supset Y$)
Proposition 4.2 in the same way as it was done for $\mathop{\rm
dim}Y=2$, we obtain the inequality
$$
\mathop{\rm mult}\nolimits_Y(D^+_1|_{R^+}\circ
D^+_2|_{R^+})>4\left(3-\frac{\nu}{n}\right)n^2.
$$
On the left in brackets we have the self-intersection of the
movable system $\Sigma^+_R$, which breaks into two natural
components:
$$
(D^+_1|_{R^+} \circ D^+_2|_{R^+})=Z^+_R+Z^{(1)}_R,
$$
where $Z^+_R$ is the strict transform of the cycle $Z_R=Z|_R$ on
$R^+$ and the support of the cycle $Z^{(1)}_R$ is contained in
$E_R$. The line $Y$ is a component of the effective 1-cycle
$Z^{(1)}_R$.

On the other hand, for the self-intersection of the movable linear
system $\Sigma^+$ we get
$$
(D^+_1\circ D^+_2)=Z^++Z_1,
$$
where the support of the cycle $Z_1$ is contained in $E$. From the
genericity of $R$ it follows that outside the line $Y$ the cycles
$Z^{(1)}_R$ and $Z_1|_{R^+}$ coincide, whereas for $Y$ we get the
equality
$$
\mathop{\rm mult}\nolimits_YZ^{(1)}_R=\mathop{\rm
mult}\nolimits_YZ^++\mathop{\rm mult}\nolimits_YZ_1.
$$
However, $\mathop{\rm mult}_YZ_1\leq\mathop{\rm deg}Z_1$, so that
$$
\mathop{\rm mult}\nolimits_oZ+\mathop{\rm mult}\nolimits_YZ^+=
$$
$$
=\nu^2+\mathop{\rm deg}Z_1+\mathop{\rm mult}\nolimits_YZ^+\geq
$$
$$
\geq\nu^2+\mathop{\rm mult}\nolimits_YZ^{(1)}_R>8n^2,
$$
which is what we need. This completes the case $\nu\geq 2n$.

Note that the key point in this argument is that the pair
(\ref{d1}) is effective. For $\nu<2n$  inversion of adjunction can
not be applied (as it was done in [39]). The additional arguments
in [34-36], proving inversion of adjunction specially for this
pair for $\nu<2n$, are faulty, see [37].\vspace{0.3cm}


{\bf 4.2. The case $\nu<2n$.} Consider again the pair (\ref{d2})
for a generic germ $R\ni o$. Let $y=Y\cap R^+$ be the point of
(transversal) intersection of the line $Y$ and the variety $R^+$.
Since $a(E_R,R)=2$, the non log canonicity of the pair (\ref{d2})
at the point $y$ implies the non log canonicity of the pair
$$
\left(R,\frac{1}{n}\Sigma_R\right)
$$
at the point $o$, whereas the centre of some non log canonical
(that is, log maximal) singularity on $R^+$ is a point $y$.

Now the $8n^2$-inequality comes from the following
fact.\vspace{0.1cm}

{\bf Lemma 4.2.} {\it The following inequality holds:
$$
\mathop{\rm mult}\nolimits_oZ_R+\mathop{\rm
mult}\nolimits_yZ^+_R>8n^2,
$$
where $Z_R$ is the self-intersection of a movable linear system
$\Sigma_R$ and $Z^+_R$ is its strict transform on
$R^+$.}\vspace{0.1cm}

{\bf Proof.} Consider the resolution of the maximal singularity of
the system $\Sigma_R$, the centre of which on $R^+$ is the point
$y$:
$$
\begin{array}{ccc}
R_i & \stackrel{\psi_i}{\to} & R_{i-1}\\
\cup & & \cup\\
E_i & & B_{i-1},\\
\end{array}
$$
where $B_{i-1}$ is the centre of the singularity on $R_{i-1}$,
$R_0=R$, $R_1=R^+$, $E_i=\psi^{-1}_i(B_{i-1})$ is the exceptional
divisor, $B_0=o$, $B_1=y\in E_1$, $i=1,\dots,N$, where the first
$L$ blow ups correspond to points, for $i\geq L+1$ curves are
blown up. Since
$$
\mathop{\rm mult}\nolimits_o\Sigma_R=\mathop{\rm
mult}\nolimits_o\Sigma<2n,
$$
we get $L<N$, $B_L\subset E_L\cong{\mathbb P}^2$ is a line and for
$i\geq L+1$
$$
\mathop{\rm deg}[\psi_i|_{B_i}\colon B_i\to B_{i-1}]=1,
$$
that is, $B_i\subset E_i$ is a section of the ruled surface $E_i$.

Consider the graph of the sequence of blow ups
$\psi_i$.\vspace{0.1cm}

{\bf Lemma 4.3.} {\it The vertices $L+1$ and $L-1$ are not
connected by an arrow:}
$$
L+1\nrightarrow L-1.
$$
\vspace{0.1cm}

{\bf Proof.} Assume the converse: $L+1\to L-1$. This means that
$$
B_L=E_L\cap E^L_{L-1}
$$
is the exceptional line on the surface $E^L_{L-1}$ and the map
$$
E^{L+1}_{L-1}\to E^L_{L-1}
$$
is an isomorphism. As usual, set
$$
\nu_i=\mathop{\rm mult}\nolimits_{B_{i-1}}\Sigma^{i-1}_R,
$$
$i=1,\dots,N$. Let us restrict the movable linear system
$\Sigma^{L+1}_R$ onto the surface $E^{L+1}_{L-1}$ (that is, onto
the plane $E_{L-1}\cong{\mathbb P}^2$ with the blown up point
$B_{L-1}$). We obtain a non-empty (but, of course, not necessarily
movable) linear system, which is a subsystem of the complete
linear system
$$
\left|\nu_{L-1}(-E_{L-1}|_{E_{L-1}})-(\nu_L+\nu_{L+1})B_L\right|.
$$
Since $(-E_{L-1}|_{E_{L-1}})$ is the class of a line on the plane
$E_{L-1}$, this implies that
$$
\nu_{L-1}\geq\nu_L+\nu_{L+1}>2n,
$$
so that the more so $\nu_1=\nu>2n$. A contradiction. Q.E.D. for
the lemma.\vspace{0.1cm}

Set, as usual,
$$
m_i=\mathop{\rm mult}\nolimits_{B_{i-1}}(Z_R)^{i-1},
$$
$i=1,\dots,L$, so that, in particular,
$$
m_1=\mathop{\rm mult}\nolimits_oZ_R\quad\mbox{and}\quad
m_2=\mathop{\rm mult}\nolimits_yZ^+_R.
$$
Let $p_i\geq 1$ be the number of paths in the graph of the
sequence of blow ups $\psi_i$ from the vertex $N$ to the vertex
$i$, and $p_N=1$ by definition. By what we proved,
$$
p_N=p_{N-1}=\dots=p_L=p_{L-1}=1,
$$
and the number of paths $p_i$ for $i\leq L$ is the number of paths
from the vertex $L$ to the vertex $i$. By the technique of
counting multiplicities [1,40], we get the inequality
$$
\sum^L_{i=1}p_im_i\geq\sum^N_{i=1}p_i\nu^2_i
$$
and, besides, the Noether-Fano inequality holds:
$$
\sum^N_{i=1}p_i\nu_i>
n\left(2\sum^L_{i=1}p_i+\sum^N_{i=L+1}p_i\right).
$$
(In fact, a somewhat stronger inequality holds, the {\it log}
Noether-Fano inequality, but we do not need that.) From the last
two estimates one obtains in the standard way [1,40] the
inequality
$$
\sum^L_{i=1}p_im_i>\frac{(2\Sigma_0+\Sigma_1)^2}
{\Sigma_0+\Sigma_1}n^2,
$$
where $\Sigma_0=\sum\limits^L\limits_{i=1}p_i$ and
$\Sigma_1=\sum\limits^N\limits_{i=L+1}p_i=N-L$. Taking into
account that for $i\geq 2$ we get
$$
m_i\leq m_2
$$
and the obvious inequality
$(2\Sigma_0+\Sigma_1)^2>4\Sigma_0(\Sigma_0+\Sigma_1)$, we obtain
the following estimate
$$
p_1m_1+(\Sigma_0-p_1)m_2>4n^2\Sigma_0.
$$
Now assume that the claim of the lemma is false:
$$
m_1+m_2\leq 8n^2.
$$
\vspace{0.1cm}

{\bf Lemma 4.4.} {\it The following inequality holds:}
$\Sigma_0\geq 2p_1$.\vspace{0.1cm}

{\bf Proof.} By definition,
$$
p_1=\sum_{i\to 1}p_i,
$$
however, by Lemma 4.3 from $i\to 1$ it follows that $i\leq L$, so
that $p_1\leq\Sigma_0-p_1$, which is what we need. Q.E.D. for the
lemma.\vspace{0.1cm}

Now, taking into account that $m_2\leq m_1$, we obtain
$$
p_1m_1+(\Sigma_0-p_1)m_2=p_1(m_1+m_2)+(\Sigma_0-2p_1)m_2\leq
$$
$$
\leq 8p_1n^2+(\Sigma_0-2p_1)\cdot 4n^2=4n^2\Sigma_0.
$$
This is a contradiction. Q.E.D. for Lemma 4.2.\vspace{0.1cm}

Proof of Proposition 4.1 is complete.\vspace{0.3cm}

{\bf Remark 4.1.} As it follows from the technique of counting
multiplicities, the graph of the sequence of blow ups $\{\psi_i\}$
can be modified in such a way that all applications still hold,
namely, one can delete all the arrows going from the vertices
$$
L+1,\dots, N
$$
of the upper part of the graph to the vertices
$$
1,\dots, L-1
$$
of the lower part (and both the Noether-Fano inequality and the
estimate for the multiplicities of the self-intersection of the
linear system are intact). The graph, modified in this way,
satisfies the property of Lemma 4.3, which makes it possible to
complete the proof of Lemma 4.2, not using Lemma 4.3 at all.
\vspace{0.3cm}


{\bf 4.3. A local inequality for a surface.} Let $o\in X$ be a
germ of a smooth surface, $C\ni o$ a smooth curve and $\Sigma$ a
movable linear system on $X$. Let, furthermore, $Z=(D_1\circ D_2)$
be the self-intersection of the linear system $\Sigma$, that is,
an effective 0-cycle. Since the situation is local, we may assume
that the support of the cycle $Z$ is the point $o$, that is,
$$
\mathop{\rm deg}Z=(D_1\cdot D_2)_o.
$$
\vspace{0.1cm}

{\bf Proposition 4.2.} {\it Assume that for some real number $a<1$
the pair
\begin{equation}\label{d3}
\left(X,\frac{1}{n}\Sigma+aC\right)
\end{equation}
is not log canonical (that is, for a general divisor $D\in\Sigma$
the pair $(X,\frac{1}{n}D+aC))$ is not log canonical, where $n>0$
is a positive number. Then the estimate holds}
\begin{equation}\label{d4}
\mathop{\rm deg}Z>4(1-a)n^2.
\end{equation}
\vspace{0.1cm}

{\bf Proof.} The original argument see in [31]. We will show that
the inequality (\ref{d4}) follows directly from some well known
facts on the infinitely near singularities of curves on a
non-singular surface [32,33]. Assume that the sequence of blow ups
$$
\varphi_{i,i-1}\,\colon X_i\to X_{i-1},
$$
$i=1,\dots,N$, where $X_0=X$, resolves the non log canonical
singularity of the pair (\ref{d3}). We use the standard notations
and conventions: the centre of the blow up $\varphi_{i,i-1}$ is
the point $x_{i-1}\in X_{i-1}$, its exceptional line is
$$
E_i=\varphi^{-1}_{i,i-1}(x_{i-1})\subset X_i,
$$
the first point to be blown up is $o=x_0$, the blown up points
$x_i$ lie over each other: $x_i\in E_i$. The last exceptional line
$E_N$ realizes the non log canonical singularity of the pair
(\ref{d3}), that is the log Noether-Fano inequality holds:
\begin{equation}\label{d5}
\sum^N_{i=1}\nu_ip_i+an\sum_{x_{i-1}\in
C^{i-1}}p_i>n\left(\sum^N_{i=1}p_i+1\right),
\end{equation}
where $\nu_i=\mathop{\rm mult}_{x_{i-1}}\Sigma^{i-1}$, the symbols
$\Sigma^i$ and $C^i$ stand for the strict transforms on $X_i$ and
$p_i$ is the number of paths in the graph $\Gamma$ of the
constructed sequence of blow ups from the vertex $E_N$ to the
vertex $E_i$, see [1,40]. Assume that
$$
x_{i-1}\in C^{i-1}
$$
for $i=1,\dots,k\leq N$, then the inequality (\ref{d5}) takes the
form
\begin{equation}\label{d6}
\sum^N_{i=1}\nu_ip_i>n\left(\sum^k_{i=1}(1-a)p_i
+\sum^N_{i=k+1}p_i+1\right).
\end{equation}

{\bf Lemma 4.5.} {\it The following inequality holds:}
\begin{equation}\label{d7}
\mathop{\rm deg}Z\geq\sum^N_{i=1}\nu^2_i.
\end{equation}

{\bf Proof:} this is obvious.\vspace{0.1cm}

{\bf Lemma 4.6.} {\it For each $i\in\{1,\dots,N-1\}$ the estimate
\begin{equation}\label{d8}
\nu_i\geq\sum_{j\to i}\nu_j
\end{equation}
holds.}\vspace{0.1cm}

{\bf Proof.} This is a very well known property of multiplicities
of curves at infinitely near points on a non-singular
surface.\vspace{0.1cm}

{\bf Lemma 4.7.} {\it The following estimate is true:
$$
\sum^N_{i=1}\nu^2_i>\frac{\Delta^2}{q}n^2,
$$
where
$$
\Delta=1+(1-a)\sum\limits^k\limits_{i=1}p_i+
\sum\limits^N\limits_{i=k+1}p_i
$$
and $q=\sum\limits^N\limits_{i=1}p^2_i$ (so that $n\Delta$ is the
right-hand side of the inequality (\ref{d6}))}.\vspace{0.1cm}

{\bf Proof.} The minimum of the quadratic form in the right-hand
side of the inequality (\ref{d7}) under the restrictions
(\ref{d8}) and
\begin{equation}\label{d9}
\sum^N_{i=1}\nu_ip_i=\Delta n
\end{equation}
is attained at $\nu_i=p_i\theta$, where $\theta=\frac{\Delta
n}{q}$ is computed from (\ref{d9}). Q.E.D. for the
lemma.\vspace{0.1cm}

Now the claim of Proposition 4.2 follows from a purely
combinatorial fact about the graph $\Gamma$, which we will now
prove.\vspace{0.1cm}

{\bf Lemma 4.8.} {\it Assume that the starting segment of the
graph $\Gamma$ with the vertices $1,\dots,k$ is a chain. Then the
estimate
\begin{equation}\label{d10}
\Delta^2\geq 4(1-a)q
\end{equation}
holds.}\vspace{0.1cm}

{\bf Proof} will be given by induction on the number $N$ of
vertices of the graph $\Gamma$. If $N=1$, then the inequality
(\ref{d10}) holds in a trivial way:
$$
(2-a)^2\geq 4(1-a).
$$
Consider the inequality (\ref{d10}) as a claim on the
non-negativity of a quadratic function of the argument $a$:
$$
a^2\left(\sum^k_{i=1}p_i\right)^2+2a\left(2q-\left(\sum^k_{i=1}p_i\right)
\left(\sum^N_{i=1}p_i+1\right)\right)+
\left(\left(\sum^N_{i=1}p_i+1\right)^2-4q\right) \geq 0
$$
on the interval $a\leq 1$. Since for $a\to\pm\infty$ this function
is positive, it is sufficient to check that its minimum is
non-negative. Elementary computations show that, up to an
inessential positive factor, this minimum is given by the formula
\begin{equation}\label{d11}
\left(\sum^k_{i=1}p_i\right)\left(\sum^N_{i=k+1}p_i+1\right)
-\sum^N_{i=1}p^2_i.
\end{equation}
Non-negativity of the latter expression we will prove by induction
on the number of vertices $N$. Recall that the only assumption,
restricting the choice of the number $k\geq 1$, is that there are
no arrows $i\to j$ for $i\geq j+2$ and $i\leq k$.

Consider first the case $k=1$. Assume that $l\geq 1$ vertices are
connected by arrows with 1, that is,
$$
2\to 1,\dots,l+1\to 1,\quad \mbox{but}\quad l+2\nrightarrow 1.
$$
In this case $p_1=p_2+\dots+p_{l+1}$ and the subgraph of the graph
$\Gamma$ with the vertices $\{2,\dots,l+1\}$ either consists of
one vertex or is a chain. The expression (\ref{d11}) transforms to
the formula
$$
\left(\sum^{l+1}_{i=2}p_i\right)
\left(\sum^N_{i=l+2}p_i+1\right)-\sum^N_{i=2}p^2_i,
$$
so that one can apply the induction hypothesis to the subgraph
with the vertices $\{2,\dots,N\}$. This completes the case $k=1$.

Now let $k\geq 2$. The following key fact is true.\vspace{0.1cm}

{\bf Lemma 4.9.} {\it The following inequality holds:}
\begin{equation}\label{d12}
p_i\leq\sum^N_{j=i+2}p_i+1.
\end{equation}
\vspace{0.1cm}

{\bf Proof:} this is Lemma 1.6 in [33].\vspace{0.1cm}

By the lemma that we have just proved, we get the inequality
$$
p_1=p_2=\dots=p_{k-1}\leq\sum^N_{i=k+1}p_i+1.
$$
For this reason, for $k\geq 2$ the expression (\ref{d11}) is
bounded from below by the number
$$
\left(\sum^k_{i=2}p_i\right)\left(\sum^N_{i=k+1}p_i+1\right)-
\sum^N_{i=2}p^2_i.
$$
Now, applying the induction hypothesis to the subgraph with the
vertices $\{2,\dots,N\}$ we complete the proof of Lemma 4.8 and
Proposition 4.2.\vspace{0.3cm}


\section{The technique of counting multiplicities}

In this section we give a stronger version of the technique of
counting multiplicities for the self-intersection of a movable
linear system. We obtain the result that, together with the
$8n^2$-inequality, forms the technical basis for the exclusion of
maximal singularities, the centre of which is of codimension
$\geq 4$. The notations in this section are independent of other
parts of this paper.\vspace{0.3cm}

{\bf 5.1. Set up of the problem.} Let $o\in X$ be a germ of a
smooth three-dimensional variety, $\varphi\colon\widetilde{X}\to
X$ a birational morphism, $E\subset\widetilde{X}$ an irreducible
exceptional divisor over the point $o$, that is, $\varphi(E)=o$.
Consider the resolution [1,40] of the discrete valuation $\nu_E$,
that is, the sequence of blow ups
$$
\varphi_{i,i-1}\,\,\colon X_i\to X_{i-1},
$$
$i=1,\dots,N$, where $X_0=X,\varphi_{i,i-1}$ blows up an
irreducible subvariety $B_{i-1}\subset X_{i-1}$ (a point or a
curve),
$$
E_i=\varphi^{-1}_{i,i-1}(B_{i-1})\subset X_i
$$
is the exceptional divisor, where $B_i$ is uniquely defined by
the conditions $B_0=o$ and for $i=1,\dots,N-1$
$$
B_i=\mathop{\rm centre}(E,X_i),
$$
and, finally, the geometric discrete valuations
$$
\nu_E\quad\mbox{and} \quad\nu_{E_N}
$$
of the field of rational functions of the variety $X$ coincide.
Geometrically this means that the birational map
$$
\varphi^{-1}_{N,0}\circ\varphi\colon\widetilde{X} \dashrightarrow
X_N
$$
is biregular at the generic point of the divisor $E$ and maps $E$
onto $E_N$ (see the details in [1,40]). Here
$$
\varphi_{N,0}=\varphi_{1,0}\circ\dots\circ\varphi_{N,N-1} \colon
X_N\to X_0,
$$
and more generally set for $i>j$
$$
\varphi_{i,j}=\varphi_{j+1,j}\circ\dots\circ\varphi_{i,i-1}\colon
X_i\to X_j.
$$
The strict transform of an irreducible subvariety (by linearity,
also of an effective cycle) $Y\subset X_j$ on $X_i$ we denote, as
usual, by adding the upper index $i$: we write $Y^i$.

Assume that for $i=1,\dots,L\leq N$ the centres $B_{i-1}$ of blow
ups are points, for $i\geq L+1$ they are curves. Let $\Gamma$ be
the graph of the constructed resolution, that is, an oriented
graph with the vertices $1,\dots,N$, where an oriented edge
(arrow) joins $i$ and $j$ for $i>j$ (notation: $i\to j$), if and
only if
$$
B_{i-1}\subset E^{i-1}_j.
$$
In particular, by construction always $i+1\to i$.

Let us describe the obvious combinatorial properties of the graph
$\Gamma$.\vspace{0.1cm}

{\bf Lemma 5.1.} {\it Let $i<j<k$ be three distinct vertices. If
$k\to i$, then $j\to i$.}\vspace{0.1cm}

{\bf Proof.} By definition, $k\to i$ means that $B_{k-1}\subset
E^{k-1}_i$. By construction of the resolution of singularities,
we get
$$
\varphi_{k-1,j-1}(B_{k-1})=B_{j-1}
$$
(the centres of blow ups with higher numbers cover the centres of
previous blow ups) and, besides, obviously
$$
\varphi_{k-1,j-1}(E^{k-1}_i)=E^{j-1}_i.
$$
From here the lemma follows directly. Q.E.D.\vspace{0.1cm}

{\bf Definition 5.1.} We say that the vertex $i$ of the graph
$\Gamma$ is of {\it class} $e\geq 1$ (notation:
$\varepsilon(i)=e$), if precisely $e$ arrows come out of it, that
is,
$$
\sharp\{j|i\to j\}=e.
$$
We say, furthermore, that the graph $\Gamma$ is of class $e\geq
1$, if for each vertex $i$ we have $\varepsilon(i)\leq
e$.\vspace{0.1cm}

For instance, a graph of class 1 is a chain:
$$
1\longleftarrow 2\longleftarrow\dots\longleftarrow N
$$
(no other arrows but $i+1\to i$). The graph of a sequence of blow
ups of points on a non-singular surface is of class
2.\vspace{0.1cm}

{\bf Lemma 5.2.} {\it The graph of the resolution of the
valuation $\nu_E$ is of class 3. If for some vertex $i$ we have
$\varepsilon(i)=3$, then $i\leq L$, that is, $B_{i-1}$ is a
point.}\vspace{0.1cm}

{\bf Proof.} By definition,
$$
B_{i-1}\subset E(i)=\bigcup_{\{j|i\to j\}}E^{i-1}_j,
$$
and moreover, at the general point $B_{i-1}$ is a smooth variety
and $E(i)$ is a normal crossings divisor, each component of which
contains $B_{i-1}$. Q.E.D. for the lemma.\vspace{0.1cm}

{\bf Remark 5.1.} Word for word the same arguments show that in
the case of arbitrary dimension $\mathop{\rm dim}X$ the graph of
the resolution of any valuation is of class at most $\mathop{\rm
dim}X$ and if
$$
\sharp\{j|i\to j\}=a,
$$
then $\mathop{\rm codim}B_{i-1}\geq a$.

Let $\Sigma$ be a germ of a movable (that is, free from fixed
components) linear system on $X$, $\Sigma^i$ its strict transform
on $X_i$,
$$
\nu_i=\mathop{\rm mult}\nolimits_{B_{i-1}}\Sigma^{i-1},
$$
so that for a general divisor $D\in\Sigma$ we have
$$
D^i=D-\sum^i_{j=1}\nu_jE_j,
$$
where we write $D$ instead of $\varphi^*_{i,0}D$ and similarly
for the exceptional divisors $E_i$.

Consider a pair of general divisors $D_1,D_2\in\Sigma$ and
construct the {\it self-intersections} of linear systems
$\Sigma^i$, which are (not uniquely determined) effective 1-cycles
$$
Z_i=(D^i_1\circ D^i_2)
$$
on $X_i$. These cycles admit the natural decomposition
$$
\begin{array}{ccl}
Z_1 & = & Z^1_0+Z_{1,1},\\
Z_2 & = & Z^2_1+Z_{2,2}=Z^2_0+Z_{1,2}+Z_{2,2},\\
& \dots &\\
Z_i & = & Z^i_{i-1}+Z_{i,i}=Z^i_0+Z_{1,i}+\dots+Z_{i,i},\\
\end{array}
$$
where $Z_{a,i}=(Z_{a,i-1})^i=\dots=Z^i_{a,a}$, $i=1,\dots,L$.
\vspace{0.1cm}

{\bf Definition 5.2.} A function $a\colon\{1,\dots, L\}\to{\mathbb
Z}_+$ is said to be {\it compatible with the graph structure}
$\Gamma$, if the inequalities
$$
a(i)\geq\sum_{j\to i}a(j)
$$
hold. \vspace{0.1cm}

By construction, $B_L\subset E_L\cong{\mathbb P}^2$ is a plane
curve. Set $\beta_L=\mathop{\rm deg}B_L$ to be its degree in
${\mathbb P}^2$ and for an arbitrary $i\geq L+1$
\begin{equation}\label{e0}
\beta_i=\beta_L\mathop{\rm deg}[\varphi_{i,L}|_{B_i}\colon B_i\to
B_L].
\end{equation}
The main computational tool of the theory of birational rigidity
is the folowing local fact.\vspace{0.1cm}

{\bf Proposition 5.1.} {\it For any function $a(\cdot)$,
compatible with the graph structure, the inequality
\begin{equation}\label{e1}
\sum^L_{i=1}a(i)m_i\geq\sum^L_{i=1}a(i)\nu^2_i+
a(L)\sum^{N}_{i=L+1}\beta_i\nu^2_i
\end{equation}
{\it holds, where} $m_i=\mathop{\rm mult}_{B_{i-1}}Z^{i-1}_0$,
$i=1,\dots L$.}\vspace{0.1cm}

{\bf Proof} is given in [1,40].\vspace{0.1cm}

The aim of this section is to prove a stronger estimate that
includes (\ref{e1}) as a particular case.\vspace{0.1cm}

{\bf Definition 5.3.} A vertex $i\in\{4,\dots,L\}$ of the graph
$\Gamma$ is said to be {\it complex}, if precisely three arrows
come out of this vertex, that is, for three distinct vertices
$i_1<i_2<i_3$ we have
$$
i\to i_1,\,\,i\to i_2,\,\,i\to i_3.
$$
\vspace{0.1cm}

Note that in the notations of Definition 5.3, by Lemma 5.1 we
always have
$$
i_2\to i_1,\quad i_3\to i_1,\quad\mbox{and}\quad i_3\to i_2.
$$
If in the graph $\Gamma$ there are no complex vertices, then it
is of class $\leq 2$.\vspace{0.1cm}

{\bf Definition 5.4.} A {\it simplification} of the graph $\Gamma$
is the oriented graph $\Gamma^*$ of class $\leq 2$ with the set
of vertices
$$
1,\dots,L,
$$
the arrows in which join the vertices in accordance with the
following rule: if
$$
\sharp\{j|i\to j\}\leq 2,
$$
then $i\to j$ in $\Gamma^*$ if and only if $i\to j$ in $\Gamma$;
however, if $i$ is a complex vertex, then in the notations of
Definition 5.3 in $\Gamma^*$ there are two arrows coming out
ofthe vertex $i$,
$$
i\to i_2\quad\mbox{and}\quad i\to i_3,
$$
that is, the arrow $i\to i_1$ is deleted.\vspace{0.1cm}

Thus the graph $\Gamma^*$ is obtained from $\Gamma$ by means of
deleting some arrows (and the vertices that correspond to the
blow ups of curves).\vspace{0.1cm}

{\bf Proposition 5.2.} {\it For any function
$a\colon\{1,\dots,L\}\to{\mathbb Z}_+$, compatible with the
structure of the graph $\Gamma^*$, the inequality (\ref{e1})
holds.}\vspace{0.3cm}


{\bf 5.2. Proof of the improved inequality.} If the function
$a(\cdot)$ is compatible with the structure of $\Gamma$, then, the
more so, it is compatible with $\Gamma^*$, so that Proposition 5.2
implies Proposition 5.1. Following the general scheme of the proof
of the inequality (\ref{e1}) in [1,40], set for $i=1,\dots,L$
$$
d_i=\mathop{\rm deg}Z_{i,i},
$$
where $Z_{i,i}\subset E_i\cong{\mathbb P}^2$ is a plane curve and
for $i<j\leq L$
$$
m_{i,j}=\mathop{\rm mult}\nolimits_{B_{j-1}}Z_{i,j-1}.
$$
We obtain the following system of equalities
\begin{equation}\label{e2}
\begin{array}{ccl}
\nu^2_1+d_1 & = & m_1,\\
\nu^2_2+d_2 & = & m_2+ m_{1,2},\\
& \dots &\\
\nu^2_i+d_i & = &m_i+m_{1,i}+\dots+m_{i-1,i},\\
& \dots &\\
\nu^2_L+d_L & = &m_L+m_{1,L}+\dots+m_{L-1,L}.\\
\end{array}
\end{equation}
Besides, we get the obvious inequality
$$
d_L\geq\sum^{N}_{i=L+1}\beta_i\nu^2_i,
$$
where the numbers $\beta_i$ are defined by the formula (\ref{e0}).
Let us multiply the $i$-th equality in (\ref{e2}) by $a(i)$ and
put together the resulting equalities. In the left hand side we
get
$$
\sum^L_{i=1}a(i)\nu^2_i+\sum^L_{i=1}a(i)d_i.
$$
In the right hand side we get
$$
\sum^L_{i=1}a(i)m_i+\sum^{L-1}_{i=1}
\left(\sum^L_{j=i+1}a(j)m_{i,j}\right).
$$
Thus Proposition 5.2 follows immediately from the following
claim.\vspace{0.1cm}

{\bf Lemma 5.3.} {\it For any $i=1,\dots,L-1$ the inequality}
\begin{equation}\label{e3}
a(i)d_i\geq\sum^L_{j=i+1}a(j)m_{i,j}
\end{equation}
{\it holds.}\vspace{0.1cm}

{\bf Proof.} Up to this moment our argument repeated word for word
the corresponding arguments in [1,40]. If the function $a(\cdot)$
is compatible with the structure of the graph $\Gamma$ (and not
$\Gamma^*$), then, taking into account that the inequality
$m_{i,j}>0$ is possible only if $j\to i$, and that always
$$
d_i\geq m_{i,j},
$$
we obtain (\ref{e3}) from the inequality
$a(i)\geq\sum\limits_{j\to i}a(j)$ (this proves Proposition 5.1).
However, this is not sufficient for the proof of Proposition 5.2,
since the fucntion $a(\cdot)$ is compatible with the structure of
the graph $\Gamma^*$ only, and the latter has, generally speaking,
less arrows than $\Gamma$ does.

To prove (\ref{e3}), recall, first of all, that the integer-valued
functions
$$
d_i=\mathop{\rm deg}Z_{i,i}\quad\mbox{и}\quad m_{i,j}=\mathop{\rm
mult}\nolimits_{B_{j-1}}Z_{i,j-1}
$$
are linear functions of effective 1-cycles $Z_{i,i}$ on the
exceptional plane $E_i\cong{\mathbb P}^2$. Since the inequality
(\ref{e3}) is also linear, the claim of Lemma 5.3 follows from a
simpler fact.\vspace{0.1cm}

{\bf Lemma 5.4.}  {\it For any irreducible curve $C\subset E_i$
the following inequality holds:}
\begin{equation}\label{e4}
a(i)\mathop{\rm deg}C\geq\sum^L_{j=i+1}a(j)\mathop{\rm
mult}\nolimits_{B_{j-1}}C^{j-1}.
\end{equation}
\vspace{0.1cm}

{\bf Proof.} Set
$$
d=\mathop{\rm deg}C,\,\,\mu_j=\mathop{\rm
mult}\nolimits_{B_{j-1}}C^{j-1},\,\,j=i+1,\dots,L.
$$
As we pointed out above, if $\mu_j>0$, then $j\to i$, so that it
is necessary to prove the inequality
$$
a(i)d\geq\sum_{j\to i}a(j)\mu_j.
$$
This inequality is a claim about singularities of plane curves.
Let us consider two cases:\vspace{0.1cm}

1) $C\subset E_i$ is a line in $E_i\cong{\mathbb P}^2$,

2) $C$ is a curve of degree $d\geq 2$.\vspace{0.1cm}

\noindent In the case 1) define the integer $k\geq 1$ by the
condition
\begin{equation}\label{e5}
\{j|B_{j-1}\in C^{j-1}\}=\{i+1,\dots,i+k\}.
\end{equation}
In order to distinguish between the arrows in the graphs $\Gamma$
and $\Gamma^*$, we write $a\stackrel{*}{\to}b$, if the vertices
$a,b$ are joined by an arrow in $\Gamma^*$, leaving the usual
arrow for $\Gamma$.

The following fact is of key importance.\vspace{0.1cm}

{\bf Lemma 5.5.} {\it For each $e$, $1\leq e\leq k$, we have:}
$$
i+e\stackrel{*}{\to} i.
$$
\vspace{0.1cm}

{\bf Proof.} Assume that $(i+e)$ is a complex vertex of the graph
$\Gamma$ (otherwise $i+e\stackrel{*}{\to}i$ by definition). Recall
that the simplification procedure removes, of the three arrows,
coming out of $i+e$, the one that goes to the {\it lowest} vertex.
Therefore, we may assume that $e\geq 3$. However, the points
$$
B_i,B_{i+1},B_{i+2},\dots,B_{i+k-1}
$$
lie on the strict transform of the smooth curve $C$, and therefore
the subgraph
$$
i+1\longleftarrow i+2\longleftarrow
i+3\longleftarrow\dots\longleftarrow 1+k
$$
is a chain, that is, between the vertices $i+1,\dots,i+k$ in the
original graph $\Gamma$ there are no arrows, except for the
consecutive. In any case there are two arrows coming out of the
vertex $i+e$:
$$
i+e\to i+e-1\quad\mbox{and}\quad i+e\to i.
$$
What has been said implies that if a third arrow comes out the
vertex $i+e$, that is, $i+e\to j$, then inevitably
$$
j\leq i-1.
$$
It is this arrow that the simplification procedure deletes.
Therefore, the arrow $i+e\to i$ will not be deleted. Q.E.D. for
the lemma.\vspace{0.1cm}

Let is come back to the case 1). The inequality (\ref{e4}) takes
the form of the estimate
\begin{equation}\label{e6}
a(i)\geq\sum^{i+k}_{j=i+1}a(j).
\end{equation}
By the lemma we have just proved, (\ref{e6}) is true, because the
function $a(\cdot)$ is compatible with the structure of the graph
$\Gamma^*$. Q.E.D. for Lemma 5.4 in the case 1).\vspace{0.1cm}

Let us consider the case 2). Again let us define $k\geq 1$ by the
condition (\ref{e5}). If $k=1$, then there is nothing to prove,
since
$$
i+1\stackrel{*}{\to} i
$$
and $\mu_j\leq d$ for any $j$. For this reason we assume that
$k\geq 2$.\vspace{0.1cm}

{\bf Lemma 5.6.} {\it The inequality $d\geq\mu_1+\mu_2$ holds.}
\vspace{0.1cm}

{\bf Proof.} Let $L\subset E_i$ be the line, passing through the
point $B_i$ in the direction of the infinitely near point
$B_{i+1}\in E^{i+1}_i$. By assumption, $C\neq L$. Then for the
intersection number on the surface $E^{i+2}_i$ we get
$$
0\leq(C^{i+2}\cdot L^{i+2})=d-\mu_1-\mu_2,
$$
which is what we need. Q.E.D. for the lemma.\vspace{0.1cm}

Let $\Gamma_C$ be the subgraph of the graph $\Gamma$ with the
vertices $\{i+1,\dots,i+k\}$ (and the same arrows), $\Gamma^*_C$
the subgraph of the graph  $\Gamma^*$ with the same set of
vertices. By definition of the simplification procedure, the arrow
$$
i+a\to i+b,\quad a>b\leq 1,
$$
can not be deleted, because in $\Gamma$ there is an arrow that
goes to a lower vertex:
\begin{equation}\label{e7}
i+a\to i.
\end{equation}
This implies, that $\Gamma_C=\Gamma^*_C$. Furthermore, the arrow
(\ref{e7}) is deleted by the simplification procedure if and only
if $i+a$ is a complex vertex and the arrow (\ref{e7}) is the
lowest. In the language of the graph $\Gamma_C$ this means that
the vertex $i+a$ is of class 2 (as a vertex of that graph). The
following claim is obvious.\vspace{0.1cm}

{\bf Lemma 5.7.} {\it The integer-valued function
$$
i+a\longmapsto\mu_{i+a}\in{\mathbb Z}_+
$$
is compatible with the structure of the graph
$\Gamma_C$.}\vspace{0.1cm}

{\bf Proof:} this is a well known property of multiplicities of a
curve at infinitely near points (on a non-singular surface).
\vspace{0.1cm}

By what has been said, the proof of the inequality (\ref{e4}) in
the case 2) is reduced to the following combinatorial fact. Let
$\Delta$ be a subgraph of class 2 with the set of vertices
$\{1,\dots,k\}$ and $\varepsilon_{\Delta}(\cdot)\in\{0,1,2\}$ the
function of the class of a vertex $(\varepsilon_{\Delta}(1)=0)$.
Let $\mu(\cdot)$ and $a(\cdot)$ be ${\mathbb Z}_+$-valued
functions, compatible with the structure of the graph
$\Delta$.\vspace{0.1cm}

{\bf Lemma 5.8.} {\it The following inequality holds:}
\begin{equation}\label{e8}
(\mu(1)+\mu(2))\sum_{\varepsilon_{\Delta}(j)\leq
1}a(j)\geq\sum^k_{j=1}\mu(j)a(j).
\end{equation}
\vspace{0.1cm}

{\bf Proof} is given by induction on the number of vertices $k\geq
2$. If $k=2$, then $\varepsilon_{\Delta}(1)=0$,
$\varepsilon_{\Delta}(2)=1$ and the inequality (\ref{e8}) is of
the form
$$
(\mu(1)+\mu(2))(a(1)+a(2))\geq\mu(1)a(1)+\mu(2)a(2),
$$
so that there is nothing to prove.

Assume that $k\geq 3$ and the vertices 3 and 1 are not joined by
an arrow: $3\nrightarrow 1$, that is, $\varepsilon_{\Delta}(3)=1$.
In that case we apply the induction hypothesis to the graph
$\Delta_1$ with the vertices $\{2,\dots,k\}$ and the same arrows
as in $\Delta$: for that graph the inequality (\ref{e8}) takes the
form
$$
(\mu(2)+\mu(3))\sum_{\varepsilon_{\Delta}(j)=
1}a(j)\geq\sum^k_{j=2}\mu(j)a(j),
$$
whence, taking into account the inequality
$\mu(1)\geq\mu(2)\geq\mu(3)$, we obtain the required inequality
(\ref{e8}) for $\Delta$.

Assume that $k\geq 3$ and there are arrows
$$
3\to 1,\dots,2+l\to 1,
$$
where $l\geq 1$. In that case let us write down the left hand side
of (\ref{e8}) as
$$
(\mu(1)+\mu(2))\left(a(1)+a(2)+\sum_{\begin{array}{l}j\geq l+3,\\
\varepsilon_{\Delta}(j)=1\\
\end{array}}a(j)\right),
$$
and the right hand side of (\ref{e8}) as
$$
\mu(1)a(1)+\mu(2)a(2)+\sum^{l+2}_{j=3}
\mu(j)a(j)+\sum^k_{j=l+3}\mu(j)a(j).
$$
Since the functions $\mu(\cdot)$ and $a(\cdot)$ are compatible
with the structure of the graph, we get the inequality
$$
\mu(2)a(1)\geq\mu(2)(a(2)+a(3)+\dots+a(l+2))
$$
and the symmetric inequality for $\mu(1)a(2)$. Applying the
induction hypothesis to the subgraph $\Delta_{l+3}$ with the
vertices $\{l+3,\dots,k\}$, we complete the proof of the lemma.

Q.E.D. for Proposition 5.2.\vspace{0.3cm}


{\bf 5.3. Counting multiplicities of the self-intersection for a
non log canonical singularity.} Let $o\in X$ be a non-singular
three-dimensional germ, $\Sigma$ a movable linear system, such
that
$$
\mathop{\rm mult}\nolimits_o\Sigma\leq 2n,
$$
but the point $o$ is an isolated centre of non log canonical
singularities of the pair
\begin{equation}\label{e9}
\left(X,\frac{1}{n}\Sigma\right),
\end{equation}
where $n>0$ is some number. Let $Z=(D_1\circ D_2)$ be the
self-intersection of the system $\Sigma$ (an effective 1-cycle on
$X$). Let
$$
\varphi_{i,i-1}\colon X_i\to X_{i-1},
$$
$i=1,\dots,K$, be the sequence of blow ups of the centres of some
non log canonical singularity of the pair (\ref{e9}), $X_0=X$,
$\varphi_{i,i-1}$ blows up an irreducible subvariety
$B_{i-1}\subset X_{i-1}$, a point or a curve, $B_0=o$,
$$
E_i=\varphi^{-1}_{i,i-1}(B_{i-1})\subset X_i
$$
is the exceptional divisor, finally, $E_K\subset X_K$ realizes a
non log canonical singularity of the pair (\ref{e9}), that is,
the log Noether-Fano inequality
\begin{equation}\label{e10}
\sum^K_{i=1}p_{Ki}\nu_i>n\left(\sum^K_{i=1}p_{Ki}\delta_i+1\right)
\end{equation}
holds, where $\delta_i=\mathop{\rm codim}B_{i-1}-1$, $p_{K_i}$ is
the number of paths from $E_K$ to $E_i$ and, as usual,
$$
\nu_i=\mathop{\rm mult}\nolimits_{B_{i-1}}\Sigma^{i-1}\leq 2n,
$$
$\Sigma^i$ is the strict transform of the movable system $\Sigma$
on $X_i$. Set
$$
\{1,\dots,L\}=\{j\,|\,\mathop{\rm dim}B_{j-1}=0\}.
$$
Since by assumption $\nu_i\leq 2n$, (\ref{e10}) implies that
$K\geq L+1$, that is, among the centres of blow ups there is at
least one curve.

Denote by the symbol $\Gamma$ the oriented graph of the sequence
of blow ups $\varphi_{i,i-1}$, $i=1,\dots,K$, and by the symbol
$\Gamma^*$ its simplification, a graph of class $\leq 2$. The
number of paths from the vertex $i$ to the vertex $j$, $i>j$, in
the graph $\Gamma^*$ we denote by the symbol $p^*_{ij}$. By
definition, $p^*_{ii}=1$. Set also
$$
m_i=\mathop{\rm mult}\nolimits_{B_{i-1}}Z^{i-1},
$$
$i=1,\dots,L$. The following fact improves the classical
$4n^2$-inequality for the case of a non {\it log} canonical
singularity of the pair (\ref{e9}).\vspace{0.1cm}

{\bf Proposition 5.3.} {\it The following estimate holds:}
$$
\sum^L_{i=1}p^*_{Li}m_i>4n^2\left(\sum^L_{i=1}p^*_{Li}+1\right).
$$
\vspace{0.1cm}

{\bf Proof.} Set
$$
N=\mathop{\rm
min}\left\{e\,\left|\,\sum^e_{i=1}p_{ei}(\nu_i-\delta_in)>0\right.\right\},
$$
that is, $E_N$ is the non canonical singularity of the pair
(\ref{e9}) with the minimal number. In particular, the pair
(\ref{e9}) is canonical at $E_1,\dots,E_{N-1}$. It follows easily
from the inequality $\nu_1\leq 2n$ (see \S 4), that the segment of
the graph $\Gamma$ with the vertices
$$
L-1, L, L+1, \dots, N
$$
is a chain. Define the number $a$, $0\leq a\leq 1$, by the
equality
$$
a=\frac{1}{n}\sum^{N-1}_{i=1}p_{Ni}(\delta_in-\nu_i)
$$
(the numbers of paths $p_{N_i}$ and $p_{N-1,i}$ coincide by what
was said above). If $a=1$, then it is easy to see that
$$
\nu_1=\dots=\nu_N=2n,
$$
whence the claim of the proposition follows directly. For that
reason, we assume that $a<1$.

First assume that $N=K$. In that case the technique of counting
multiplicities together with the inequality (\ref{e10}) gives:
$$
\sum^L_{i=1}p^*_{Li}m_i>\frac{(2\Sigma^*_0+\Sigma^*_1+1)^2}
{\Sigma^*_0+\Sigma^*_1}n^2=
4(\Sigma^*_0+1)n^2+\frac{(\Sigma^*_1-1)^2}{\Sigma^*_0+\Sigma^*_1}n^2,
$$
where
$$
\Sigma^*_{2-\alpha}=\sum\limits_{\delta_i=\alpha}p^*_{Ni},\quad
\alpha=1,2,
$$
which was required (note that the log Noether-Fano inequality can
get only stronger when we replace $p_{Ki}$ by $p^*_{Ki}$, since
those coefficients are changed only for $i=1,\dots,L-1$ and
$\nu_i\leq 2n$ for these values of $i$, this is a well known
fact).

Thus for $N=K$ Proposition 5.3 holds.

Assume that $K\geq N+1$.\vspace{0.1cm}

{\bf Lemma 5.9.} {\it The pair
$$
(X_{N-1},\frac{1}{n}\Sigma^{N-1}-aE_{N-1})
$$
is not log canonical at the curve $B_{N-1}\subset
E_{N-1}$.}\vspace{0.1cm}

{\bf Proof.} Taking into account the standard properties of the
numbers $p_{ij}$, rewrite the inequality (\ref{e10}) in the form
$$
\left(\sum_{\alpha\to
N-1}p_{K\alpha}\right)\left(\sum^{N-1}_{i=1}p_{Ni}(\nu_i-\delta_in)\right)
+\sum^K_{i=N}p_{Ki}(\nu_i-\delta_in)>n
$$
(obviously, for $i\leq N-1$ the equality $p_{Ni}=p_{N-1,i}$
holds). Now the claim of the lemma follows directly from the
definition of the number $a$ and the fact that $N\to N-1$. Q.E.D.
for the lemma.

Now setting
$$
\Sigma^*_0=\sum^L_{i=1}p^*_{Li},\quad\Sigma^*_1
=\sum^{N-1}_{i=L+1}p^*_{Ni},
$$
by Propositions 4.2 and 5.2 we get
$$
\sum^L_{i=1}p^*_{Li}m_i>\sum^{N-1}_{i=1}p^*_{Ni}\nu^2_i+4(1+a)n^2,
$$
whence, taking into account the Noether-Fano inequality, we get
\begin{equation}\label{e11}
\sum^L_{i=1}p^*_{Li}m_i>\left[\frac{(2\Sigma^*_0+\Sigma^*_1-a)^2}
{\Sigma^*_0+\Sigma^*_1}+4(1+a)\right]n^2=
4(\Sigma^*_0+1)n^2+\frac{(\Sigma^*_1+a)^2}{\Sigma^*_0+\Sigma^*_1}n^2.
\end{equation}

Q.E.D. for Proposition 5.3.\vspace{0.3cm}


\section{Exclusion of infinitely near\\ maximal singularities}

In this section we complete the proof of Proposition 0.3: under
the assumption that the system $\Sigma$ has no maximal subvariety
of the form $\sigma^{-1}(P)$, where $P\subset{\mathbb P}$ is a
linear subspace of codimension two, we prove that $\Sigma$ has no
maximal singularities with the centre $B$ of codimension $\geq
4$.\vspace{0.3cm}

{\bf 6.1. The centre of the singularity is not contained in the
ramification divisor.} By Propositions 2.1 and 3.1 we may assume
that the centre $B$ of the maximal singularity of the linear
system $\Sigma$ is of codimension $\geq 4$ (and the same is true
for any other maximal singularity). Let $P\subset{\mathbb P}$ be
a generic linear subspace of dimension $\mathop{\rm codim}B$ and
$o\in\sigma^{-1}(P)\cap B$ is some point. Set
$V_P=\sigma^{-1}(P)$.

Consider first the case when $\sigma(B)\not\subset W$. In this
case $p=\sigma(o)\not\in W$. The variety $V_P$ is smooth,
$$
\sigma_P=\sigma|_{V_P}\colon V_P\to P
$$
is the double cover, branched over $W_P=W\cap P$,
$$
\mathop{\rm Pic}V_P={\mathbb Z}H_P,
$$
where $H_P=H|_{V_P}$. Let $\Sigma_P$ be the restriction of the
system $\Sigma$ onto $V_P$. This is a movable linear system and
the pair
$$
(V_P,\frac{1}{n}\Sigma_P)
$$
is not log canonical, where the point $o$ is an isolated centre
of a non log canonical singularity. Let
$$
\varphi\colon V^+_P\to V_P
$$
be the blow up of the point $o$, $E=\varphi^{-1}(o)$ the
exceptional divisor. Set
$$
Z_P=(D_1\circ D_2)
$$
to be the self-intersection of the linear system $\Sigma_P$ and
$Z^+_P$ its strict transform on $V^+_P$. By Proposition 4.1, for
some plane $\Pi\subset E$ of codimension two the inequality
\begin{equation}\label{f1}
\mathop{\rm mult}\nolimits_oZ_P+\mathop{\rm
mult}\nolimits_{\Pi}Z^+_P>8n^2
\end{equation}
holds. Now let
$$
\varphi_P\colon P^+\to P
$$
be the blow up of the point $p=\sigma(o)$ and
$E_P=\varphi^{-1}_P(p)$ the exceptional divisor, $E$ identifies
naturally with $E_P$. Let $\Lambda\subset P$ be the unique plane
of codimension two, containing the point $p$ and cutting out $\Pi$
on $E_P=E$:
$$
\Lambda^+\cap E_P=\Pi,
$$
where $\Lambda^+\subset P^+$ is the strict transform. The
subvariety $Q=\sigma^{-1}_P(\Lambda)\subset V_P$ is irreducible,
of codimension two (with respect to $V_P$), and moreover,
\begin{equation}\label{f2}
\mathop{\rm deg}Q=2,\quad\mathop{\rm mult}\nolimits_oQ=\mathop{\rm
mult}\nolimits_{\Pi}Q^+=1,
\end{equation}
where $Q^+\subset V^+_P$ is the strict transform. Since the cycle
$Z_P$ satisfies the inequality
$$
\mathop{\rm mult}\nolimits_oZ_P+\mathop{\rm
mult}\nolimits_{\Pi}Z^+_P>\mathop{\rm deg}Z_P=8n^2,
$$
writing down
$$
Z_P=aQ+Z^{\sharp}_P,
$$
where $a\in{\mathbb Z}_+$ and $Z^{\sharp}_P$ does not contain $Q$
as a component, we obtain
$$
\mathop{\rm mult}\nolimits_oZ^{\sharp}_P+\mathop{\rm
mult}\nolimits_{\Pi}(Z^{\sharp}_P)^+>\mathop{\rm deg}Z_P,
$$
$(Z^{\sharp}_P)^+$ is the strict transform. Finally, let $R$ be
the $\sigma$-preimage of a generic hyperplane in $P$, containing
the point $p$ and cutting out $\Pi$ on $E_P=E$, that is,
$$
\sigma(R)^+\cap E_P=\Pi.
$$
The divisor $R$ contains no component of the effective cycle
$Z^{\sharp}_P$, so that for the scheme-theoretic intersection
$$
Z^{\sharp}_R=(Z^{\sharp}_P\circ R)
$$
we obtain the inequality
$$
\mathop{\rm mult}\nolimits_oZ^{\sharp}_R\geq\mathop{\rm
mult}\nolimits_oZ^{\sharp}_P+\mathop{\rm
mult}\nolimits_{\Pi}(Z^{\sharp}_P)^+>\mathop{\rm deg}Z^{\sharp}_R,
$$
which is impossible. This contradiction proves Proposition 0.3 in
the case when $\sigma(B)\not\subset W$.\vspace{0.3cm}


{\bf 6.2. The centre of the singularity is contained in the
ramification divisor: the simple case.} Consider, finally, the
last case when $\sigma(B)\subset W$. Once again, we work on the
variety $V_P=\sigma^{-1}(P)$, the linear system $\Sigma_P$ is
movable and the point $o$ is an isolated centre of log maximal
singularity of this linear system. Set $p=\sigma(o)$. For the
blow ups
$$
\varphi\colon V^+_P\to V_P\quad\mbox{and}\quad\varphi\colon P^+\to
P
$$
of the points $o$ and $p$, respectively, with the exceptional
divisors
$$
E=\varphi^{-1}(o)\quad\mbox{and}\quad E_P=\varphi^{-1}_P(p),
$$
the double cover $\sigma_P\colon V_P\to P$ does not extend to an
isomorphism of the exceptional divisors $E$ and $E_P$. Let
$T_P=T_pW_P$ be the tangent hyperplane to the branch divisor
$W_P=W\cap P$ at the point $p$, ${\mathbb T}_P$ the corresponding
hyperplane in $E_P$. It is easy to see that there exist a
hyperplane ${\mathbb T}\subset E$ and a point $\xi\in
E\backslash{\mathbb T}$ such that $\sigma_P$ induces an
isomorphism of ${\mathbb T}$ and ${\mathbb T}_P$ and the rational
map
$$
\sigma_E\colon E\dashrightarrow E_P
$$
is the composition of the projection $\mathop{\rm pr}_{\xi}\colon
E\dashrightarrow{\mathbb T}$ from the point $\xi$ and the
isomorphism ${\mathbb T}\cong{\mathbb T}_P$. In particular,
$\sigma_E(E)={\mathbb T}_P$ (all these facts are easy to check in
suitable local coordinates $z_1,\dots,z_k$ at the point $p$ on
$P$, in which $V_P$ is given by the local equation $y^2=z_1$).

Let, as above, $Z_P$ and $Z^+_P$ be the self-intersection of the
linear system $\Sigma_P$ and its strict transform on $V^+_P$,
respectively. Let $\Pi\subset E$ be the plane of codimension two,
satisfying the inequality (\ref{f1}),
$$
\Pi_P=\sigma_E(\Pi)\subset{\mathbb T}_P\subset E_P
$$
the image of the plane $\Pi$. Obviously, $\Pi_P$ is a linear
subspace in ${\mathbb T}_P$ of codimension 1 or 2. In the latter
case for a generic hyperplane $R\ni p$, $R\subset P$, such that
$R^+\supset\Pi_P$, we get: none of the components of the
effective cycle $(\sigma_P)_*Z_P$ of codimension two is not
contained in $R$. By the inequality (\ref{f1}) we get
$$
\mathop{\rm mult}\nolimits_o(\sigma^{-1}(R)\circ
Z_P)\geq\mathop{\rm mult}\nolimits_oZ_P+\mathop{\rm
mult}\nolimits_{\Pi}Z^+_P>
$$
$$
>8n^2=\mathop{\rm deg}Z_P=\mathop{\rm
deg}(\sigma^{-1}(R)\circ Z_P),
$$
which is impossible. Therefore, $\Pi_P$ is a hyperplane in
${\mathbb T}_P$. Now we can not argue in the same way as in the
case when the point $o$ does not lie on the ramification divisor:
let $\Lambda\subset P$ be the unique plane of codimension two,
such that $\Lambda\ni p$ and $\Lambda^+\cap E_P=\Pi_P$. The
subvariety $Q=\sigma^{-1}_P(\Lambda)$ is irreducible, however,
$$
\mathop{\rm mult}\nolimits_oQ=2,\quad\mathop{\rm
mult}\nolimits_{\Pi}Q^+=1,
$$
so that the aruments, similar to the case when $\sigma(o)\not\in
W$, do not work. To exclude this case, we need the improved
technique of counting multiplicities (\S 5). We complete the
proof in all details for $M\geq 6$.\vspace{0.3cm}


{\bf 6.3. The centre of the singularity is contained in the
ramification divisor: the hard case.} The following claim is true.
\vspace{0.1cm}

{\bf Lemma 6.1.} {\it We have the inequality} $\nu=\mathop{\rm
mult}_o\Sigma\leq 2n$.\vspace{0.1cm}

{\bf Proof.} Consider the divisor $T=\sigma^{-1}(T_pW)$, where
$T_pW\subset{\mathbb P}$ is the tangent hyperplane. The system
$\Sigma$ is movable, so that the effective cycle $(D\circ T)$ is
well defined, where $D\in\Sigma$ is a general divisor. Now we get
$$
2\nu\leq\mathop{\rm mult}\nolimits_o(D\circ T)\leq\mathop{\rm
deg}(D\circ T)=\mathop{\rm deg}D=4n,
$$
as we claimed.

Let $\Delta\subset P$, $\Delta\ni p$, be a generic 3-plane, so
that $V_{\Delta}=\sigma^{-1}(\Delta)$ is a smooth variety,
$\sigma_{\Delta}\colon V_{\Delta}\to \Delta$ a double cover,
$\Sigma_{\Delta}$ the restriction of the system $\Sigma$ onto
$V_{\Delta}$. The pair
\begin{equation}\label{f2a}
(V_{\Delta},\frac{1}{n}\Sigma_{\Delta})
\end{equation}
is not log canonical at the point $o$, and, moreover, $o$ is an
isolated centre of non log canonical singularities of this pair.

Set $C=Q\cap V_{\Delta}=\sigma^{-1}(L)$, where $L=\Delta\cap
\Lambda$ is the line, passing through the point $p$ and tangent to
$W$ at that point (for the definition of the plane $\Lambda$, see
above in the end of Sec. 6.2). Set also
$$
y=\Pi\cap V^+_{\Delta},
$$
where $V^+_{\Delta}$ is the strict transform of $V_{\Delta}$ on
$V^+$, that is,
$$
\varphi_{\Delta}=\varphi|_{V^+_{\Delta}}\colon V^+_{\Delta}\to
V_{\Delta}
$$
is the blow up of the point $o$ with the exceptional plane
$E_{\Delta}$, $y=\Pi\cap E_{\Delta}$. There is a non log
canonical singularity of the pair (\ref{f2a}), the centre of
which on $V^+_{\Delta}$ is the point $y$.

Consider the self-intersection
$$
Z_{\Delta}=(D_1\circ D_2)=Z|_{V_{\Delta}}
$$
of the movable linear system $\Sigma_{\Delta}$ and write down
$$
Z_{\Delta}=bC+Z_1,
$$
where $b\in{\mathbb Z}_+$, the effective 1-cycle $Z_1$ does not
contain the curve $C$ as a component and for this reason
satisfies the inequality
\begin{equation}\label{f3}
\mathop{\rm mult}\nolimits_oZ_1+\mathop{\rm
mult}\nolimits_yZ^+_1\leq\mathop{\rm deg}Z_1=8n^2-2b.
\end{equation}

Assume now that at the point $o$ the curve $C$ has two distinct
branches:
$$
C^+\cap E_{\Delta}=\{y,y^*\},
$$
where $y,y^*\in E_{\Delta}$ are distinct points. This assumption
is justified for $M\geq 6$ by the conditions of general position,
since for any subspace ${\cal U}\subset{\mathbb P}$ of
codimension two, ${\cal U}\ni p$, the quadratic point
$o\in\sigma^{-1}({\cal U})$ is of rank at least two, see
Proposition 0.6, (iii), so that the same is true for the
quadratic point $o\in\sigma^{-1}(L)=C$, either, since $L={\cal
U}\cap\Delta$, where $\Delta$ is a generic 3-plane, containing
the point $p$. For $M=5$ we must consider also the case when
$o\in C$ is a simple cuspidal singularity, see below.

In the notations of Sec. 5.3 let
$$
\varphi_{i,i-1}\colon X_i\to X_{i-1},
$$
$i=1,\dots,N$, $X_0=V_{\Delta}$, be the resolution of the non log
canonical singularity, the centre of which on $X_1=V^+_{\Delta}$
is the point $B_1=y$. Set
$$
\{1,\dots,k\}=\{i\,|\,1\leq i\leq L,B_{i-1}\in C^{i-1}\}.
$$
By the assumption about the branches of the curve $C$, we get
$k\geq 2$, and, moreover, the subgraph with the vertices
$1,\dots,k$ is a chain.

Note that $b\geq 1$: otherwise the inequality
$$
\mathop{\rm mult}\nolimits_oZ_{\Delta}+\mathop{\rm
mult}\nolimits_yZ^+_{\Delta}\leq\mathop{\rm deg}Z_{\Delta}=8n^2
$$
holds and one can argue in word for word the same way as for
$\sigma(o)\not\in W$.

We have
$$
p^*_{L1}=\dots=p^*_{L,k-1}.
$$
Set, furthermore,
$$
\mu_i=\mathop{\rm mult}\nolimits_{B_{i-1}}Z^{i-1}_1
$$
for $i=1,\dots,L$. By the inequality (\ref{f3}), we get the
estimate
\begin{equation}\label{f4}
\begin{array}{c}
\sum\limits^L_{i=1}p^*_{Li}m_i=b\left(\sum\limits^k_{i=1}p^*_{Li}+p^*_{L1}\right)
+\sum\limits^L_{i=1}p^*_{Li}\mu_i\leq\\
\leq
b\left(\sum\limits^k_{i=1}p^*_{Li}+p^*_{L1}\right)+\frac12\mathop{\rm
deg}Z_1\sum\limits^L_{i=1}p^*_{Li}.
\end{array}
\end{equation}

Since $\Gamma^*$ is a graph of class $\leq 2$, we get the estimate
$$
p^*_{L1}=\dots=p^*_{L,k-1}\leq 1+\sum\limits^L_{i=k+1}p^*_{Li}
$$
(see Lemma 4.9). Therefore the right hand side of the inequality
(\ref{f4}) is bounded from above by the number
$$
b+(b+\frac12\mathop{\rm
deg}Z_1)\sum^L_{i=1}p^*_{Li}=b+4n^2\sum\limits^L_{i=1}p^*_{Li}.
$$
Now Proposition 5.3 implies the estimate
$$
b>4n^2,
$$
which is impossible. This contradiction completes the proof of
Proposition 0.3 for $M\geq 6$.\vspace{0.1cm}

For $M=5$ to complete the proof it remains to consider the case
when the curve $C$ has a simple cuspidal singularity at the point
$o$, so that $C^+$ is tangent to the exceptional divisor $E$ at
the point $y$ (the tangency is simple). If $L=2$, then the
previous arguments work. If $L\geq 3$ and $B_2\in C^2$, then the
vertices 3 and 1 are joined by an arrow, $3\to 1$, so that
$$
p^*_{L1}\geq p^*_{L2}+p^*_{L3}
$$
and the input of the first vertex of the graph into the sum
$\sum\limits^L_{i=1}p^*_{Li}m_i$ is not compensated by the number
$$
b\left(1+\sum\limits^L_{i=k+1}p^*_{Li}\right),
$$
as above. However, the previous estimates can be improved in the
following way. Set
$$
\mu=\mathop{\rm mult}\nolimits_C\Sigma\leq n
$$
(if $\mu>n$, then the system $\Sigma$ has a maximal subvariety of
codimension two, which is what we need). Let us restrict $\Sigma$
onto the $\sigma$-preimage $S$ of a generic 2-plane in ${\mathbb
P}$, containing the line $L$. Obviously, $\Sigma_S=\Sigma|_S$ is a
non-empty linear system of curves with a single fixed component
$C$ of multiplicity $\mu$. For a generic curve $G\in\Sigma_S$ we
get the inequality
$$
((G-\mu C)\cdot C)\geq 2\mathop{\rm mult}\nolimits_o(G-\mu C)
+\sum\limits^k_{i=2}\mathop{\rm mult}\nolimits_{B_{i-1}}(G-\mu C),
$$
where
$$
\{2,\dots,k\}=\{i\,|\,B_{i-1}\in C^{i-1}\},
$$
$k\geq 3$. There is a standard technique (used in [41, \S 8], also
in [12,15,32]), which makes it possible to derive from this
estimate and the Noether-Fano inequality, that the ``upper'' sum
$\Sigma^*_1$ is high compared with the ``lower'' one $\Sigma^*_0$,
whence a considerable improvement of the quadratic inequality
(\ref{e11}) is obtained. That improved inequality is already
sufficient to exclude the cuspidal case. We omit the details; they
will be given in another paper. This completes the proof of
Proposition 0.3 for $M=5$.\vspace{0.3cm}


\section{Double spaces of general position}

In this section we prove Propositions 0.4-0.6.\vspace{0.3cm}

{\bf 7.1. Lines on the variety $V$.} Let us prove Proposition 0.4.
The non-trivial part of that claim is that through every point
there are at most finitely many lines; the fact that any (not
necessarily generic) double space of index two is swept out by
lines, is almost obvious. It is easy to see that the image
$L=\sigma(C)$ of a line $C\subset V$ on ${\mathbb P}$ is a line in
the usual sense and
$$
\sigma|_C\colon C\to L\subset{\mathbb P}
$$
is an isomorphism. Thus there are two possible cases: either
$L\not\subset W$, so that $\sigma^{-1}(L)=C\cup C^*$ is a part of
smooth rational curves (permuted by the Galois involution of the
double cover $\sigma$), or $L\subset W$ is contained entirely in
the branch divisor, that is, $\sigma^{-1}(L)=C$. The converse is
also true: if a line $L\subset{\mathbb P}$ is such that the curve
$\sigma^{-1}(L)$ is reducible or $L\subset W$, then
$\sigma^{-1}(L)$ consists of two or one lines on $V$,
respectively. An easy dimension count shows that on a generic
hypersurface in ${\mathbb P}$ of degree $2(M-1)$ there are no
lines, so that the second option does not take place. Furthermore,
the double cover $\sigma^{-1}(L)\to L$ is reducible if and only if
the divisor $W\,|\,_L$ on $L={\mathbb P}^1$ is divisible by 2,
that is,
$$
\frac12(W\,|\,_L)\in\mathop{\rm Div}L
$$
is an integral divisor. Thus Proposition 0.4 follows immediately
from the following fact.\vspace{0.1cm}

{\bf Proposition 7.1.} {\it For a generic smooth hypersurface
$W\subset{\mathbb P}$ of degree $2(M-1)$ through every point
$x\in{\mathbb P}$ there are finitely many lines $L$ such that}
$W\,|\,_L\in 2\mathop{\rm Div}L$.\vspace{0.1cm}

{\bf Proof.} Let us denote by the symbol ${\cal P}_k({\mathbb
P}^l)$ the space of homogeneous polynomials of degree $k$ on the
projective space ${\mathbb P}^l$ (that is, $H^0({\mathbb
P}^l,{\cal O}_{{\mathbb P}^l}(k))$), considered as an affine
algebraic variety of dimension ${k+l}\choose l$. Let
$$
\begin{array}{cccc}
\mathop{\rm sq}\colon & {\cal P}_k({\mathbb P}^l) & \to & {\cal
P}_{2k}({\mathbb P}^l),\\
\mathop{\rm sq}\colon & f & \mapsto & f^2
\end{array}
$$
be the map of taking the square. Its image
$$
\mathop{\rm sq}({\cal P}_k({\mathbb P}^l))\subset{\cal
P}_{2k}({\mathbb P}^l)
$$
will be denoted by the symbol $[{\cal P}_k({\mathbb P}^l)]^2$.
Consider the space of pairs
$$
\Pi={\mathbb P}\times{\cal P}_{2(M-1)}({\mathbb P})
$$
and set $\Pi(x)=\{x\}\times{\cal P}_{2(M-1)}({\mathbb P})$ for an
arbitrary point $x\in{\mathbb P}$. Set
$$
Y(x)\subset\Pi(x)
$$
to be the closed algebraic subset of pairs $(x,F)$, $F\in{\cal
P}_{2(M-1)}({\mathbb P})$, defined by the condition \vspace{0.1cm}

$(+)$ the set of lines $L\subset{\mathbb P}$, $L\ni x$, for which
$F\,|\,_L\in[{\cal P}_{M-1}(L)]^2$, is of positive
dimension.\vspace{0.1cm}

It is easy to see that the closure
$$
\overline{\bigcup_{x\in{\mathbb P}}Y(x)}\subset\Pi
$$
is a closed algebraic subset of dimension $\leq M+\mathop{\rm
dim}Y(x)$. Therefore, Proposition 7.1, in its turn, is implied by
the following fact.\vspace{0.1cm}

{\bf Proposition 7.2.} {\it The codimension of the closed set
$Y(x)$ in $\Pi(x)\cong{\cal P}_{2(M-1)}({\mathbb P}$) is at least
$M+1$.}\vspace{0.1cm}

In fact, as we will see from the proof, a much stronger estimate
for the codimension of the set $Y(x)$ holds. In particular, the
claim of Proposition 7.1 remains true for double spaces of index
two with elementary singularities (quadratic points). However, we
do not need it here.\vspace{0.1cm}

{\bf Proof of Proposition 7.2.} Let $z_1,\dots,z_M$ be a system of
affine coordinates on ${\mathbb P}$ with the origin at the point
$x=(0,\dots,0)$. We write the polynomial $F\in{\cal
P}_{2(M-1)}({\mathbb P})$ in the form
$$
F=q_0+q_1(z_1,\dots,z_M)+\dots+q_{2(M-1)}(z_1,\dots,z_M),
$$
where $q_i(z_*)$ is a homogeneous polynomial of degree $i$. The
line $L\ni x$ corresponds to a set of homogeneous coordinates
$$
(a_1:\dots:a_M)\in{\mathbb P}^{M-1},
$$
$L=\{t(a_1,\dots,a_M)\,|\,t\in{\mathbb C}\}$. Obviously, $F|_L\in
[{\cal P}_{M-1}(L)]^2$ if and only if the polynomial
$$
q_0+tq_1(a_*)+\dots+t^{2(M-1)}q_{2(M-1)}(a_*)\in{\mathbb C}[t]
$$
is a full square in ${\mathbb C}[t]$.

For each $k=0,1,\dots,2(M-1)$ we define the set
$Y_k(x)\subset\Pi(x)$ by the condition\vspace{0.1cm}

$+_k$ there exists an irreducible closed subset $Z\subset{\mathbb
P}^{M-1}$ of positive dimension such that for a general line $L\in
Z$ we have $F_L\in[{\cal P}_{M-1}(L)]^2$, and moreover, $F|_L$ has
a zero of order $k$ at the point $x\in L$.\vspace{0.1cm}

(Recall that we identify the points ${\mathbb P}^{M-1}$ with the
lines in ${\mathbb P}$, passing through the point $x$.) Obviously,
for an odd $k\not\in 2{\mathbb Z}$ we have $Y_k=\emptyset$ and
$$
Y(x)=\bigcap^{M-1}_{i=0}Y_{2i}(x).
$$
Therefore, it is sufficient to prove the estimate of Proposition
7.2 for each of the (constructive) sets $Y_{2i}(x)$,
$i=0,\dots,M-1$. Consider first the set $Y_0(X)$ (corresponding to
the lines on $V$, passing through the point outside the branch
divisor). For $F\in Y_0(x)$ we have $q_0\neq 0$ and we may assume
that $q_0=1$.\vspace{0.1cm}

{\bf Lemma 7.1.} {\it For any $m\geq 1$ there exists a set of
quasi-homogeneous polynomials
$$
A_{m,i}(s_1,\dots,s_m)\in{\mathbb Q}[s_1,\dots,s_m]
$$
of degree $\mathop{\rm deg}A_{m,i}=i\in\{m+1,\dots,2m\}$, where
the weight of the variable $s_j$ is $\mathop{\rm wt}(s_j)=j$, such
that the polynomial
$$
1+b_1t+\dots+b_{2m}t^{2m}
$$
is a full square in ${\mathbb Q}[t]$ for
$b_1,\dots,b_{2m}\in{\mathbb C}$ if and only if the following
system of equalities is satisfied:}
$$
b_i=A_{m,i}(b_1,\dots,b_m),
$$
$i=m+1,\dots,2m$.\vspace{0.1cm}

{\bf Proof.} Consider the equality
$$
1+s_1t+\dots+s_{2m}t^{2m}=(1+r_1t+\dots+r_mt^m)^2.
$$
Equating the coefficients at the same powers of $t$, we find $r_i$
as polynomials in $s_1,\dots,s_i$ for $i\leq m$ (with coefficients
in ${\mathbb Z}[\frac12]$). The equality of coefficients at
$t^{m+1},\dots,t^{2m}$, gives the required system of equations.
Q.E.D.\vspace{0.1cm}

By the lemma, the restriction of $F$ onto the line $\{t(a_*)\}$ is
a full square if and only if the system of equations
\begin{equation}\label{g1}
q_i(a_*)=A_{M-1,i}(q_1(a_*),\dots,q_{M-1}(a_*)),
\end{equation}
$i=M,\dots,2(M-1)$, is satisfied. This is a system of $(M-1)$
polynomial homogeneous equations in $(a_1:\dots:a_M)$ of degrees
$M,\dots$, $2(M-1)$, respectively (which, in particular, implies
immediately that through every point $x\in V$ there are at least
two lines, and in the case of general position through $x\in V$
there are
$$
2\cdot(M\cdot(M+1)\cdot\dots\cdot 2(M-1))=2\frac{(2M-2)!}{(M-1)!}
$$
lines). As can be seen from (\ref{g1}), the coefficients of the
right hand side depend polynomially on the coefficients of the
polynomials $q_1,\dots,q_{M-1}$. Therefore, $Y_0(x)$ consists of
such polynomials $F\in{\cal P}_{2(M-1)}({\mathbb P})$, for which
the system (\ref{g1}) defines an algebraic set of positive
dimension.

Now the codimension $\mathop{\rm codim}Y_0(x)$ can be estimated by
the method of [4]. Assuming the polynomials $q_1,\dots,q_{M-1}$ to
be fixed, we obtain the equality
$$
\mathop{\rm codim}Y_0(x)=\mathop{\rm codim}Y^*_0(x),
$$
where the closed set  $Y^*_0\subset{\cal P}_M({\mathbb
P}^{M-1})\times\dots\times{\cal P}_{2(M-1)}({\mathbb P}^{M-1})$
(which consists of sets $q^*_M,\dots,q^*_{2(M-1)}$ of homogeneous
polynomials of the corresponding degrees) is defined by the
condition: the system of equations
$$
q^*_M=\dots=q^*_{2(M-1)}=0
$$
has a positive-dimensional set of solutions. Repeating the proof
of Lemma 3.3 in [1] (this argument is well known, it was applied
and published many times), define the subsets $Y^*_{0,j}\subset
Y^*_0$ for $j=M,\dots,2(M-1)$, fixing the first ``incorrect''
codimension:
$$
Y^*_{0,j}=\{(q^*_M,\dots,q^*_{2(M-1)})|\mathop{\rm
codim}\{q^*_M=\dots=q^*_j=0\}=j-M\},
$$
so that
$$
Y^*_0=\bigvee^{2(M-1)}_{j=M}Y^*_{0,j},
$$
where the symbol $\bigvee$ stands for a disjoint union (for
instance, $Y^*_{0,M}$ consists of the sets $(q^*_*)$ with
$q^*_M\equiv 0$). For $(q^*_M,\dots,q^*_{2(M-1)})\in Y^*_{0,j}$
there is an irreducible component
$$
B\subset\{q^*_M=\dots=q^*_{j-1}=0\}
$$
of codimension precisely $j-M$, on which $q^*_j$ vanishes
identically. Let
$$
\pi\colon B\to{\mathbb P}^{{\rm dim}B}\subset{\mathbb P} ^{M-1}
$$
be a generic linear projection onto a generic $\mathop{\rm
dim}B$-dimensional plane. Since the $\pi$-pull back of a non-zero
homogeneous polynomial on ${\mathbb P}^{{\rm dim}B}$ does not
vanish on $B$ identically, we get the estimate
$$
\mathop{\rm codim}Y^*_{0,j}\geq\mathop{\rm dim}{\cal P}_j({\mathbb
P}^{2M-1-j})={{2M-1}\choose j}
$$
(since $\mathop{\rm dim}B=2M-1-j$, $j\in \{M,\dots,2M-2\}$), so
that
$$
\mathop{\rm codim}Y^*_0\geq\mathop{\rm
min}\left\{\left.{{2M-1}\choose
j}\,\right|\,j=M,\dots,2M-2\right\}=2M-1.
$$
Thus $\mathbb{\rm codim}Y_0(x)\geq 2M-1\geq M+1$, as required.

Let us consider now the problem of estimating the codimension of
the set $Y_k(x)$, $k=2e\geq 2$. For $F\in Y_k(x)$ there exists a
set $Z_F\subset{\mathbb P}^{M-1}$ of positive dimension, on which
identically vanish the polynomials
$$
q_0,q_1(z_*),\dots,q_{k-1}(z_*)
$$
and for a point of general position $(a_1:\dots:a_M)\in Z_F$ the
polynomial
$$
t^kq_k(a_*)+\dots+t^{2(M-1)}q_{2(M-1)}(a_*)\in{\mathbb C}[t]
$$
is a full square, and moreover, $q_k(a_*)\neq 0$. Applying Lemma
7.1, we obtain the system of equalities
$$
\frac{q_{i}(a_*)}{q_k(a_*)}=A_{M-e-1,i}\left(
\frac{q_{k+1}(a_*)}{q_k(a_*)},\dots,\frac{q_{M+e-1}(a_*)}{q_k(a_*)}\right),
$$
$i=M-e,\dots,2(M-e-1)$, or, after multiplying by $q_k(a_*)^i$,
$$
q_{i}(a_*)q_k(a_*)^{i-1}=
A^+_{M-e-1,i}(q_k(a_*),q_{k+1}(a_*),\dots,q_{M+e-1}(a_*)),
$$
$i=M-e,\dots,2(M-e-1)$, where $A^+(\cdot)$ is the appropriately
modified polynomial. We obtain a system of $M+e-1$ homogeneous
equations on ${\mathbb P}^{M-1}$, and the set $Y_k(x)$ consists of
those polynomials $F$, for which that system has a set of
solutions of positive dimension, on which $q_k$ does not vanish
identically. Here $q_0$ is a non-zero constant by definition of
the set $Y_k(x)$ for $k=2e\geq 2$.

Now we argue as above: we assume the polynomials
$q_k,\dots,q_{M+e-1}$ to be fixed, so that the codimension of the
set $Y_k(x)$ is the codimension of the subset
$$
Y^*_k\subset{\mathbb C}\times{\cal P}_1({\mathbb P}^{M-1})
\times\dots\times{\cal P}_{k-1}({\mathbb P}^{M-1})\times{\cal
P}_{M+e}({\mathbb P}^{M-1})\times\dots\times{\cal
P}_{2(M-1)}({\mathbb P}^{M-1}),
$$
defined by the condition: the sequence
$$
(q_0,q_1,\dots,q_{k-1},q_{M+e},\dots,q_{2(M-1)})\in Y^*_k
$$
if and only if $q_0=0$ and the system of equations
$$
q_1=\dots=q_{k-1}=q_{M+e}=\dots=q_{2(M-1)}=0
$$
defines an algebraic set that has a component of positive
dimension, on which $q_k\not\equiv 0$. Now we may forget about the
latter condition.

Now we estimate the codimension of the set $Y^*_k$ by the method
of [4] (see [1]) in precisely the same way, as it was done above
for $k=0$: we fix the first ``incorrect'' codimension, when the
next polynomial $q_i$ in the sequence that was written out above
vanishes on an irreducible component, defined by the previous
equations (which is of ``correct'' codimension). We get the worst
estimate at the first step: the condition $q_1\equiv 0$ (together
with the condition $q_0=0$) gives the codimension
$$
\mathop{\rm codim}Y^*_2=M+1,
$$
and this estimate is optimal. Indeed, $q_1\equiv 0$ means that the
branch divisor is singular at the point $x$, and then through this
point there is a one-dimensional family of lines. In all other
cases the estimate for $\mathop{\rm codim}Y^*_k$ is considerably
stronger (we omit the elementary computations). Q.E.D. for
Proposition 0.4.\vspace{0.3cm}


{\bf 7.2. Isolated singular points.} The proof of Proposition 0.5
is elementary and we just point out its main steps. Assume that
for some subspace $P\subset{\mathbb P}$ of codimension two the
intersection $P\cap W$ has a whole curve $C$ of simgular points.
It is convenient to consider the pair $p\in C$, where $p$ is an
arbitrary point, so that in any case
$$
P\subset T_pW.
$$
There is a $(2M-3)$-dimensional family of pairs $(p,P\ni p)$,
satisfying this condition. It is sufficient to show that the
number of independent condition, which are imposed on the
(non-homogeneous) polynomial
$$
f|_P=f(z_1,\dots,z_{M-2})
$$
of degree $2M-2$ by the condition that the hypersurface
$\{f|_P=0\}=W\cap P$ contains a curve $C$ of singular points,
passing through $p=(0,\dots,0)$, is at least $2M-2$. There are
three possible cases:\vspace{0.1cm}

--- $C$ is a line,

--- $C$ is a plane curve, $C\subset \Lambda\subset P$, where $\Lambda$
is some 2-plane,

--- the linear span of the curve $C$ is a $k$-plane, where
$k\geq 3$.\vspace{0.1cm}

In the first case one can compute the number of independent
conditions precisely, this is an elementary exercise.

In the second case the plane curve $\{f|_{\Lambda}=0\}$ has an
irreducible component $C$ of degree $\geq 2$ and multiplicity
$\geq 2$, which gives an estimate from below for the number of
independent conditions (which is essentially stronger than we
need).

In the third case we choose on the curve $C$ $3(M-1)$ points in
general position (neither three lie on a line and neither four lie
in the same plane). It is easy to check (considering hypersurfaces
that are unions of singular quadrics), that being singular at
these points imposes on $f$ independent conditions, which
completes the proof of Proposition 0.5. (In fact, the codimension
of the set of hypersurfaces with a whole curve of singular points
is much higher, but we do not need that.) The details are left to
the reader.\vspace{0.3cm}


{\bf 7.3. The rank of quadratic singularities.} Let us prove
Proposition 0.6. We will show the claim (i). The claims (ii) and
(iii) are proved in a similar way. It is easy to check that the
planes $P\subset{\mathbb P}$ of codimension two that are tangent
to $W$ at at least one point (that is, $\mathop{\rm Sing}P\cap
W\neq\emptyset$) form a $(2M-3)$-dimensional family. So it is
sufficient to prove the following fact.

Let $P\subset{\mathbb P}$ be a fixed plane of codimension two,
$p\in P$ a fixed point,
$$
{\cal W}={\mathbb P}(H^0({\mathbb P},{\cal O}_{\mathbb P}(2M-2)))
$$
the space of hypersurfaces of degree $2M-2$. Let us define the
subset ${\cal W}_P\subset{\cal W}$ by the conditions:

\begin{itemize}
\item a hypersurface $W\in{\cal W}_P$ is non-singular at the point $p$,
\item the tangent hyperplane $T_pW$ contains $P$,
\item the point $p$ is an isolated singular point of the
intersection $P\cap W$.
\end{itemize}
For $W\in{\cal W}_P$ set $\sigma_W\colon V_W\to{\mathbb P}$ to be
the double cover, branched over $W$, $o=\sigma^{-1}_W(p)$ is a
singular point, $R=\sigma^{-1}_W(P)$,
$$
\varphi\colon V^+_W\to V_W
$$
the blow up of the subvariety $R$. On the variety $V^+_W$ there is
a unique singular point $o^+\in\varphi^{-1}(o)$. Let us define the
closed subset $Y\subset W_P$ by the condition that for $W\in Y$
the rank of the quadratic singularity $o^+\in V^+_W$ is at most 3.
Now Proposition 0.6, (i) follows immediately from the estimate
\begin{equation}\label{g2}
\mathop{\rm codim}(Y\subset W_P)\geq 2M-2
\end{equation}
for $M\geq 6$.

The proof of the inequality (\ref{g2}) is obtained by simple local
computations which we will just describe. Let $(z_1,\dots,z_M)$ be
affine coordinates at the point $p$, where the plane $P$ is
defined by the system of equations $z_1=z_2=0$, and the tangent
hyperplane to $W$ is $z_1=0$. The local equation of the double
cover $V_W$ at the point $o=\sigma^{-1}_W(p)$ is of the form
$$
y^2=z_1+q_2(z_1,\dots,z_M)+q_3(z_*)+\dots\,\,,
$$
and the local equation of the blow up $V^+_W$ at the point $o^+$
is of the form
$$
u^2=u_1u_2+q_2(0,u_2,\dots,u_M)+\dots\,\,.
$$
This implies that the condition that the rank of the quadratic
point $o^+$ is at most $M+1-k$, imposes
$$
\frac{k(k+1)}{2}
$$
independent conditions on the coefficients of the equation of the
hypersurface $W$. If $M+1-k\leq 3$, then we get at least
$$
\frac{(M-2)(M-1)}{2}
$$
independent conditions, which for $M\geq 6$ is strictly higher
than $2M-3$. Q.E.D. for Proposition 0.6, (i). The claims (ii) and
(iii) are shown in a similar way.\vspace{0.3cm}


\section*{References}

{\small

\noindent 1. Pukhlikov A.V., Birationally rigid varieties. I. Fano
varieties. Russian Math. Surveys. {\bf 62} (2007), no. 5, 857-942.
\vspace{0.3cm}

\noindent 2. Graber T., Harris J. and Starr J. Families of
rationally connected varieties. J. Amer. Math. Soc. {\bf 16}
(2002), no. 1, 57-67. \vspace{0.3cm}

\noindent 3. Iskovskikh V.A. and Manin Yu.I., Three-dimensional
quartics and counterexamples to the L\" uroth problem, Math. USSR
Sb. {\bf 86} (1971), no. 1, 140-166. \vspace{0.3cm}

\noindent 4. Pukhlikov A.V., Birational automorphisms of Fano
hypersurfaces, Invent. Math. {\bf 134} (1998), no. 2, 401-426.
\vspace{0.3cm}

\noindent 5. Pukhlikov A.V., Birationally rigid Fano complete
intersections, Crelle J. f\" ur die reine und angew. Math. {\bf
541} (2001), 55-79. \vspace{0.3cm}

\noindent 6. Pukhlikov A.V., Maximal singularities on the Fano
variety $V^3_6$. Moscow Univ. Math. Bull. {\bf 44} (1989), no. 2,
70-75.\vspace{0.3cm}

\noindent 7. Iskovskikh V.A. and Pukhlikov A.V., Birational
automorphisms of multi-dimensional algebraic varieties, J. Math.
Sci. {\bf 82} (1996), 3528-3613. \vspace{0.3cm}

\noindent 8. Cheltsov I.A. and Grinenko M.M., Birational rigidity
is not an open property. 2006. arXiv:math/0612159, 26
p.\vspace{0.3cm}

\noindent 9. Fano G., Nuove ricerche sulle varieta algebriche a
tre dimensioni a curve-sezioni canoniche, Comm. Rent. Ac. Sci.
{\bf 11} (1947), 635-720. \vspace{0.3cm}

\noindent 10. Iskovskikh V.A., Birational automorphisms of
three-dimensional algebraic varieties, J. Soviet Math. {\bf 13}
(1980), 815-868.\vspace{0.3cm}

\noindent 11. Khashin S.I., Birational automorphisms of the
Veronese double cone of dimension three. Moscow Univ. Math. Bull.
1984, no. 1, 13-16.\vspace{0.3cm}

\noindent 12. Grinenko M.M., Mori structures on Fano threefold of
index 2 and degree 1. Proc. Steklov Inst. Math. {\bf 246} (2004,
no. 3), 103-128.\vspace{0.3cm}

\noindent 13. Pukhlikov A.V., Birational automorphisms of
three-dimensional algebraic varieties with a pencil of del Pezzo
surfaces, Izvestiya: Mathematics {\bf 62} (1998), no. 1,
115-155.\vspace{0.3cm}

\noindent 14. Iskovskih V.A., Algebraic threefolds with special
regard to the problem of rationality. Proc. of Intern. Cong. of
Math. 1983. V. 1-2, 733-746.\vspace{0.3cm}

\noindent 15. Grinenko M.M., New Mori structures on a double space
of index 2. Russian Math. Surveys. {\bf 59} (2004), no. 3,
573-574.\vspace{0.3cm}

\noindent 16. Grinenko M.M., Fibrations into del Pezzo surfaces,
Russian Math. Surveys.  {\bf 61} (2006), no. 2, 255-300.
\vspace{0.3cm}

\noindent 17. Clemens H., Griffiths Ph. The intermediate Jacobian
of the cubic three-fold. Ann. of Math. {\bf 95} (1972), 73-100.
\vspace{0.3cm}

\noindent 18. Tyurin A.N., Five lectures on threefolds. Russian
Math. Surveys. {\bf 27} (1972), no. 5.\vspace{0.3cm}

\noindent 19. Tyurin A.N., Intermediate Jacobian of
three-dimensional varieties. Cont. Probl. Math. V. 12 (1979).
\vspace{0.3cm}

\noindent 20. Koll\' ar J., Nonrational hypersurfaces, J. Amer.
Math. Soc. {\bf 8} (1995), 241-249.\vspace{0.3cm}

\noindent 21. Tikhomirov A. S. The intermediate Jacobian of the
double covering of $P^{3}$ branched at a quartic. Math.
USSR-Izvestiya. {\bf 17} (1981), no. 3, 523-566.\vspace{0.3cm}

\noindent 22. Tikhomirov A. S. Singularities of the theta-divisor
of the intermediate Jacobian of a double cover of $P^{3}$ of index
two. Math. USSR-Izvestiya. {\bf 21} (1983), no. 2,
355-373.\vspace{0.3cm}

\noindent 23. Tikhomirov A. S. Letter to the editors of the
journal ``Izvestiya AN SSSR. Seriya Matematicheskaya'', Math.
USSR-Izvestiya.  {\bf 27} (1986), no. 1, 201.\vspace{0.3cm}

\noindent 24. Call F. and Lyubeznik G., A simple proof of
Grothendieck's theorem on the parafactoriality of local rings,
Contemp. Math. {\bf 159} (1994), 15-18.\vspace{0.3cm}

\noindent 25. Pukhlikov A.V., Birational automorphisms of a double
space and a double quadric, Math. USSR Izv. {\bf 32} (1989),
233-243. \vspace{0.3cm}

\noindent 26. Pukhlikov A.V., Birational geometry of Fano direct
products, Izvestiya: Mathematics, {\bf 69} (2005), no. 6,
1225-1255. \vspace{0.3cm}

\noindent 27. Pukhlikov A.V., Explicit examples of birationally
rigid Fano varieties. Mosc. Math. J. {\bf 7} (2007), no. 3,
543-560. \vspace{0.3cm}

\noindent 28. Cheltsov I., Park J. and Won J., Log canonical
thresholds of certain Fano hypersurfaces. 2007. arXiv:0706.075, 32
p.\vspace{0.3cm}

\noindent 29. Koll{\'a}r J., et al., Flips and Abundance for
Algebraic Threefolds, Asterisque {\bf 211}, 1993. \vspace{0.3cm}

\noindent 30. Koll\' ar J., Singularities of pairs, in: Algebraic
Geometry, Santa Cruz 1995, 221-286.\vspace{0.3cm}

\noindent 31. Corti A., Singularities of linear systems and 3-fold
birational geometry, in: ``Explicit Birational Geometry of
Threefolds'', London Mathematical Society Lecture Note Series {\bf
281} (2000), Cambridge University Press, 259-312. \vspace{0.3cm}

\noindent 32. Pukhlikov A.V., Birational automorphisms of a
three-dimensional quartic with an elementary singularity, Math.
USSR Sb. {\bf 63} (1989), 457-482. \vspace{0.3cm}

\noindent 33. Pukhlikov A.V., Birationally rigid iterated Fano
double covers. Izvestiya: Mathematics. {\bf 67} (2003), no. 3,
555-596.  \vspace{0.3cm}

\noindent 34. Cheltsov I.A., Local inequalities and the birational
superrigidity of Fano varieties. Izvestiya: Mathematics. {\bf 70}
(2006), no. 3, 605-639.\vspace{0.3cm}

\noindent 35. Cheltsov I., Non-rationality of a four-dimensional
smooth complete intersection of a quadric and a quartic, not
containing a plane, Sbornik: Mathematics, {\bf 194} (2003),
1679-1699.\vspace{0.3cm}

\noindent 36. Cheltsov I., Double cubics and double quartics,
Math. Z. {\bf 253} (2006), no. 1, 75-86. \vspace{0.3cm}

\noindent 37. Pukhlikov A.V., On the $8n^2$-inequality. 2008.
arXiv:0811.0183, 8 p.\vspace{0.3cm}

\noindent 38. Shokurov, V. V. Three-dimensional log flips.
Izvestiya: Mathematics. {\bf 40} (1993), no. 1, 95-202.
\vspace{0.3cm}

\noindent 39. Iskovskikh V.A., Birational rigidity of Fano
hypersurfaces in the framework of Mori theory. Russian Math.
Surveys. {\bf 56} (2001), no. 2, 207-291.\vspace{0.3cm}

\noindent 40. Pukhlikov A.V., Essentials of the method of maximal
singularities, in ``Explicit Birational Geometry of Threefolds'',
London Mathematical Society Lecture Note Series {\bf 281} (2000),
Cambridge University Press, 73-100. \vspace{0.3cm}

\noindent 41. Pukhlikov A.V., Birational isomorphisms of
four-dimensional quintics, Invent. Math. {\bf 87} (1987), 303-329.
\vspace{0.3cm}

}

\begin{flushleft}
e-mail: {\it pukh@liv.ac.uk}, {\it pukh@mi.ras.ru}
\end{flushleft}

\end{document}